\documentclass[thmsa,onecolumn,10pt]{article}
\usepackage{amssymb}
\usepackage{amsfonts}
\usepackage{graphicx}
\usepackage{amsmath}

\newtheorem{theorem}{Theorem}

\newtheorem{corollary}[theorem]{Corollary}

\newtheorem{lemma}[theorem]{Lemma}

\newtheorem{proposition}[theorem]{Proposition}

\newenvironment{proof}[1][Proof]{\textbf{#1.} }{\ \rule{0.5em}{0.5em}}

\begin{document}

\title{Stable convergence of multiple Wiener-It\^{o} integrals}
\author{Giovanni PECCATI\thanks{%
Laboratoire de Statistique Th\'{e}orique et Appliqu\'{e}e, Universit\'{e}
Paris VI, France. E-mail: \texttt{giovanni.peccati@gmail.com}}, and Murad
TAQQU\thanks{%
Boston University, Departement of Mathematics, 111 Cummington Road, Boston
(MA), USA.  E-mail: murad@math.bu.edu.}\ \thanks{%
This research was partially supported by the NSF Grant DNS-050547
at Boston University.}}
\date{April 25, 2006}
\maketitle

\begin{abstract}
We prove sufficient conditions, ensuring that a sequence of multiple
Wiener-It\^{o} integrals (with respect to a general Gaussian process)
converges stably to a mixture of normal distributions. Our key tool is an
asymptotic decomposition of contraction kernels, realized by means of
increasing families of projection operators. We also use an
infinite-dimensional Clark-Ocone formula, as well as a version of the
correspondence between \textquotedblleft abstract\textquotedblright\ and
\textquotedblleft concrete\textquotedblright\ filtered Wiener spaces, in a
spirit similar to \"{U}st\"{u}nel and Zakai (1997).


\end{abstract}

\section{Introduction}

Let $X$ be a centered Gaussian process and, for $d\geq 2$ and $n\geq 1$, let
$I_{d}^{X}\left( f_{n}\right) $ be a multiple Wiener-It\^{o} stochastic
integral, of order $d$, of some symmetric and square-integrable kernel $%
f_{n} $ with respect to $X$. The aim of this paper is to establish general
sufficient conditions on the kernels $f_{n}$, ensuring that the sequence $%
I_{d}^{X}\left( f_{n}\right) $ converges stably to a mixture of Gaussian
probability laws. The reader is referred e.g. to \cite[Chapter 4]{JacSh},
\cite{PeTaq} and Section 2.3 below, for an exhaustive characterization of
stable convergence. Here, we shall recall that such a convergence is
stronger than the convergence in law, and can be used in particular to
explain several non-central limit results for functionals of independently
scattered random measures; see for instance \cite{PeTaq}. Our starting point
is the following Central Limit Theorem (CLT).

\begin{theorem}[{see \protect\cite[Theorem 1 and Proposition 3]{NuPe}}]
\label{T : intro}If the variance of $I_{d}^{X}\left( f_{n}\right) $\
converges to 1 ($n\rightarrow +\infty $) the following three conditions are
equivalent: (i)\ $I_{d}^{X}\left( f_{n}\right) $\ converges in law to a
standard Gaussian random variable $N\left( 0,1\right) $, (ii)\ $E\left[
I_{d}^{X}\left( f_{n}\right) ^{4}\right] \rightarrow 3$, (iii)\ for every $%
r=1,...,d-1$, the contraction kernel $f_{n}\otimes _{d-r}f_{n}$\ converges
to $0$.
\end{theorem}

Although the implication (ii) $\Rightarrow $ (i) is rather striking, several
recent applications of Theorem \ref{T : intro} (see \cite{PT}, \cite{HuNu},
\cite{DPY} or \cite{Nua Corc W}) have shown that condition (iii) is easier
to verify than (ii), since in general there is no manageable formula for the
fourth moment of a non-trivial multiple Wiener-It\^{o} integral. Also, the
implication (iii) $\Rightarrow $ (i) (which can be regarded as a
simplification of the \textit{method of diagrams}---see e.g. \cite{Sur})
suggests that the asymptotic study of the contraction kernels associated to
the sequence $f_{n}$ may lead to more general convergence results. In
particular, in this paper we address the following problem. Let $Y\geq 0$ be
a non-constant random variable having the (finite) chaotic representation $%
Y=1+I_{2}^{X}\left( g_{2}\right) +\cdot \cdot \cdot +I_{2\left( d-1\right)
}^{X}\left( g_{2\left( d-1\right) }\right) $, let $N\sim N\left( 0,1\right) $
be independent of $Y$, and suppose that the sequence $I_{d}^{X}\left(
f_{n}\right) $ satisfies adequate normalization conditions; is it possible
to associate to each $f_{n}$ and each $r=1,...,d-1$, two \textit{generalized
contraction kernels}, say $f_{n}\otimes _{d-r}^{\ast }f_{n}$ and $%
f_{n}\otimes _{d-r}^{\ast \ast }f_{n}$, in such a way that the two relations%
\begin{equation}
f_{n}\otimes _{d-r}^{\ast }f_{n}\underset{n\rightarrow +\infty }{\rightarrow
}g_{2r}\text{ \ \ and \ \ }f_{n}\otimes _{d-r}^{\ast \ast }f_{n}\underset{%
n\rightarrow +\infty }{\rightarrow }0,\text{ \ }\forall r=1,...,d-1\text{,}
\label{CVintro}
\end{equation}%
imply that $I_{d}^{X}\left( f_{n}\right) $ converges stably to $\sqrt{Y}%
\times N$? This kind of non-central phenomena (convergence towards
non-trivial mixtures of Gaussian laws) appears regularly, for instance in
the analysis of the power variations of fractional processes (see e.g. \cite%
{Nua Corc W}), or in the study of non linear functionals of stationary
Gaussian sequences (see e.g. \cite[Theorems 7-9]{GirSur}). Although there
exists a panoply of results characterizing the stable convergence in a
semi-martingale setting (see \cite{LiSh}, \cite{Feigin} or \cite[Ch. 4]%
{JacSh}), none of them can be directly applied to the case of a Gaussian
process for which there is no explicit (semi)martingale structure (this is
true, in particular, for fractional processes). In this paper, we aim at
providing results in this direction for multiple integrals with respect to
general Gaussian processes, by using some ancillary devices borrowed from
continuous-time martingale calculus (in a spirit similar to \cite{NuPe}), as
well as a part of the theory of filtrations on general Wiener spaces, as
developed e.g. in \cite{Wu} and \cite{UZ} (see also \cite{PeTaq} for some
related results in a non-Gaussian framework).

\bigskip

Now let $\mathfrak{H}$ be a separable Hilbert space, and suppose that the
process $X=X\left( \mathfrak{H}\right) =\left\{ X\left( h\right) :h\in
\mathfrak{H}\right\} $ is a centered Gaussian measure (also called an
isonormal Gaussian process) over $\mathfrak{H}$ (see e.g. \cite[Ch. 1]%
{Nualart}, or Section 2.2 below). Then, $f_{n}\ $is a symmetric element of $%
\mathfrak{H}^{\otimes d}$ (i.e., the $d$th tensor product of $\mathfrak{H}$)
for every $n$, and $f_{n}\otimes _{d-r}f_{n}\in \mathfrak{H}^{\otimes 2r}$, $%
\forall r=1,...,d-1$. In what follows (see Theorem \ref{T : Main} and
formulae (\ref{backtoIntro1}) and (\ref{backtoIntro2}) below) we construct
the two kernels $f_{n}\otimes _{d-r}^{\ast }f_{n}$ and $f_{n}\otimes
_{d-r}^{\ast \ast }f_{n}$ appearing in (\ref{CVintro}), by using \textit{%
resolutions of the identity}. These objects are defined as continuous and
non-decreasing families of orthogonal projections $\pi =\left\{ \pi
_{t}:t\in \left[ 0,1\right] \right\} $ over $\mathfrak{H}$, indexed by $%
\left[ 0,1\right] $ and such that $\pi _{0}=0$ and $\pi _{1}=$ \textsc{Id.}.
Each resolution $\pi $ induces a time structure on the Gaussian field $%
X\left( \mathfrak{H}\right) $, and generates the canonical filtration $%
\mathcal{F}_{t}^{\pi }=\sigma \left\{ X\left( \pi _{t}h\right) :h\in
\mathfrak{H}\right\} $, $t\in \left[ 0,1\right] $ (note that $\mathcal{F}%
_{1}^{\pi }=\sigma \left( X\right) $ for every $\pi $). In particular, the
infinite dimensional process $t\mapsto \left\{ X\left( \pi _{t}h\right)
:h\in \mathfrak{H}\right\} \triangleq X\left( \pi _{t}\mathfrak{H}\right) $
can be seen as an infinite collection of possibly correlated Gaussian $%
\mathcal{F}_{t}^{\pi }$-martingales. As proved e.g. in \cite{Wu} and \cite%
{UZ} in the framework of abstract Wiener spaces, each $\sigma \left(
X\right) $-measurable and square integrable random variable (as Skorohod and
multiple Wiener-It\^{o} integrals with respect to $X$) is therefore the
terminal value of a $\mathcal{F}_{t}^{\pi }$-martingale, which is in turn a
\textquotedblleft generalized adapted stochastic integral\textquotedblright\
with respect to the infinite dimensional process $t\mapsto X\left( \pi _{t}%
\mathfrak{H}\right) $. Since every real-valued $\mathcal{F}_{t}^{\pi }$%
-martingale can be shown to be continuous, it follows that the stable
convergence of $\sigma \left( X\right) $-measurable random variables can be
studied by means of the theory of stable convergence for continuous local
martingales (see e.g. \cite[Ch. 4]{JacSh}). In particular, our starting
point in the construction of the two contraction operators appearing in (\ref%
{CVintro}) is a stable convergence result, proved in Proposition \ref{P :
stopping nest} below, involving the quadratic variation of continuous local
martingales, as well as a stochastic time-change result known as the
Dambis-Dubins-Schwarz Theorem (DDS Theorem) (see e.g. \cite[Ch. V]{ReYor}).
Observe that our Proposition \ref{P : stopping nest} is reminiscent of the
stable convergence results proved by Feigin in \cite{Feigin}. See \cite%
{vZ2000} for similar results involving the stable convergence of
multi-dimensional martingales, and \cite{PeTaq} for an alternative approach
based on a decoupling technique, known as the \textquotedblleft principle of
conditioning\textquotedblright .

\bigskip

We recall that the use of the DDS Theorem has already been crucial in the
proof of Theorem \ref{T : intro} and its generalizations, as stated in \cite%
{NuPe} and \cite{PT}. However, we shall stress that the proofs of the main
results of the present paper (in particular, Theorem \ref{T : Main} and
Theorem \ref{T : Main Double} below) are considerably more complicated.
Indeed, when no resolution of the identity is involved -- as it is the case
for Theorem \ref{T : intro} -- all infinite dimensional Gaussian spaces are
trivially isomorphic. It follows that every relevant element of the proof of
Theorem \ref{T : intro} is contained in the case of $X\left( \mathfrak{H}%
\right) $ being the Gaussian space generated by a standard one-dimensional
Brownian motion on $\left[ 0,1\right] $ (that is, $\mathfrak{H}=L^{2}\left( %
\left[ 0,1\right] \right) $), and the extension to general Gaussian measures
can be achieved by elementary considerations (see for instance \cite[Section
2.2]{NuPe}). However, in the present paper the filtrations $\mathcal{F}%
_{t}^{\pi }=\sigma \left\{ X\left( \pi _{t}h\right) :h\in \mathfrak{H}%
\right\} $ play a prominent role, and the complexity of these objects may
considerably vary, depending on the structure of the resolution $\pi $ (in
particular, depending on the \textit{rank} of $\pi $--see Section 2.1
below). We shall therefore use a notion of \textit{equivalence} between
pairs $\left( \mathfrak{H},\pi \right) $, where $\mathfrak{H}$ is a Hilbert
space and $\pi $ is a resolution, instead of the usual notion of isomorphism
between Hilbert spaces. The use of this equivalence relation implies that,
if the rank of $\pi $ equals $q$ ($q=1,...,+\infty $), then $\mathcal{F}%
_{t}^{\pi }$ has roughly the structure of the filtration generated by a $q$%
-dimensional Brownian motion. As a consequence, our first step will be the
proof of our main results in the framework of an infinite-dimensional
Brownian motion, and the extension to the general case will be realized by
means of rather delicate arguments involving the previously described
equivalence relation (see Lemma \ref{L : opIdentities} below). As will
become clear later on, our techniques can be regarded as a ramification of
the theory of \textit{concrete representations }for abstract Wiener spaces,
a concept introduced in \cite[Section \ 5]{UZ}. The reader is also referred
to \cite{PeTaq} for some related results in a non-Gaussian context.

\bigskip

The remainder of the paper is organized as follows. In Section 2.1, we
formally introduce the notion of resolution of the identity and discuss some
of its basic properties. In Section 2.2 some notions from stochastic
analysis and Skorohod integration are recalled. Sections 2.3-2.5 contain the
statements and the proofs of some useful stable convergence result for
Skorohod integrals. Section 3 is devoted to the proof of our main
convergence results. We also discuss some relations with the theory of
abstract Wiener spaces. An Appendix contains the proof of a technical lemma.

\section{Preliminary definitions and results}

Throughout the paper, the following conventions are in order: all random
objects are supposed to be defined on the same probability space $\left(
\Omega ,\mathcal{F},\mathbb{P}\right) $; all $\sigma $-fields are assumed to
be complete; the symbol $\overset{\mathbb{P}}{\rightarrow }$ stands for
convergence in probability; $\mathbb{R}$ is the set of real numbers.

\subsection{Hilbert spaces and resolutions of the identity}

Let $\mathfrak{H}$ be a real separable Hilbert space. The symbol $\left(
\cdot ,\cdot \cdot \right) _{\mathfrak{H}}$ indicates the inner product on $%
\mathfrak{H}$, and $\left\Vert \cdot \right\Vert _{\mathfrak{H}}=\left(
\cdot ,\cdot \right) _{\mathfrak{H}}^{1/2}$ as usual. The space $\mathfrak{H}
$ is always endowed with the Borel $\sigma $-field generated by the open
sets of the canonical distance associated to $\left\Vert \cdot \right\Vert _{%
\mathfrak{H}}$. As already done in \cite{PeTaq}, we first study the
convergence of Skorohod integrals by means of increasing families of
orthogonal projections, known as \textit{resolutions of the identity}.

\bigskip

\textbf{Definition I \ -- }A \textit{continuous} \textit{resolution of the
identity}, is a family $\pi =\left\{ \pi _{t}:t\in \left[ 0,1\right]
\right\} $ of orthogonal projections satisfying:

\begin{description}
\item[\textbf{(I-a)}] $\pi _{0}=0$, and $\pi _{1}=$ \textsc{Id.}$;$

\item[\textbf{(I-b)}] $\forall 0\leq s<t\leq 1$, $\pi _{s}\mathfrak{H}%
\subseteq \pi _{t}\mathfrak{H};$

\item[\textbf{(I-c)}] $\forall t_{0}\in \left[ 0,1\right] $, $\forall h\in
\mathfrak{H}$, $\lim_{t\rightarrow t_{0}}\left\Vert \left( \pi _{t}-\pi
_{t_{0}}\right) h\right\Vert _{\mathfrak{H}}=0.$
\end{description}

A subset $F$ of $\mathfrak{H}$ is said to be $\pi $-\textit{reproducing} if
the linear span of the set $\left\{ \pi _{t}f:f\in F\text{, }t\in \left[ 0,1%
\right] \right\} $ is dense in $\mathfrak{H}$. The \textit{rank} of $\pi $
is the smallest of the dimensions of all the subspaces generated by the $\pi
$-reproducing subsets of $\mathfrak{H}$. A $\pi $-reproducing subset $F$ of $%
\mathfrak{H}$ is \textit{fully orthogonal }if $\left( \pi _{t}f,g\right) _{%
\mathfrak{H}}=0$ for every $t\in \left[ 0,1\right] $ and every $f,g\in F$.
The collection of all $\pi $ verifying properties (\textbf{I-a})-(\textbf{I-c%
}) is noted $\mathcal{R}\left( \mathfrak{H}\right) $.

\bigskip

The reader is referred to \cite{Brod} or \cite{Yoshida} for further
properties and characterizations of the class $\mathcal{R}\left( \mathfrak{H}%
\right) $. In particular, we shall use the following consequence of \cite[%
Lemma 23.2]{Brod}, that can be proved by a standard Gram-Schmidt
orthogonalization.

\bigskip

\begin{lemma}
\label{L : Fully O}Let $\pi \in \mathcal{R}\left( \mathfrak{H}\right) $ and
let $F$ be a $\pi $-reproducing subset of $\mathfrak{H}$ such that $\dim
\left( \overline{F}\right) =rank\left( \pi \right) $, where $\overline{A}$
stands for the closure of the vector space generated by a given set $A$.
Then, there exists a $\pi $-reproducing \emph{and}\textit{\ }\emph{fully
orthogonal} subset $F^{\prime }$ of $\mathfrak{H}$, such that $\dim \left(
\overline{F^{\prime }}\right) =\dim \left( \overline{F}\right) $.
\end{lemma}

\bigskip

We will sometimes need to work with elements of $\mathcal{R}\left( \mathfrak{%
H}\right) $ that are not only continuous, but also absolutely continuous.

\bigskip

\textbf{Definition II }-- A resolution $\pi =\left\{ \pi _{t}:t\in \left[ 0,1%
\right] \right\} \in \mathcal{R}\left( \mathfrak{H}\right) $ is said to be
\textit{absolutely continuous} if, for every $f,g\in \mathfrak{H}$, the
function $t\mapsto \left( \pi _{t}g,f\right) _{\mathfrak{H}}$, $t\in \left[
0,1\right] $, is absolutely continuous with respect to the Lebesgue measure
on $\left[ 0,1\right] $. The class of absolutely continuous resolutions in $%
\mathcal{R}\left( \mathfrak{H}\right) $ is noted $\mathcal{R}_{AC}\left(
\mathfrak{H}\right) .$

\bigskip

The elements of $\mathcal{R}_{AC}\left( \mathfrak{H}\right) $ are used in
\cite[Section 5]{UZ} to prove a remarkable bijection between \textit{%
abstract }and \textit{concrete }filtered Wiener spaces. More details will be
given in Section 3, were we establish a similar result for isonormal
Gaussian processes as a step to prove stable convergence criteria for
multiple integrals. With the next result we point out that, up to a
\textquotedblleft change of time\textquotedblright , every $\pi \in \mathcal{%
R}\left( \mathfrak{H}\right) $ can be represented in terms of some element
of $\mathcal{R}_{AC}\left( \mathfrak{H}\right) $.

\bigskip

\begin{lemma}
\label{L : AC}For any $\pi =\left\{ \pi _{t}:t\in \left[ 0,1\right] \right\}
\in \mathcal{R}\left( \mathfrak{H}\right) $, there exists a non decreasing
function
\begin{equation*}
\psi =\left\{ \psi \left( t\right) :t\in \left[ 0,1\right] \right\}
\end{equation*}%
such that $\psi \left( 0\right) =0$ and the monotone family of projections%
\begin{equation*}
\widetilde{\pi }_{t}\triangleq \pi _{\psi \left( t\right) }\text{, \ \ }t\in %
\left[ 0,1\right] \text{,}
\end{equation*}%
is an element of $\mathcal{R}_{AC}\left( \mathfrak{H}\right) $.
\end{lemma}

\begin{proof}
Let $q=rank\left( \pi \right) $ ($q$ is possibly infinite) and let $F_{\pi
}=\left\{ f_{j}:1\leq j\leq q\right\} $ be a $\pi $-reproducing subset of $%
\mathfrak{H}$, normalized in such a way that $\sum_{j=1}^{q}\left\Vert
f_{j}\right\Vert _{\mathfrak{H}}^{2}=1$. Define moreover the increasing
function $\phi \left( t\right) =\sum_{j=1}^{q}\left\Vert \pi
_{t}f_{j}\right\Vert _{\mathfrak{H}}^{2}$, $t\in \left[ 0,1\right] $, and
set $\psi \left( t\right) =\inf \left\{ a:\phi \left( a\right) =t\right\} $.
Then, $\psi $ is non decreasing, $\psi \left( 0\right) =0$, and the family
of projections
\begin{equation*}
\widetilde{\pi }_{t}\triangleq \pi _{\psi \left( t\right) }\text{, \ \ }t\in %
\left[ 0,1\right] \text{,}
\end{equation*}%
is a resolution of the identity verifying $\sum_{j=1}^{q}\left\Vert
\widetilde{\pi }_{t}f_{j}\right\Vert _{\mathfrak{H}}^{2}=t$, for every $t\in %
\left[ 0,1\right] $. Since $F_{\pi }$ is also $\widetilde{\pi }$%
-reproducing, we deduce from \cite[Lemma 23.1]{Brod} that $\widetilde{\pi }$
is absolutely continuous.
\end{proof}

\subsection{Gaussian processes, Malliavin operators and representation
theorems}

Throughout the following, we write%
\begin{equation*}
X=X\left( \mathfrak{H}\right) =\left\{ X\left( f\right) :f\in \mathfrak{H}%
\right\}
\end{equation*}%
to indicate an \textit{isonormal Gaussian process}, or a \textit{Gaussian
measure}, over the Hilbert space $\mathfrak{H}$. This means that $X$ is a
centered Gaussian family, indexed by the elements of $\mathfrak{H}$ and
satisfying the isomorphic relation%
\begin{equation}
\mathbb{E}\left[ X\left( f\right) X\left( g\right) \right] =\left(
f,g\right) _{\mathfrak{H}}\text{, \ \ for every }f,g\in \mathfrak{H}
\label{basicISO}
\end{equation}%
(the notation $X\left( \mathfrak{H}\right) $ is adopted exclusively when the
role of $\mathfrak{H}$ is relevant to the discussion).

\bigskip

As in \cite{Wu} or \cite{UZ}, to every $\pi \in \mathcal{R}\left( \mathfrak{H%
}\right) $ we associate the collection of $\sigma $-fields%
\begin{equation}
\mathcal{F}_{t}^{\pi }\left( X\right) =\sigma \left\{ X\left( \pi
_{t}f\right) :f\in \mathfrak{H}\right\} \text{, \ \ }t\in \left[ 0,1\right]
\text{,}  \label{resfiltration}
\end{equation}%
and we observe that, for every $\pi \in \mathcal{R}\left( \mathfrak{H}%
\right) $, $t\mapsto \mathcal{F}_{t}^{\pi }\left( X\right) $ defines a
continuous filtration (see \cite[p. 14]{UZ}). Also, for every $f\in
\mathfrak{H}$, the process $t\mapsto X\left( \pi _{t}f\right) $, $t\in \left[
0,1\right] $, is a centered and continuous $\mathcal{F}_{t}^{\pi }\left(
X\right) $-martingale such that, for every $\eta >0$, the increment $X\left(
\left( \pi _{t+\eta }-\pi _{t}\right) f\right) =X\left( \pi _{t+\eta
}f\right) -X\left( \pi _{t}f\right) $ is independent of $\mathcal{F}%
_{t}^{\pi }\left( X\right) $ (see e.g. \cite[Corollary 2.1]{UZ}).

\bigskip

As in \cite{PeTaq}, we write $L^{2}\left( \mathbb{P},\mathfrak{H},X\right)
=L^{2}\left( \mathfrak{H},X\right) $ to indicate the set of $\sigma \left(
X\right) $-measurable and $\mathfrak{H}$-valued random variables $Y$ such
that $\mathbb{E}\left[ \left\Vert Y\right\Vert _{\mathfrak{H}}^{2}\right]
<+\infty $. The class $L^{2}\left( \mathfrak{H},X\right) $ is a Hilbert
space, with inner product given by $\left( Y,Z\right) _{L^{2}\left(
\mathfrak{H},X\right) }$ $=\mathbb{E}\left[ \left( Y,Z\right) _{\mathfrak{H}}%
\right] $. Following \cite{UZ}, we associate to every $\pi \in \mathcal{R}%
\left( \mathfrak{H}\right) $ the subspace $L_{\pi }^{2}\left( \mathfrak{H}%
,X\right) $ of $\pi $-\textit{adapted} elements of $L^{2}\left( \mathfrak{H}%
,X\right) $, that is: $Y\in L_{\pi }^{2}\left( \mathfrak{H},X\right) $ if,
and only if, $Y\in L^{2}\left( \mathfrak{H},X\right) $ and, for every $t\in %
\left[ 0,1\right] $ and every $h\in \mathfrak{H}$,%
\begin{equation}
\left( Y,\pi _{t}h\right) _{\mathfrak{H}}\in \mathcal{F}_{t}^{\pi }\left(
X\right) \text{.}  \label{adaptation}
\end{equation}

For any resolution $\pi \in \mathcal{R}\left( \mathfrak{H}\right) $, $L_{\pi
}^{2}\left( \mathfrak{H},X\right) $ is a closed subspace of $L^{2}\left(
\mathfrak{H},X\right) $. We may occasionally write $\left( u,z\right)
_{L_{\pi }^{2}\left( \mathfrak{H}\right) }$ instead of $\left( u,z\right)
_{L^{2}\left( \mathfrak{H}\right) }$, when both $u$ and $z$ are in $L_{\pi
}^{2}\left( \mathfrak{H},X\right) $. Now, for $\pi \in \mathcal{R}\left(
\mathfrak{H}\right) $, define $\mathcal{E}_{\pi }\left( \mathfrak{H}%
,X\right) $ to be the space of \textit{elementary} \textit{elements} of $%
L_{\pi }^{2}\left( \mathfrak{H},X\right) $, that is, $\mathcal{E}_{\pi
}\left( \mathfrak{H},X\right) $ is the collection of those elements of $%
L_{\pi }^{2}\left( \mathfrak{H},X\right) $ that are linear combinations of $%
\mathfrak{H}$-valued random variables of the type%
\begin{equation}
h=\Phi \left( t_{1}\right) \left( \pi _{t_{2}}-\pi _{t_{1}}\right) f\text{,}
\label{elementaryad}
\end{equation}%
where $t_{2}>t_{1}$, $f\in \mathfrak{H}$ and $\Phi \left( t_{1}\right) $ is
a $\mathcal{F}_{t_{1}}^{\pi }\left( X\right) $-measurable, real-valued and
square-integrable\ random variable. A proof of the following useful result
can be found in \cite[Lemma 3]{PeTaq} or \cite[Lemma 2.2]{UZ}.

\begin{lemma}
\label{L : adapted}For every $\pi \in \mathcal{R}\left( \mathfrak{H}\right) $%
, the span of the set $\mathcal{E}_{\pi }\left( \mathfrak{H},X\right) $ of
adapted elementary elements is dense in $L_{\pi }^{2}\left( \mathfrak{H}%
,X\right) $.
\end{lemma}

\bigskip

In what follows, we shall apply to the Gaussian measure $X$ some standard
notions and results from Malliavin calculus (the reader is again referred to
\cite{Nualart} and \cite{Nualart2} for any unexplained notation or
definition). For instance, $D=D_{X}$ and $\delta =\delta _{X}$ stand,
respectively, for the usual Malliavin derivative and Skorohod integral with
respect to the Gaussian measure $X$ (the dependence on $X$ will be dropped,
when there is no risk of confusion); for $k\geq 1$, $\mathbb{D}_{X}^{k,2}$
is the space of $k$ times differentiable functionals of $X$, endowed with
the norm $\left\Vert \cdot \right\Vert _{k,2}$ (see \cite[Chapter 1]{Nualart}
for a definition of this norm); $dom\left( \delta _{X}\right) $ is the
domain of the operator $\delta _{X}$. Note that $D_{X}$ is an operator from $%
\mathbb{D}_{X}^{k,2}$ to $L^{2}\left( \mathfrak{H},X\right) $, and also that
$dom\left( \delta _{X}\right) \subset L^{2}\left( \mathfrak{H},X\right) $.
For every $d\geq 1$, we define $\mathfrak{H}^{\otimes d}$ and $\mathfrak{H}%
^{\odot d}$ to be, respectively, the $d$th tensor product and the $d$th
\textit{symmetric} tensor product of $\mathfrak{H}$. For $d\geq 1$ we will
denote by $I_{d}^{X}$ the isometry between $\mathfrak{H}^{\odot d}$ equipped
with the norm $\sqrt{d!}\left\Vert \cdot \right\Vert _{\mathfrak{H}^{\otimes
d}}$ and the $d$th Wiener chaos of $X$. Given $g\in \mathfrak{H}^{\otimes d}$%
, we note $\left( g\right) _{s}$ the symmetrization of $g$, and%
\begin{equation*}
I_{d}^{X}\left( g\right) =I_{d}^{X}\left( \left( g\right) _{s}\right)
\end{equation*}

Plainly, for $f,g\in \mathfrak{H}^{\otimes d}$, $I_{d}^{X}\left( f+g\right)
=I_{d}^{X}\left( \left( f\right) _{s}+\left( g\right) _{s}\right) =$ $%
I_{d}^{X}\left( f\right) +I_{d}^{X}\left( g\right) $. Recall that, when $%
\mathfrak{H}=L^{2}\left( Z,\mathcal{Z},\nu \right) $, $\left( Z,\mathcal{Z}%
\right) $ is a measurable space, and $\nu $ is a $\sigma $-finite measure
with no atoms, then $\mathfrak{H}^{\odot d}=L_{s}^{2}\left( Z^{d},\mathcal{Z}%
^{\otimes d},\nu ^{\otimes d}\right) $, where $L_{s}^{2}\left( Z^{d},%
\mathcal{Z}^{\otimes d},\nu ^{\otimes d}\right) $ is the space of symmetric
and square integrable functions on $Z^{d}$. Moreover, for $f\in \mathfrak{H}%
^{\odot d}$, $I_{d}^{X}\left( f\right) $ coincides with the \textit{multiple
Wiener-It\^{o} integral} (of order $d$) of $f$ with respect to $X$, as
defined e.g. in \cite[Section 1.1.2]{Nualart}.

\bigskip

To establish the announced stable convergence results, we use the elements
of $\mathcal{R}\left( \mathfrak{H}\right) $ to represent random variables of
the type $\delta _{X}\left( u\right) $, $u\in dom\left( \delta _{X}\right) $%
, in terms of continuous-time martingales. In particular, we will use the
fact that (i) for any $\pi \in \mathcal{R}\left( \mathfrak{H}\right) $, $%
L_{\pi }^{2}\left( \mathfrak{H},X\right) \subseteq dom\left( \delta
_{X}\right) $, and (ii) for any $u\in L_{\pi }^{2}\left( \mathfrak{H}%
,X\right) $ the random variable $\delta _{X}\left( u\right) $ can be
regarded as the terminal value of a real-valued $\mathcal{F}_{t}^{\pi }$%
-martingale, where $\mathcal{F}_{t}^{\pi }$ is given by (\ref{resfiltration}%
). A proof of the following result can be found in \cite[Lemme 1]{Wu} and
\cite[Corollary 2.1]{UZ}

\bigskip

\begin{proposition}
\label{P : SkoAd}Let the assumptions of this section prevail. Then:

\begin{enumerate}
\item $L_{\pi }^{2}\left( \mathfrak{H},X\right) \subseteq dom\left( \delta
_{X}\right) $, and for every $h_{1},h_{2}\in L_{\pi }^{2}\left( \mathfrak{H}%
,X\right) $%
\begin{equation}
\mathbb{E}\left[ \delta _{X}\left( h_{1}\right) \delta _{X}\left(
h_{2}\right) \right] =\left( h_{1},h_{2}\right) _{L_{\pi }^{2}\left(
\mathfrak{H},X\right) }.  \label{isoskorohod}
\end{equation}

\item If $h\in \mathcal{E}_{\pi }\left( \mathfrak{H},X\right) $ has the form
$h=\sum_{i=1}^{n}h_{i}$, where $n\geq 1$, and $h_{i}\in \mathcal{E}_{\pi
}\left( \mathfrak{H},X\right) $ is s.t.%
\begin{equation*}
h_{i}=\Phi _{i}\times \left( \pi _{t_{2}^{\left( i\right) }}-\pi
_{t_{1}^{\left( i\right) }}\right) f_{i}\text{, \ \ }f_{i}\in \mathfrak{H}%
\text{, \ \ }i=1,...,n,
\end{equation*}%
with $t_{2}^{\left( i\right) }>t_{1}^{\left( i\right) }$ and $\Phi _{i}$
square integrable and $\mathcal{F}_{t_{1}^{\left( i\right) }}^{\pi }\left(
X\right) $-measurable, then
\begin{equation}
\delta _{X}\left( h\right) =\sum_{i=1}^{n}\Phi _{i}\times \left[ X\left( \pi
_{t_{2}^{\left( i\right) }}f_{i}\right) -X\left( \pi _{t_{1}^{\left(
i\right) }}f_{i}\right) \right] .  \label{simpleSkorohod}
\end{equation}

\item For every $u\in L_{\pi }^{2}\left( \mathfrak{H},X\right) $, the
process
\begin{equation*}
t\mapsto \delta _{X}\left( \pi _{t}u\right) \text{, \ \ }t\in \left[ 0,1%
\right] \text{,}
\end{equation*}%
is a continuous $\mathcal{F}_{t}^{\pi }\left( X\right) $-martingale
initialized at zero, with quadratic variation equal to
\begin{equation*}
\left\{ \left\Vert \pi _{t}u\right\Vert _{\mathfrak{H}}^{2}:t\in \left[ 0,1%
\right] \right\} .
\end{equation*}
\end{enumerate}
\end{proposition}

\bigskip

In the terminology of \cite{Wu}, relation (\ref{isoskorohod}) implies that $%
L_{\pi }^{2}\left( \mathfrak{H},X\right) $ is a closed subspace of the
\textit{isometric subset of }$dom\left( \delta _{X}\right) $, defined as the
collection of those $h\in dom\left( \delta _{X}\right) $ such that
\begin{equation}
\mathbb{E}\left( \delta _{X}\left( h\right) ^{2}\right) =\left\Vert
h\right\Vert _{L^{2}\left( \mathfrak{H},X\right) }^{2}.  \label{isosubset}
\end{equation}

Note that, in general, this isometric subset is not a vector space -- see
e.g. \cite[p. 170]{Wu}. The next result is partly a consequence of the
continuity of $\pi $. It is an abstract version of the \textit{Clark-Ocone
formula} (see \cite{Nualart}), and can be proved along the lines of \cite[Th%
\'{e}or\`{e}me 1, formula (2.4) and Th\'{e}or\`{e}me 3]{Wu}. Observe that,
in \cite{Wu}, such a result is proved in the context of abstract Wiener
spaces. However, such a proof uses exclusively isometric properties such as (%
\ref{isosubset}), and the role of the underlying probability space is
immaterial. It follows that the extension to general isonormal Gaussian
processes is standard: see e.g. \cite[Section 1.1]{Nualart2}. The reader is
also referred to \cite{MWZ} for a general Clark-Ocone formula concerning
Banach space valued Wiener functionals.

\bigskip

\begin{proposition}[Abstract Clark-Ocone formula]
\label{P : AbsCO}Under the above notation and assumptions (in particular, $%
\pi \in \mathcal{R}\left( \mathfrak{H}\right) $), for every $F\in \mathbb{D}%
_{X}^{1,2}$,
\begin{equation}
F=\mathbb{E}\left( F\right) +\delta \left( proj\left\{ D_{X}F\mid L_{\pi
}^{2}\left( \mathfrak{H},X\right) \right\} \right) \text{,}  \label{ACO}
\end{equation}%
where $D_{X}F$ is the Malliavin derivative of $F$, and $proj\left\{ \cdot
\mid L_{\pi }^{2}\left( \mathfrak{H},X\right) \right\} $ is the orthogonal
projection operator on $L_{\pi }^{2}\left( \mathfrak{H},X\right) $.
\end{proposition}

\bigskip

\textbf{Remarks -- }(a) The right-hand side of (\ref{ACO}) is well defined,
since $D_{X}F\in L^{2}\left( \mathfrak{H},X\right) $ by definition, and
therefore
\begin{equation*}
proj\left\{ D_{X}F\mid L_{\pi }^{2}\left( \mathfrak{H},X\right) \right\} \in
L_{\pi }^{2}\left( \mathfrak{H},X\right) \subseteq dom\left( \delta
_{X}\right) ,
\end{equation*}%
where the last inclusion is stated in Proposition \ref{P : SkoAd}.

(b) Since $\mathbb{D}_{X}^{1,2}$ is dense in $L^{2}\left( \mathbb{P}\right) $
and $\delta _{X}\left( L_{\pi }^{2}\left( \mathfrak{H},X\right) \right) $ is
an isometry (due to (\ref{isoskorohod})), formula (\ref{ACO}) yields that
every $F\in L^{2}\left( \mathbb{P},\sigma \left( X\right) \right) $ admits a
unique \textquotedblleft predictable\ representation\textquotedblright\ of
the form%
\begin{equation}
F=\mathbb{E}\left( F\right) +\delta _{X}\left( u\right) \text{, \ \ }u\in
L_{\pi }^{2}\left( \mathfrak{H},X\right) \text{;}  \label{Apredictable}
\end{equation}%
see also \cite[Remarque 2, p. 172]{Wu}.

\bigskip

In the next section, we present a general criterion (Theorem \ref{T : Gen CV}%
), ensuring the stable convergence of a sequence of Skorohod integrals
towards a mixture of Gaussian distributions. The result has been proved in
\cite{PeTaq}, by using a general convergence criteria for functionals of
independently scattered random measures. Here we present an alternative
proof (partly inspired by some arguments contained in \cite{vZ2000}), which
is based on a time-change technique for continuous-time martingales.

\subsection{Stable convergence of Skorohod integrals}

We first present a standard definition of the classes $\mathbf{M}$ and $%
\widehat{\mathbf{M}}$ of random probability measures and random fourier
transform.

\bigskip

\textbf{Definition III}\ -- Let $\mathcal{B}\left( \mathbf{%
\mathbb{R}
}\right) $ denote the Borel $\sigma $-field on $\mathbf{%
\mathbb{R}
}$.

\begin{description}
\item[\textbf{(III-a)}] A map $\mu \left( \cdot ,\cdot \right) $, from $%
\mathcal{B}\left( \mathbf{%
\mathbb{R}
}\right) \times \Omega $ to $\mathbf{%
\mathbb{R}
}$ is called a \textit{random probability }(on $\mathbf{%
\mathbb{R}
}$) if, for every $C\in \mathcal{B}\left( \mathbf{%
\mathbb{R}
}\right) $, $\mu \left( C,\cdot \right) $ is a random variable and, for $%
\mathbb{P}$-a.e. $\omega $, the map $C\mapsto \mu \left( C,\omega \right) $,
$C\in \mathcal{B}\left( \mathbf{%
\mathbb{R}
}\right) $, defines a probability measure on $\mathbf{%
\mathbb{R}
}$. The class of all random probabilities is noted $\mathbf{M}$, and, for $%
\mu \in \mathbf{M}$, we write $\mathbb{E}\mu \left( \cdot \right) $ to
indicate the (deterministic) probability measure
\begin{equation}
\mathbb{E}\mu \left( C\right) \triangleq \mathbb{E}\left[ \mu \left( C,\cdot
\right) \right] \text{, \ \ }C\in \mathcal{B}\left( \mathbf{%
\mathbb{R}
}\right) .  \label{mixt}
\end{equation}

\item[\textbf{(III-b)}] For a measurable map $\phi \left( \cdot ,\cdot
\right) $, from $\mathbf{%
\mathbb{R}
}\times \Omega $ to $\mathbb{C}$, we write $\phi \in \widehat{\mathbf{M}}$
if there exists $\mu \in \mathbf{M}$ such that
\begin{equation}
\phi \left( \lambda ,\omega \right) =\widehat{\mu }\left( \lambda \right)
\left( \omega \right) \text{, \ \ }\forall \lambda \in \mathbf{%
\mathbb{R}
}\text{, for }\mathbb{P}\text{-a.e. }\omega \text{,}  \label{randomfou}
\end{equation}%
where $\widehat{\mu }\left( \cdot \right) $ is defined as
\begin{equation}
\widehat{\mu }\left( \lambda \right) \left( \omega \right) =\left\{
\begin{array}{ll}
\int \exp \left( i\lambda x\right) \mu \left( dx,\omega \right) & \text{if }%
\mu \left( \cdot ,\omega \right) \text{ is a probability measure} \\
1 & \text{otherwise.}%
\end{array}%
\right. \text{, \ \ }\lambda \in \mathbf{%
\mathbb{R}
}.  \label{ftr}
\end{equation}
\end{description}

\bigskip

For every $\omega \in \Omega $, $\widehat{\mu }\left( \lambda \right) \left(
\omega \right) $ is of course a continuous function of $\lambda $, and the
probability $\mathbb{E}\mu \left( \cdot \right) =\int_{\Omega }\mu \left(
\cdot ,\omega \right) d\mathbb{P}\left( \omega \right) $ defined in (\ref%
{mixt}) is often called a \textit{mixture }of probability measures. The
notion of \textit{stable convergence}, which is the content of the next
definition,\textit{\ }extends the usual notion of convergence in law.

\bigskip

\textbf{Definition IV }(see e.g. \cite[Chapter 4]{JacSh}) -- Let $\mathcal{F}%
^{\ast }\subseteq \mathcal{F}$ be a $\sigma $-field, and let $\mu \in
\mathbf{M}$. A sequence of real valued r.v.'s $\left\{ Z_{n}:n\geq 1\right\}
$ is said to \textit{converge }$\mathcal{F}^{\ast }$-\textit{stably }to $%
\mathbb{E}\mu \left( \cdot \right) $, written $X_{n}\rightarrow _{\left( s,%
\mathcal{F}^{\ast }\right) }\mathbb{E}\mu \left( \cdot \right) $, if, for
every $\lambda ,\gamma \in \mathbf{%
\mathbb{R}
}$ and every $\mathcal{F}^{\ast }$-measurable r.v. $Z$,%
\begin{equation}
\lim_{n\rightarrow +\infty }\mathbb{E}\left[ \exp \left( i\gamma Z\right)
\times \exp \left( i\lambda X_{n}\right) \right] =\mathbb{E}\left[ \exp
\left( i\gamma Z\right) \times \widehat{\mu }\left( \lambda \right) \right]
\text{,}  \label{stable cond}
\end{equation}%
$\widehat{\mu }\in \widehat{\mathbf{M}}$ is given by (\ref{ftr}).

\bigskip

If $X_{n}$ converges $\mathcal{F}^{\ast }$-stably, then the conditional
distributions $\mathcal{L}\left( X_{n}\mid A\right) $ converge for any $A\in
\mathcal{F}^{\ast }$ such that $\mathbb{P}\left( A\right) >0$ (the reader is
referred e.g. to \cite[Proposition 5.33]{JacSh} for an exhaustive
characterization of stable convergence). By setting $Z=0$, we obtain that if
$X_{n}\rightarrow _{\left( s,\mathcal{F}^{\ast }\right) }\mathbb{E}\mu
\left( \cdot \right) $, then the law of the $X_{n}$'s converges weakly to $%
\mathbb{E}\mu \left( \cdot \right) $. Observe also that, if a sequence of
random variables $\left\{ U_{n}:n\geq 0\right\} $ is such that $\left(
U_{n}-Z_{n}\right) \rightarrow 0$ in $L^{1}\left( \mathbb{P}\right) $ and $%
X_{n}\rightarrow _{\left( s,\mathcal{F}^{\ast }\right) }\mathbb{E}\mu \left(
\cdot \right) $, then $U_{n}\rightarrow _{\left( s,\mathcal{F}^{\ast
}\right) }\mathbb{E}\mu \left( \cdot \right) $.

\bigskip

In what follows, $\mathfrak{H}_{n}$, $n\geq 1$, is a sequence of real
separable Hilbert spaces, whereas, for each $n\geq 1$, $X_{n}=X_{n}\left(
\mathfrak{H}_{n}\right) =\left\{ X_{n}\left( g\right) :g\in \mathfrak{H}%
_{n}\right\} $, is an isonormal Gaussian process over $\mathfrak{H}$. The
following theorem already appears in \cite{PeTaq}, where it is proved by
using a decoupling technique known as the \textquotedblleft principle of
conditioning\textquotedblright . In Section 2.4 we shall present an
alternative proof based exclusively on continuous-time martingale arguments.

\begin{theorem}
\label{T : Gen CV}Under the previous notation and assumptions, for $n\geq 1$%
, let $\pi _{n}=\left\{ \pi _{n,t}:t\in \left[ 0,1\right] \right\} \in
\mathcal{R}\left( \mathfrak{H}_{n}\right) $ and $u_{n}\in
L_{\pi_{n}}^{2}\left( \mathfrak{H}_{n},X_{n}\right) $. Suppose that there
exists a sequence $\left\{ t_{n}:n\geq 1\right\} $ $\subset \left[ 0,1\right]
$ and $\sigma $-fields $\left\{ \mathcal{U}_{n}:n\geq 1\right\} $, such that
\begin{equation}
\left\Vert \pi _{n,t_{n}}u_{n}\right\Vert _{\mathfrak{H}_{n}}^{2}\overset{%
\mathbb{P}}{\rightarrow }0  \label{f1}
\end{equation}%
and
\begin{equation}
\mathcal{U}_{n}\subseteq \mathcal{U}_{n+1}\cap \mathcal{F}_{t_{n}}^{\pi
_{n}}\left( X_{n}\right) .  \label{f2}
\end{equation}%
If%
\begin{equation}
\left\Vert u_{n}\right\Vert _{\mathfrak{H}}^{2}\overset{\mathbb{P}}{%
\rightarrow }Y\text{,}  \label{f3}
\end{equation}%
for some $Y\in L^{2}\left( \mathbb{P}\right) $ such that $Y\neq 0$, $Y\geq 0$
and $Y\in \mathcal{U}^{\ast }\triangleq \vee _{n}\mathcal{U}_{n}$, then, as $%
n\rightarrow +\infty $,
\begin{equation*}
\delta _{X_{n}}\left( u_{n}\right) \rightarrow _{\left( s,\mathcal{U}^{\ast
}\right) }\mathbb{E}\mu \left( \cdot \right) \text{,}
\end{equation*}%
where $\mu \in \mathbf{M}$ verifies $\widehat{\mu }\left( \lambda \right)
=\exp \left( -\frac{\lambda ^{2}}{2}Y\right) $.
\end{theorem}

\bigskip

\textbf{Remark -- }Condition (\ref{f2}), that already appears in the
statement of the main results of \cite{PeTaq}, can be seen as a weak version
of the \textit{nesting condition} used e.g. in \cite{Feigin} to establish
sufficient conditions for the stable convergence of semimartingales.

\bigskip

By using the Clark-Ocone formula stated in Proposition \ref{P : AbsCO}, we
deduce from\ Theorem \ref{T : Gen CV} a criterion for the stable convergence
of (Malliavin) differentiable functionals. It is the key to prove the main
results of the paper.

\begin{corollary}
\label{C : Diff Conv}Let $\mathfrak{H}_{n}$, $X_{n}\left( \mathfrak{H}%
_{n}\right) $, $\pi _{n}$, $t_{n}$ and $\mathcal{U}_{n}$, $n\geq 1$, satisfy
the assumptions of Theorem \ref{T : Gen CV}, and consider a sequence of
random variables $\left\{ F_{n}:n\geq 1\right\} $, such that $\mathbb{E}%
\left( F_{n}\right) =0$ and $F_{n}\in \mathbb{D}_{X_{n}}^{1,2}$ for every $n$%
. Then, a sufficient condition to have that
\begin{equation*}
F_{n}\rightarrow _{\left( s,\mathcal{U}^{\ast }\right) }\mathbb{E}\mu \left(
\cdot \right)
\end{equation*}%
where $\mathcal{U}^{\ast }\triangleq \vee _{n}\mathcal{U}_{n}$, $\widehat{%
\mu }\left( \lambda \right) =\exp \left( -\frac{\lambda ^{2}}{2}Y\right) $, $%
\forall \lambda \in \mathbb{R}$, and $Y\geq 0$ is such that $Y\in \mathcal{U}%
^{\ast }$, is
\begin{equation}
\left\Vert \pi _{n,t_{n}}proj\left\{ D_{X_{n}}F_{n}\mid L_{\pi
_{n}}^{2}\left( \mathfrak{H}_{n},X_{n}\right) \right\} \right\Vert _{%
\mathfrak{H}_{n}}^{2}\overset{\mathbb{P}}{\rightarrow }0\text{ \ \ and \ \ }%
\left\Vert proj\left\{ D_{X_{n}}F_{n}\mid L_{\pi _{n}}^{2}\left( \mathfrak{H}%
_{n},X_{n}\right) \right\} \right\Vert _{\mathfrak{H}_{n}}^{2}\overset{%
\mathbb{P}}{\underset{n\rightarrow +\infty }{\rightarrow }}Y\text{.}
\label{f4}
\end{equation}
\end{corollary}

\subsection{Martingale proof of Theorem \protect\ref{T : Gen CV}}

In this section, we provide a proof of Theorem \ref{T : Gen CV} involving
exclusively continuous martingale arguments. It is based on the following
general result.

\begin{proposition}
\label{P : stopping nest}Fix $T\leq +\infty $. For $n\geq 1$, let $\left\{
W_{t}^{n}:t\in \left[ 0,T\right) \right\} $ be a Brownian motion with
respect to a filtration $\mathcal{H}^{n}=\left\{ \mathcal{H}_{t}^{n}:t\in %
\left[ 0,T\right) \right\} $ of the space $\left( \Omega ,\mathcal{F},%
\mathbb{P}\right) $ (satisfying the usual conditions) and suppose that there
exists a sequence of random variables $\left\{ \tau _{n}:n\geq 1\right\} $
such that $\tau _{n}$ is a $\mathcal{H}^{n}$-stopping time with values in $%
\left[ 0,T\right) \ $and $\tau _{n}\overset{\mathbb{P}}{\rightarrow }0$ as $%
n\rightarrow +\infty $. Set moreover, for $n\geq 1$, $V_{t}^{n}=W_{\tau
_{n}+t}^{n}-W_{\tau _{n}}^{n}$, $t\in \left[ 0,T\right) $. Then,

\begin{enumerate}
\item $V^{n}-W^{n}\overset{\text{law}}{\Rightarrow }0;$

\item if there exists a sequence of $\sigma $-fields $\left\{ \mathcal{U}%
_{n}:n\geq 1\right\} $ such that
\begin{equation*}
\mathcal{U}_{n}\subseteq \mathcal{U}_{n+1}\cap \mathcal{H}_{\mathcal{\tau }%
_{n}}^{n}
\end{equation*}%
then for every random element $X$ defined on $\left( \Omega ,\mathcal{F},%
\mathbb{P}\right) $, with values in some Polish space $\left( S,\mathcal{S}%
\right) $ and measurable with respect to $\mathcal{U}^{\ast }\triangleq \vee
_{n}\mathcal{U}_{n}$,
\begin{equation}
\left( W^{n},X\right) \overset{\text{law}}{\Rightarrow }\left( W,X\right)
\text{ \ \ and \ \ }\left( V^{n},X\right) \overset{\text{law}}{\Rightarrow }%
\left( W,X\right)  \label{A1}
\end{equation}%
where $W$ is a standard Brownian motion independent of $X$.
\end{enumerate}
\end{proposition}

\begin{proof}
The proof is partly inspired by that of \cite[Theorem 3.1]{vZ2000}. Since $%
\tau _{n}\overset{\mathbb{P}}{\rightarrow }0$ by assumption, Point 1 in the
statement is a direct consequence of the Continuous Mapping Theorem (see
e.g. \cite{Bill}). Moreover, (\ref{A1}) is proved, once it is shown that $%
\left( V^{n},X\right) \overset{\text{law}}{\Rightarrow }\left( W,X\right) $.
To do this, observe first that, for every $n$, $V_{n}$ is a standard
Brownian motion, started from zero and independent of $\mathcal{H}_{\mathcal{%
\tau }_{n}}^{n}$. We shall now show that, as $n\rightarrow +\infty $, for
every $A\in \mathcal{B}\left( C\left[ 0,T\right) \right) $ ($\mathcal{B}%
\left( C\left[ 0,T\right) \right) $ is the Borel $\sigma $-field of the
class $C\left[ 0,T\right) $ of the continuous functions on $\left[
0,T\right) $) and every $B\in \mathcal{S}$,%
\begin{equation}
\mathbb{P}\left[ V^{n}\in A,X\in B\right] -\mathbb{P}\left[ W\in A\right]
\mathbb{P}\left[ X\in B\right] \rightarrow 0.  \label{A2}
\end{equation}%
As a matter of fact, since $\mathcal{U}_{n}\subseteq \mathcal{H}_{\mathcal{%
\tau }_{n}}^{n}$, and thanks to the martingale convergence theorem and the
fact that $X\in \vee _{n}\mathcal{U}_{n}$,%
\begin{eqnarray*}
\mathbb{E}\left\vert \mathbb{P}\left[ X\in B\mid \mathcal{H}_{\mathcal{\tau }%
_{n}}^{n}\right] -\mathbf{1}_{X\in B}\right\vert &\leq &\mathbb{E}\left\vert
\mathbb{P}\left[ X\in B\mid \mathcal{H}_{\mathcal{\tau }_{n}}^{n}\right] -%
\mathbb{P}\left[ X\in B\mid \mathcal{U}_{n}\right] \right\vert \\
&&\mathbf{+}\mathbb{E}\left\vert \mathbf{1}_{X\in B}-\mathbb{P}\left[ X\in
B\mid \mathcal{U}_{n}\right] \right\vert \\
&\leq &2\mathbb{E}\left\vert \mathbf{1}_{X\in B}-\mathbb{P}\left[ X\in B\mid
\mathcal{U}_{n}\right] \right\vert \underset{n\rightarrow +\infty }{%
\rightarrow }0,
\end{eqnarray*}%
thus implying that $\mathbb{P}\left[ X\in B\mid \mathcal{H}_{\mathcal{\tau }%
_{n}}^{n}\right] \overset{L^{1}}{\rightarrow }\mathbf{1}_{X\in B}$, and
therefore%
\begin{eqnarray*}
\mathbb{P}\left[ V^{n}\in A,X\in B\right] -\mathbb{P}\left[ W\in A\right]
\mathbb{P}\left[ X\in B\right] &=&\mathbb{P}\left[ V^{n}\in A,X\in B\right] -%
\mathbb{E}\left[ \mathbf{1}_{V_{n}\in A}\mathbb{P}\left[ X\in B\mid \mathcal{%
H}_{\mathcal{\tau }_{n}}^{n}\right] \right] \\
&&+\mathbb{E}\left[ \mathbf{1}_{V_{n}\in A}\mathbb{P}\left[ X\in B\mid
\mathcal{H}_{\mathcal{\tau }_{n}}^{n}\right] \right] -\mathbb{P}\left[ W\in A%
\right] \mathbb{P}\left[ X\in B\right] \\
&=&\mathbb{P}\left[ V^{n}\in A,X\in B\right] -\mathbb{E}\left[ \mathbf{1}%
_{V_{n}\in A}\mathbb{P}\left[ X\in B\mid \mathcal{H}_{\mathcal{\tau }%
_{n}}^{n}\right] \right] \underset{n\rightarrow +\infty }{\rightarrow }0,
\end{eqnarray*}%
where the last equality follows from the independence of $V_{n}$ and $%
\mathcal{H}_{\mathcal{\tau }_{n}}^{n}$.\ Since (\ref{A2}) implies that $%
\left( V^{n},X\right) \overset{\text{law}}{\Rightarrow }\left( W,X\right) $,
with $W$ and $X$ independent, the proof is concluded.
\end{proof}

\bigskip

\textbf{Proof of Theorem \ref{T : Gen CV} -- }According to Proposition \ref%
{P : SkoAd}-3, for each $n$ the process $t\mapsto \delta _{X_{n}}\left( \pi
_{n,t}u_{n}\right) $ is a continuous square-integrable $\mathcal{F}_{t}^{\pi
_{n}}\left( X_{n}\right) $-martingale with quadratic variation $t\mapsto
\left\Vert \pi _{n,t}u_{n}\right\Vert _{\mathfrak{H}_{n}}^{2}\triangleq \psi
_{n}\left( t\right) $, $t\in \left[ 0,1\right] $. Now define
\begin{equation*}
\mathcal{G}_{s}^{n}\triangleq \mathcal{F}_{\rho _{n,s}}^{\pi _{n}}\text{, \
\ }s\geq 0\text{, \ \ where \ \ }\rho _{n,s}=\inf \left\{ x\in \left[ 0,1%
\right] :\psi _{n}\left( x\right) >s\right\} ,
\end{equation*}%
with $\inf \varnothing =1$, and observe that the above definition is well
given and also that, for every $n\geq 1$, every $t\in \left[ 0,1\right] $
and $s\geq 0$, $\psi _{n}\left( t\right) =\left\Vert \pi
_{n,t}u_{n}\right\Vert _{\mathfrak{H}_{n}}^{2}$ is a $\mathcal{G}^{n}$%
-stopping time and $\rho _{n,s}$ is a $\mathcal{F}^{\pi _{n}}$-stopping
time. In particular, for every $x\geq 0$ and $t\in \left[ 0,1\right] $,
\begin{equation*}
\left\{ \psi _{n}\left( t\right) >x\right\} =\left\{ \rho _{n,x}<t\right\}
\in \mathcal{F}_{t}^{\pi _{n}}.
\end{equation*}%
According to the well known Dambis-Dubins-Schwarz Theorem (see e.g. \cite[%
Ch. V]{ReYor}), the underlying probability space can be suitably enlarged in
order to support a sequence of stochastic processes $W^{n}$ such that, for
each fixed $n$, $W^{n}$ is a $\mathcal{G}^{n}$-Brownian motion started from
zero, and also%
\begin{equation}
\delta _{X_{n}}\left( \pi _{n,t}u_{n}\right) =W_{\psi _{n}\left( t\right)
}^{\left( n\right) }\text{, \ \ }t\in \left[ 0,1\right] \text{. }  \label{*}
\end{equation}%
Since, in general, $\rho _{n,\psi _{n}\left( t\right) }\geq t$, $\mathcal{G}%
_{\psi _{n}\left( t\right) }^{n}\supseteq \mathcal{F}_{t}^{\pi _{n}}$ for
every $t\in \left[ 0,1\right] $. It follows that, for the sequence $t_{n}$
appearing in the statement of Theorem \ref{T : Gen CV},
\begin{equation*}
\mathcal{U}_{n}\subseteq \mathcal{U}_{n+1}\cap \mathcal{F}_{t_{n}}^{\pi
_{n}}\subseteq \mathcal{U}_{n+1}\cap \mathcal{G}_{\psi _{n}\left(
t_{n}\right) }^{n}.
\end{equation*}%
Thus, all conditions of Proposition \ref{P : stopping nest} are verified,
with $\mathcal{H}^{n}=\mathcal{G}^{n}$ and $\tau _{n}=\psi _{n}\left(
t_{n}\right) $, and therefore, for every $\mathcal{U}^{\ast }=\vee _{n}%
\mathcal{U}_{n}$-measurable and real-valued random variable $Z$,
\begin{equation*}
\left( W^{n},Z\right) \overset{\text{law}}{\Rightarrow }\left( W,Z\right)
\text{,}
\end{equation*}%
where $W$ is a Brownian motion independent of $Z$. Moreover, since $\psi
_{n}\left( 1\right) =\left\Vert \pi _{n,1}u_{n}\right\Vert _{\mathfrak{H}%
_{n}}^{2}\overset{\mathbb{P}}{\rightarrow }Y\in \mathcal{U}^{\ast }$ by
assumption, we conclude that, for every $Z\in \mathcal{U}^{\ast }$,
\begin{equation*}
\left( W^{n},Z,\psi _{n}\left( 1\right) \right) \overset{\text{law}}{%
\Rightarrow }\left( W,Z,Y\right) .
\end{equation*}%
Now observe that $\delta _{X_{n}}\left( u_{n}\right) =W_{\psi _{n}\left(
1\right) }^{n}$ by (\ref{*}) and also that, thanks to a further application
of the Continuous Mapping Theorem,
\begin{equation*}
\left( W_{\psi _{n}\left( 1\right) }^{n},Z\right) =\left( \delta
_{X_{n}}\left( u_{n}\right) ,Z\right) \overset{\text{law}}{\Rightarrow }%
\left( W_{Y},Z\right) \text{, }
\end{equation*}%
implying that, for every $\gamma ,\lambda \in \mathbb{R}$,
\begin{equation*}
\mathbb{E}\left[ e^{i\gamma Z}e^{i\lambda \delta _{X_{n}}\left( u_{n}\right)
}\right] \rightarrow \mathbb{E}\left[ e^{i\gamma Z}e^{i\lambda W_{Y}}\right]
=\mathbb{E}\left[ e^{i\gamma Z}e^{-\frac{\lambda ^{2}}{2}Y}\right] \text{,}
\end{equation*}%
which yields the desired conclusion. $\ \ \blacksquare $

\subsection{Further refinements}

The following result is a refinement of Theorem \ref{T : Gen CV} and
Corollary \ref{C : Diff Conv}. It will be used in the next section to
characterize the stable convergence of double Wiener-It\^{o} integrals. It
is proved in \cite[Proposition 10, Theorem 22 and formula (123)]{PeTaq}. The
setting is that of Theorem \ref{T : Gen CV}: $\mathfrak{H}_{n}$, $n\geq 1$,
is a real separable Hilbert space; $X_{n}=X_{n}\left( \mathfrak{H}%
_{n}\right) $, $n\geq 1$, is an isonormal Gaussian process over $\mathfrak{H}%
_{n}$.

\begin{theorem}
\label{T : GenConv 2}Keep the assumptions of Theorem \ref{T : Gen CV} (in
particular, $u_{n}\in L_{\pi _{n}}^{2}\left( \mathfrak{H}_{n},X_{n}\right) $
for every $n$, and (\ref{f1}), (\ref{f2}) and (\ref{f3}) are verified).
Then, as $n\rightarrow +\infty $,%
\begin{equation*}
\mathbb{E}\left[ \exp \left( i\lambda \delta _{X_{n}}\left( u_{n}\right)
\right) \mid \mathcal{F}_{t_{n}}^{\pi _{n}}\left( X_{n}\right) \right]
\overset{\mathbb{P}}{\rightarrow }\exp \left( -\frac{\lambda ^{2}}{2}%
Y\right) ,\text{ \ \ }\forall \lambda \in
\mathbb{R}
\text{.}
\end{equation*}%
Moreover, if there exists a finite random variable $C\left( \omega \right)
>0 $ such that, for some $\eta >0$,%
\begin{equation*}
\mathbb{E}\left[ \left\vert \delta _{X_{n}}\left( u_{n}\right) \right\vert
^{\eta }\mid \mathcal{F}_{t_{n}}^{\pi _{n}}\right] <C\left( \omega \right)
\text{, \ }\forall n\geq 1\text{,\ \ a.s.-}\mathbb{P},
\end{equation*}%
then, there is a subsequence $\left\{ n\left( k\right) :k\geq 1\right\} $
such that, a.s. - $\mathbb{P}$,
\begin{equation*}
\mathbb{E}\left[ \exp \left( i\lambda \delta _{X_{n}}\left( u_{n}\right)
\right) \mid \mathcal{F}_{t_{n\left( k\right) }}^{\pi _{n\left( k\right) }}%
\right] \underset{k\rightarrow +\infty }{\rightarrow }\exp \left( -\frac{%
\lambda ^{2}}{2}Y\right) \text{, \ \ }\forall \lambda \in \mathbb{R}.
\end{equation*}
\end{theorem}

\section{Main results}

Although Corollary \ref{C : Diff Conv} is quite general, the explicit
computation of the projections
\begin{equation*}
proj\left\{ D_{X_{n}}F_{n}\mid L_{\pi _{n}}^{2}\left( \mathfrak{H}%
,X_{n}\right) \right\} ,\text{ \ \ }n\geq 1\text{,}
\end{equation*}%
may be rather difficult. In this section, we prove simpler sufficient
conditions ensuring that the second asymptotic relation in (\ref{f4}) is
satisfied, when $\left( F_{n}\right) $ is a sequence of multiple Wiener-It%
\^{o} integrals of a fixed order. In particular, these conditions do not
involve any projection on the spaces $L_{\pi _{n}}^{2}\left( \mathfrak{H}%
,X_{n}\right) $. The techniques developed below can be suitably extended to
study the joint convergence of vectors of multiple Wiener-It\^{o} integrals.
This issue will be studied in a separate paper.

\subsection{Statements}

To start, fix a real separable Hilbert space $\mathfrak{H}$ and let $%
\{e_{k}:k\geq 1\}$ be a complete orthonormal system in $\mathfrak{H}$. For
every $d\geq 1$, every $p=0,...,d$ and $f\in \mathfrak{H}^{\odot d}$, we
define the \textit{contraction }of $f$ of order $p$ to be the element of $%
\mathfrak{H}^{\otimes 2(d-p)}$ given by
\begin{equation}
f\otimes _{p}f=\sum_{i_{1},\ldots ,i_{p}=1}^{\infty }\ \left\langle
f,e_{i_{1}}\otimes \cdots \otimes e_{i_{p}}\right\rangle _{\mathfrak{H}%
^{\otimes p}}\otimes \left\langle f,e_{i_{1}}\otimes \cdots \otimes
e_{i_{p}}\right\rangle _{\mathfrak{H}^{\otimes p}},  \label{op0}
\end{equation}%
and we denote by $\left( f\otimes _{p}f\right) _{s}$ its symmetrization. As
shown in \cite{NuPe} and \cite{PT}, the asymptotic behavior of the
contractions $f\otimes _{p}f$, $p=1,...,n-1$, plays a crucial role in the
proof of CLTs for multiple Wiener-It\^{o} integrals. To obtain analogous
results in the case of stable convergence, we need to define a further class
of contraction operators, constructed by means of resolutions of the
identity. To this end, fix $\pi \in \mathcal{R}\left( \mathfrak{H}\right) $,
$t\in \left[ 0,1\right] $ and $d\geq 1$, and define $\pi _{t}^{\otimes d}:%
\mathfrak{H}^{\otimes d}\mapsto \mathfrak{H}^{\otimes d}$ to be the $n$th
tensor product of $\pi _{t}$, that is, $\pi _{t}^{\otimes d}$\ is the
projection operator, from $\mathfrak{H}^{\otimes d}$ to itself, given by
\begin{equation}
\pi _{t}^{\otimes d}=\underset{d\text{ times}}{\underbrace{\pi _{t}\otimes
\pi _{t}\otimes \cdot \cdot \cdot \otimes \pi _{t}}}\text{.}  \label{op1}
\end{equation}

For every $d\geq 1$, $p=0,...,n$, $t\in \left[ 0,1\right] $ and $f\in
\mathfrak{H}^{\odot d}$, we write $f\otimes _{p}^{\pi ,t}f$ to indicate the
element of $\mathfrak{H}^{\otimes 2\left( d-p\right) }$ given by%
\begin{equation}
f\otimes _{p}^{\pi ,t}f=\sum_{i_{1},\ldots ,i_{p}=1}^{\infty }\ \left\langle
f,\left( \pi _{1}^{\otimes p}-\pi _{t}^{\otimes p}\right) e_{i_{1}}\otimes
\cdots \otimes e_{i_{p}}\right\rangle _{\mathfrak{H}^{\otimes p}}\otimes
\left\langle f,\left( \pi _{1}^{\otimes p}-\pi _{t}^{\otimes p}\right)
e_{i_{1}}\otimes \cdots \otimes e_{i_{p}}\right\rangle _{\mathfrak{H}%
^{\otimes p}},  \label{op2}
\end{equation}%
and, as before, we denote by $\left( f\otimes _{p}^{\pi ,t}f\right) _{s}$
its symmetrization. We define $f\otimes _{p}^{\pi ,t}f$ to be the \textit{%
generalized contraction kernel} of order $p$, associated to $\pi $ and $t$.
For instance, for $f\in \mathfrak{H}^{\odot d}$,%
\begin{equation}
f\otimes _{p}^{\pi ,0}f=f\otimes _{p}f\text{, \ \ }f\otimes _{p}^{\pi ,1}f=0%
\text{, \ and \ }f\otimes _{d}^{\pi ,t}f=\left\Vert \left( \pi _{1}^{\otimes
d}-\pi _{t}^{\otimes d}\right) f\right\Vert _{\mathfrak{H}^{\otimes d}}^{2}%
\text{.}  \label{op3}
\end{equation}

\bigskip

\textbf{Remark --} When $\mathfrak{H}=L^{2}\left( Z,\mathcal{Z},\nu \right) $%
, where $\nu $ is $\sigma $-finite and non-atomic, and $\pi \in \mathcal{R}%
\left( \mathfrak{H}\right) $ has the form
\begin{equation*}
\pi _{t}f\left( z\right) =f\left( z\right) \mathbf{1}_{Z_{t}}\left( z\right)
\text{, \ \ }z\in Z\text{,}
\end{equation*}%
where $Z_{t}$ is an increasing sequence in $\mathcal{Z}$ such that $%
Z_{0}=\varnothing $ and $Z_{1}=Z$, we have the following elementary
relation: for every $d\geq 1$, $p=0,...,d$, $t\in \left[ 0,1\right] $ and $%
f\in \mathfrak{H}^{\odot d}=L_{s}^{2}\left( Z^{\otimes d},\mathcal{Z}%
^{\otimes d},\nu ^{\otimes d}\right) $,
\begin{equation*}
f\otimes _{p}^{\pi ,t}f\left( z_{1},...,z_{2\left( d-p\right) }\right)
=\int_{Z^{p}\backslash Z_{t}^{p}}f\left( z_{1},...,z_{d-p},\mathbf{x}%
_{p}\right) f\left( z_{d-p+1},...,z_{2\left( d-p\right) },\mathbf{x}%
_{p}\right) \nu ^{\otimes p}\left( d\mathbf{x}_{p}\right) ,
\end{equation*}%
since, for every $h\in L_{s}^{2}\left( Z^{\otimes p},\mathcal{Z}^{\otimes
p},\nu ^{\otimes p}\right) $,%
\begin{equation*}
\left( \pi _{1}^{\otimes p}-\pi _{t}^{\otimes p}\right) h\left( \mathbf{x}%
_{p}\right) =\left( \mathbf{1}_{Z^{p}}\left( \mathbf{x}_{p}\right) -\mathbf{1%
}_{Z_{t}^{p}}\left( \mathbf{x}_{p}\right) \right) h\left( \mathbf{x}%
_{p}\right) =\mathbf{1}_{Z^{p}\backslash Z_{t}^{p}}\left( \mathbf{x}%
_{p}\right) h\left( \mathbf{x}_{p}\right) \text{.}
\end{equation*}

\bigskip

The next result, which is the main achievement of the paper, generalizes the
crucial part of the CLT stated in \cite[Theorem 1]{NuPe} to the case of the
stable convergence. Its proof is postponed to the next section.

\bigskip

\begin{theorem}
\label{T : Main}Let $\mathfrak{H}_{n}$, $X_{n}\left( \mathfrak{H}_{n}\right)
$, $\pi _{n}$, $t_{n}$ and  $\mathcal{U}_{n}$, $n\geq 1$, satisfy the
assumptions of Corollary \ref{C : Diff Conv}. Fix $d\geq 2$, and consider a
sequence of random variables $\left\{ F_{n}:n\geq 1\right\} $, such that,
for every $n$,
\begin{equation}
F_{n}=I_{d}^{X_{n}}\left( f_{n}\right) \text{,}  \label{MWIintegr}
\end{equation}%
for a certain $f_{n}\in \mathfrak{H}_{n}^{\odot d}$, and moreover%
\begin{equation}
\mathbb{E}\left[ F_{n}\mid \mathcal{F}_{t_{n}}^{\pi _{n}}\left( X_{n}\right) %
\right] =I_{d}^{X_{n}}\left( \pi _{n,t_{n}}^{\otimes d}f_{n}\right) \overset{%
L^{2}\left( \mathbb{P}\right) }{\rightarrow }0  \label{NEG}
\end{equation}%
and
\begin{equation}
\mathbb{E}\left[ F_{n}^{2}\mid \mathcal{F}_{t_{n}}^{\pi _{n}}\left(
X_{n}\right) \right] \overset{\mathbb{P}}{\rightarrow }Y\in \mathcal{U}%
^{\ast }.  \label{NORM}
\end{equation}%
Then, the following hold:

\begin{enumerate}
\item for every $n\geq 1$,
\begin{eqnarray}
\mathbb{E}\left[ F_{n}^{2}\mid \mathcal{F}_{t_{n}}^{\pi _{n}}\left(
X_{n}\right) \right] &=&d!\left\Vert f_{n}\right\Vert _{\mathfrak{H}%
_{n}^{\otimes d}}^{2}  \label{DEV} \\
&&+\sum_{r=1}^{d-1}\left( d-r\right) !\dbinom{d}{r}^{2}I_{2r}^{X_{n}}\left[
\pi _{n,t_{n}}^{\otimes 2r}\left( f_{n}\otimes _{d-r}^{\pi
_{n},t_{n}}f_{n}\right) \right] +o_{\mathbb{P}}\left( 1\right) ,  \notag
\end{eqnarray}%
where $o_{\mathbb{P}}\left( 1\right) $ stands for a sequence of random
variables converging to zero in probability;

\item if $\pi _{n}\in \mathcal{R}_{AC}\left( \mathfrak{H}_{n}\right) $ for
every $n\geq 1$ and, for every $r=1,...,d-1,$%
\begin{equation}
\left\Vert \left( \pi _{n,1}^{\otimes 2r}-\pi _{n,t_{n}}^{\otimes 2r}\right)
\left( f_{n}\otimes _{d-r}^{\pi _{n},t_{n}}f_{n}\right) \right\Vert _{%
\mathfrak{H}_{n}^{\otimes 2r}}^{2}\underset{n\rightarrow +\infty }{%
\rightarrow }0\text{,}  \label{asCV}
\end{equation}%
then%
\begin{equation}
\left\Vert proj\left\{ D_{X_{n}}F_{n}\mid L_{\pi _{n}}^{2}\left( \mathfrak{H}%
_{n},X_{n}\right) \right\} \right\Vert _{\mathfrak{H}_{n}}^{2}\overset{%
\mathbb{P}}{\underset{n\rightarrow +\infty }{\rightarrow }}Y\text{,}
\label{CVDIFF}
\end{equation}%
and therefore%
\begin{equation*}
F_{n}\rightarrow _{\left( s,\mathcal{U}^{\ast }\right) }\mathbb{E}\mu \left(
\cdot \right) \text{ \ \ and \ \ }\mathbb{E}\left[ \exp \left( i\lambda
F_{n}\right) \mid \mathcal{F}_{t_{n}}^{\pi _{n}}\left( X_{n}\right) \right]
\overset{\mathbb{P}}{\rightarrow }\exp \left( -\frac{\lambda ^{2}}{2}%
Y\right) ,\text{ \ \ }\forall \lambda \in
\mathbb{R}
\text{,}
\end{equation*}%
where $\widehat{\mu }\left( \lambda \right) =\exp \left( -\frac{\lambda ^{2}%
}{2}Y\right) $, $\forall \lambda \in \mathbb{R}$.
\end{enumerate}
\end{theorem}

\bigskip

\textbf{Remarks -- }(a) Since, due to Lemma \ref{L : AC} and for any
continuous $\pi \in \mathcal{R}\left( \mathfrak{H}\right) $, there exists a
non decreasing function $\phi $ such that $\widetilde{\pi }\triangleq \pi
_{\phi \left( \cdot \right) }$ is absolutely continuous, Theorem \ref{T :
Main} applies \textit{de facto }to any sequence $\pi _{n}\in \mathcal{R}%
\left( \mathfrak{H}_{n}\right) $, $n\geq 1$.

(b) Suppose that $X_{n}\left( \mathfrak{H}_{n}\right) =X\left( \mathfrak{H}%
\right) $ for every $n\geq 1$. Then, the random variables $F_{n}$, $n\geq 1$%
, appearing in (\ref{MWIintegr}) all belong to the same Wiener chaos, and,
due to (\ref{DEV}), the sequence $\mathbb{E}\left[ F_{n}^{2}\mid \mathcal{F}%
_{t_{n}}^{\pi _{n}}\left( X_{n}\right) \right] $, $n\geq 1$, belongs to the
same finite sum of $d$ Wiener chaoses. Recall also that, on a finite sum of
Wiener chaoses, the topology induced by convergence in probability is
equivalent to the $L^{p}$ topology, for every $p\geq 1$ (see e.g. \cite{Sch}%
). When (\ref{NEG}) and (\ref{NORM}) are verified, we therefore deduce from (%
\ref{DEV}) that $Y$ has necessarily the form
\begin{equation}
Y=\mathbb{E}\left( Y\right) +\sum_{r=1}^{d-1}I_{2r}^{X}\left( g_{r}\right)
\text{,}  \label{Yform}
\end{equation}%
for some $g_{r}\in \mathfrak{H}^{\odot r}$ and $r=1,...,d-1$. Moreover, (\ref%
{NORM}) is equivalent to the condition: as $n\rightarrow +\infty $, $%
d!\left\Vert f_{n}\right\Vert _{\mathfrak{H}^{\otimes d}}^{2}\rightarrow
\mathbb{E}\left( Y\right) $ and, for $r=1,...,d-1$,
\begin{equation}
\left( d-r\right) !\dbinom{d}{r}^{2}\times \left( \pi _{n,t_{n}}^{\otimes
2r}\left( f_{n}\otimes _{d-r}^{\pi _{n},t_{n}}f_{n}\right) \right)
_{s}\rightarrow g_{r}  \label{convEX}
\end{equation}%
in $\mathfrak{H}^{\odot 2r}$. It follows that, for $r=1,...,d-1$, the two
operators $f_{n}\otimes _{d-r}^{\ast }f_{n}$ and $f_{n}\otimes _{d-r}^{\ast
\ast }f_{n}$, from $\mathfrak{H}\ $to $\mathfrak{H}^{\otimes 2r}$, defined as%
\begin{eqnarray}
f_{n} &\mapsto &\left( d-r\right) !\dbinom{d}{r}^{2}\left( \pi
_{n,t_{n}}^{\otimes 2r}\left( f_{n}\otimes _{d-r}^{\pi
_{n},t_{n}}f_{n}\right) \right) _{s}\triangleq f_{n}\otimes _{d-r}^{\ast
}f_{n}  \label{backtoIntro1} \\
f_{n} &\mapsto &\left( \pi _{n,1}^{\otimes 2r}-\pi _{n,t_{n}}^{\otimes
2r}\right) \left( f_{n}\otimes _{d-r}^{\pi _{n},t_{n}}f_{n}\right)
\triangleq f_{n}\otimes _{d-r}^{\ast \ast }f_{n},  \label{backtoIntro2}
\end{eqnarray}%
solve the problem raised in the Introduction. Indeed, under (\ref{NEG}) and
the normalization condition $d!\left\Vert f_{n}\right\Vert _{\mathfrak{H}%
^{\otimes d}}^{2}\rightarrow \mathbb{E}\left( Y\right) $, due to Theorem \ref%
{T : Main} and (\ref{convEX}), the asymptotic relation (\ref{CVintro})
implies that $I_{d}^{X}\left( f_{n}\right) $ converges stably to $\sqrt{Y}%
\times N$, where $N$ is a centered standard Gaussian random variable
independent of $Y$.

\bigskip

We now show that the conclusions of Theorem \ref{T : Main} may be
strengthened in the case of a sequence of double Wiener-It\^{o} integrals,
i.e. in the case $d=2$. The proof of the next result\ is deferred to the
Section 3.3.

\bigskip

\begin{theorem}
\label{T : Main Double}Under the assumptions and notation of Theorem \ref{T
: Main} (in particular, (\ref{NEG}) and (\ref{NORM}) are in order), suppose
that $d=2$ and that the following implication holds:
\begin{equation}
\mathbb{E}\left( \mathbb{E}\left[ F_{n}^{4}\mid \mathcal{F}_{t_{n}}^{\pi
_{n}}\left( X_{n}\right) \right] -3Y^{2}\right) ^{2}\underset{n\rightarrow
+\infty }{\rightarrow }0\text{ \ \ if, and only if, \ \ }\mathbb{E}\left[
F_{n}^{4}\mid \mathcal{F}_{t_{n}}^{\pi _{n}}\left( X_{n}\right) \right]
\overset{\mathbb{P}}{\rightarrow }3Y^{2}.  \label{L2/PRO}
\end{equation}%
Then, the following are equivalent

\begin{description}
\item[(i)] $\left\Vert proj\left\{ D_{X_{n}}F_{n}\mid L_{\pi _{n}}^{2}\left(
\mathfrak{H}_{n},X_{n}\right) \right\} \right\Vert _{\mathfrak{H}_{n}}^{2}%
\overset{\mathbb{P}}{\underset{n\rightarrow +\infty }{\rightarrow }}Y$;

\item[(ii)] $\mathbb{E}\left[ \exp \left( i\lambda F_{n}\right) \mid
\mathcal{F}_{t_{n}}^{\pi _{n}}\left( X_{n}\right) \right] \overset{\mathbb{P}%
}{\rightarrow }\exp \left( -\frac{\lambda ^{2}}{2}Y\right) ,$ \ \ $\forall
\lambda \in
\mathbb{R}
$;

\item[(iii)] $\mathbb{E}\left[ F_{n}^{4}\mid \mathcal{F}_{t_{n}}^{\pi
_{n}}\left( X_{n}\right) \right] \overset{\mathbb{P}}{\rightarrow }3Y^{2}$;

\item[(iv)] $\left\Vert \left( \pi _{n,1}^{\otimes 2}-\pi
_{n,t_{n}}^{\otimes 2}\right) \left( f_{n}\otimes _{1}^{\pi
_{n},t_{n}}f_{n}\right) \right\Vert _{\mathfrak{H}_{n}^{\odot 2}}^{2}%
\underset{n\rightarrow +\infty }{\rightarrow }0.$
\end{description}

Moreover, if either one of conditions (i)-(iv) is satisfied, $%
F_{n}\rightarrow _{\left( s,\mathcal{U}^{\ast }\right) }\mathbb{E}\mu \left(
\cdot \right) $, where $\widehat{\mu }\left( \lambda \right) =\exp \left( -%
\frac{\lambda ^{2}}{2}Y\right) $.
\end{theorem}

\bigskip

\textbf{Remark -- }(a) Due again to the equivalence of the $L^{0}$ and $%
L^{2} $ topology on a finite sum of Wiener chaoses (see \cite{Sch}),
condition (\ref{L2/PRO}) is verified in the case $\mathfrak{H}_{n}=\mathfrak{%
H}$ and $X_{n}\left( \mathfrak{H}\right) =X\left( \mathfrak{H}\right) $, for
every $n\geq 1.$

(b) When $d=2$, the second part of Theorem \ref{T : Main} corresponds to the
implications (iv) $\Longrightarrow $ (i) $\Longrightarrow $ (ii) of Theorem %
\ref{T : Main Double}.

\bigskip

The next consequence of Theorem \ref{T : Main} is a central limit theorem,
generalizing Theorem \ref{T : intro}.

\bigskip

\begin{corollary}
Let $\mathfrak{H}_{n}$, $X_{n}\left( \mathfrak{H}_{n}\right) $, $n\geq 1$,
be defined as above, and suppose that $\pi _{n}\in \mathcal{R}_{AC}\left(
\mathfrak{H}_{n}\right) $ for each $n$. For $d\geq 2$, consider a sequence
of multiple Wiener-It\^{o} integrals $\left\{ I_{d}^{X_{n}}\left(
f_{n}\right) :n\geq 1\right\} $ s.t. $f_{n}\in \mathfrak{H}_{n}^{\odot d}$,
and
\begin{equation*}
\mathbb{E}\left[ I_{d}^{X_{n}}\left( f_{n}\right) ^{2}\right] =d!\left\Vert
f_{n}\right\Vert _{\mathfrak{H}_{n}^{\otimes d}}^{2}\underset{n\rightarrow
+\infty }{\rightarrow }1.
\end{equation*}%
Then, the following are equivalent

\begin{description}
\item[(i)] $\left\Vert proj\left\{ D_{X_{n}}I_{d}^{X_{n}}\left( f_{n}\right)
\mid L_{\pi _{n}}^{2}\left( \mathfrak{H}_{n},X_{n}\right) \right\}
\right\Vert _{\mathfrak{H}_{n}}^{2}\overset{\mathbb{P}}{\underset{%
n\rightarrow +\infty }{\rightarrow }}1$;

\item[(ii)] $\mathbb{E}\left[ \exp \left( i\lambda I_{d}^{X_{n}}\left(
f_{n}\right) \right) \right] \underset{n\rightarrow +\infty }{\rightarrow }%
\exp \left( -\frac{\lambda ^{2}}{2}\right) $,\ $\forall \lambda \in
\mathbb{R}
$, that is, $I_{d}^{X_{n}}\left( f_{n}\right) \overset{\text{law}}{%
\rightarrow }N\left( 0,1\right) $;

\item[(iii)] $\mathbb{E}\left[ I_{d}^{X_{n}}\left( f_{n}\right) ^{4}\right]
\underset{n\rightarrow +\infty }{\rightarrow }3$;

\item[(iv)] $\left\Vert f_{n}\otimes _{d-r}f_{n}\right\Vert _{\mathfrak{H}%
_{n}^{\odot 2r}}^{2}\underset{n\rightarrow +\infty }{\rightarrow }0$, $%
\forall r=1,...,d-1$.
\end{description}
\end{corollary}

\begin{proof}
The equivalence of the three conditions (ii)-(iv) is the object of Theorem %
\ref{T : intro}. That (i) implies (ii) follows from Corollary \ref{C : Diff
Conv}, in the case $t_{n}=0$ for every $n$. Finally, (iv) implies (i) thanks
to Theorem \ref{T : Main}-2, again in the case $t_{n}=0$.
\end{proof}

\subsection{Proof of Theorem \protect\ref{T : Main}}

We start by proving an auxiliary analytic result. Let $\left( A,\mathcal{A}%
\right) $ be a measurable space. For $m\geq 1$, let $\mathbf{\xi }_{m}$ be
shorthand for a vector $\mathbf{\xi }_{m}=\left( \left( a_{1},x_{1}\right)
;...,\left( a_{m};x_{m}\right) \right) \in \left( A\times \left[ 0,1\right]
\right) ^{m}$ and, for such $\mathbf{\xi }_{m}$, note $\overline{\mathbf{\xi
}}_{m}$ the maximum of $\mathbf{\xi }_{m}$ in the variables $x_{1},...,x_{m}$%
, i.e. $\overline{\mathbf{\xi }}_{m}=\overline{\left( \left(
a_{1},x_{1}\right) ;...,\left( a_{m};x_{m}\right) \right) }%
=\max_{i=1,...,n}\left( x_{i}\right) $. In what follows, $l\left( d\xi
\right) $ stands for a $\sigma $-finite positive measure on $A\times \left[
0,1\right] $, such that, for every fixed $x^{\ast }\in \left[ 0,1\right] $, $%
l\left\{ \left( a,x\right) :x=x^{\ast }\right\} =0$ (note that this implies
that $l$ is non-atomic). For $m\geq 1$, $l^{m}\left( d\mathbf{\xi }%
_{m}\right) $ is the canonical product measure on $\left( A\times \left[ 0,1%
\right] \right) ^{m}$ (with $l^{1}=l$ by convention).

\bigskip

\begin{lemma}
\label{L : Sq. Eq.}Let $m,r\geq 1$, $C\subseteq \left[ 0,1\right] ^{m}$ and $%
D\subseteq \left[ 0,1\right] ^{r}$. Then, for every symmetric function $f\in
L_{s}^{2}\left( \left( A\times \left[ 0,1\right] \right) ^{m+r},l
^{m+r}\right) \triangleq L_{s}^{2}\left( l ^{m+r}\right) $%
\begin{eqnarray}
&&\int_{A\times C}\int_{A\times C}\left[ \int_{A\times D}f\left( \mathbf{\xi
}_{m},\mathbf{\alpha }_{r}\right) f\left( \mathbf{\gamma }_{m},\mathbf{%
\alpha }_{r}\right) l^{r}\left( d\mathbf{\alpha }_{r}\right) \right]
^{2}l^{m}\left( d\mathbf{\xi }_{m}\right) l^{m}\left( d\mathbf{\gamma }%
_{m}\right)  \notag \\
&=&\int_{A\times C}\int_{A\times C}\left[ \int_{A\times D}f\left( \mathbf{%
\xi }_{m},\mathbf{\alpha }_{r}\right) f\left( \mathbf{\gamma }_{m},\mathbf{%
\alpha }_{r}\right) \mathbf{1}_{\left( \overline{\mathbf{\alpha }}_{r}<\max
\left( \overline{\mathbf{\xi }}_{m},\overline{\mathbf{\gamma }}_{m}\right)
\right) }l^{r}\left( d\mathbf{\alpha }_{r}\right) \right] ^{2}l^{m}\left( d%
\mathbf{\xi }_{m}\right) l^{m}\left( d\mathbf{\gamma }_{m}\right)
\label{SqEq} \\
&&+\int_{A\times D}\int_{A\times D}\left[ \int_{A\times C}f\left( \mathbf{%
\xi }_{m},\mathbf{\alpha }_{r}\right) f\left( \mathbf{\xi }_{m},\mathbf{%
\beta }_{r}\right) \mathbf{1}_{\left( \overline{\mathbf{\xi }}_{m}<\max
\left( \overline{\mathbf{\alpha }}_{r},\overline{\mathbf{\beta }}_{r}\right)
\right) }l ^{m}\left( d\mathbf{\xi }_{m}\right) \right] ^{2}l^{r}\left( d%
\mathbf{\alpha }_{r}\right) l^{r}\left( d\mathbf{\beta }_{r}\right) .  \notag
\end{eqnarray}
\end{lemma}

\begin{proof}
Start by writing
\begin{eqnarray}
&&\int_{A\times C}\int_{A\times C}\left[ \int_{A\times D}f\left( \mathbf{\xi
}_{m},\mathbf{\alpha }_{r}\right) f\left( \mathbf{\gamma }_{m},\mathbf{%
\alpha }_{r}\right) l^{r}\left( d\mathbf{\alpha }_{r}\right) \right]
^{2}l^{m}\left( d\mathbf{\xi }_{m}\right) l^{m}\left( d\mathbf{\gamma }%
_{m}\right)   \notag \\
&=&\int_{A\times C}\int_{A\times C}\left[ \int_{A\times D}f\left( \mathbf{%
\xi }_{m},\mathbf{\alpha }_{r}\right) f\left( \mathbf{\gamma }_{m},\mathbf{%
\alpha }_{r}\right) \right. \times   \label{nnat} \\
&&\left. \times \left( \mathbf{1}_{\left( \overline{\mathbf{\alpha }}%
_{r}<\max \left( \overline{\mathbf{\xi }}_{m},\overline{\mathbf{\gamma }}%
_{m}\right) \right) }+\mathbf{1}_{\left( \overline{\mathbf{\alpha }}%
_{r}>\max \left( \overline{\mathbf{\xi }}_{m},\overline{\mathbf{\gamma }}%
_{m}\right) \right) }\right) l^{r}\left( d\mathbf{\alpha }_{r}\right)
^{^{^{{}}}}\right] ^{2}l^{m}\left( d\mathbf{\xi }_{m}\right) l^{m}\left( d%
\mathbf{\gamma }_{m}\right)   \notag \\
&=&\int_{A\times C}\int_{A\times C}\left[ \int_{A\times D}f\left( \mathbf{%
\xi }_{m},\mathbf{\alpha }_{r}\right) f\left( \mathbf{\gamma }_{m},\mathbf{%
\alpha }_{r}\right) \mathbf{1}_{\left( \overline{\mathbf{\alpha }}_{r}<\max
\left( \overline{\mathbf{\xi }}_{m},\overline{\mathbf{\gamma }}_{m}\right)
\right) }l^{r}\left( d\mathbf{\alpha }_{r}\right) \right] ^{2}l^{m}\left( d%
\mathbf{\xi }_{m}\right) l^{m}\left( d\mathbf{\gamma }_{m}\right)   \notag \\
&&+\int_{A\times C}\int_{A\times C}\left[ \int_{A\times D}f\left( \mathbf{%
\xi }_{m},\mathbf{\alpha }_{r}\right) f\left( \mathbf{\gamma }_{m},\mathbf{%
\alpha }_{r}\right) \mathbf{1}_{\left( \overline{\mathbf{\alpha }}_{r}>\max
\left( \overline{\mathbf{\xi }}_{m},\overline{\mathbf{\gamma }}_{m}\right)
\right) }l^{r}\left( d\mathbf{\alpha }_{r}\right) \right] ^{2}l^{m}\left( d%
\mathbf{\xi }_{m}\right) l^{m}\left( d\mathbf{\gamma }_{m}\right)   \notag \\
&&+2\int_{A\times C}\int_{A\times C}\left[ \int_{A\times D}f\left( \mathbf{%
\xi }_{m},\mathbf{\alpha }_{r}\right) f\left( \mathbf{\gamma }_{m},\mathbf{%
\alpha }_{r}\right) \mathbf{1}_{\left( \overline{\mathbf{\alpha }}_{r}<\max
\left( \overline{\mathbf{\xi }}_{m},\overline{\mathbf{\gamma }}_{m}\right)
\right) }l^{r}\left( d\mathbf{\alpha }_{r}\right) \right] \times   \notag \\
&&\text{ \ \ \ \ \ \ \ \ \ \ \ \ \ \ \ \ \ \ \ \ }\times \left[
\int_{A\times D}f\left( \mathbf{\xi }_{m},\mathbf{\beta }_{r}\right) f\left(
\mathbf{\gamma }_{m},\mathbf{\beta }_{r}\right) \mathbf{1}_{\left( \overline{%
\mathbf{\beta }}_{r}>\max \left( \overline{\mathbf{\xi }}_{m},\overline{%
\mathbf{\gamma }}_{m}\right) \right) }l^{r}\left( d\mathbf{\beta }%
_{r}\right) \right] l^{m}\left( d\mathbf{\xi }_{m}\right) l^{m}\left( d%
\mathbf{\gamma }_{m}\right)   \notag \\
&\triangleq &L\left( 1\right) +L\left( 2\right) +L\left( 3\right)   \notag
\end{eqnarray}%
(note that the equality (\ref{nnat}) holds because of the assumption: $%
l\left\{ \left( a,x\right) :x=x^{\ast }\right\} =0$, $\forall x^{\ast }$).
Now, by using a standard Fubini theorem,%
\begin{equation*}
L\left( 2\right) =\int_{A\times D}\int_{A\times D}\left[ \int_{A\times
C}f\left( \mathbf{\gamma }_{m},\mathbf{\alpha }_{r}\right) f\left( \mathbf{%
\gamma }_{m},\mathbf{\beta }_{r}\right) \mathbf{1}_{\left( \overline{\mathbf{%
\gamma }}_{m}<\min \left( \overline{\mathbf{\alpha }}_{r},\overline{\mathbf{%
\beta }}_{r}\right) \right) }l^{m}\left( d\mathbf{\gamma }_{m}\right) \right]
^{2}l^{r}\left( d\mathbf{\alpha }_{r}\right) l^{r}\left( d\mathbf{\beta }%
_{r}\right) ,
\end{equation*}%
and also
\begin{eqnarray*}
L\left( 3\right)  &=&\int_{A\times C}\int_{A\times C}\left[ \int_{A\times
D}f\left( \mathbf{\xi }_{m},\mathbf{\alpha }_{r}\right) f\left( \mathbf{%
\gamma }_{m},\mathbf{\alpha }_{r}\right) \mathbf{1}_{\left( \min \left(
\overline{\mathbf{\alpha }}_{r},\overline{\mathbf{\beta }}_{r}\right) <\max
\left( \overline{\mathbf{\xi }}_{m},\overline{\mathbf{\gamma }}_{m}\right)
\right) }l^{r}\left( d\mathbf{\alpha }_{r}\right) \right] \times  \\
&&\text{ \ \ \ \ }\times \left[ \int_{A\times D}f\left( \mathbf{\xi }_{m},%
\mathbf{\beta }_{r}\right) f\left( \mathbf{\gamma }_{m},\mathbf{\beta }%
_{r}\right) \mathbf{1}_{\left( \max \left( \overline{\mathbf{\alpha }}_{r},%
\overline{\mathbf{\beta }}_{r}\right) >\max \left( \overline{\mathbf{\xi }}%
_{m},\overline{\mathbf{\gamma }}_{m}\right) \right) }l^{r}\left( d\mathbf{%
\beta }_{r}\right) \right] l^{m}\left( d\mathbf{\xi }_{m}\right) l^{m}\left(
d\mathbf{\gamma }_{m}\right)  \\
&=&\int_{A\times D}\int_{A\times D}\left[ \int_{A\times C}f\left( \mathbf{%
\xi }_{m},\mathbf{\alpha }_{r}\right) f\left( \mathbf{\xi }_{m},\mathbf{%
\beta }_{r}\right) \right. \times  \\
&&\text{ \ \ \ \ \ \ \ \ \ }\times \left. \mathbf{1}_{\left( \min \left(
\overline{\mathbf{\alpha }}_{r},\overline{\mathbf{\beta }}_{r}\right) <%
\overline{\mathbf{\xi }}_{m}<\max \left( \overline{\mathbf{\alpha }}_{r},%
\overline{\mathbf{\beta }}_{r}\right) \right) }l^{m}\left( d\mathbf{\xi }%
_{m}\right) ^{^{{}}}\right] ^{2}l^{r}\left( d\mathbf{\alpha }_{r}\right)
l^{r}\left( d\mathbf{\beta }_{r}\right) + \\
&&+2\int_{A\times D}\int_{A\times D}\left[ \int_{A\times C}f\left( \mathbf{%
\xi }_{m},\mathbf{\alpha }_{r}\right) f\left( \mathbf{\xi }_{m},\mathbf{%
\beta }_{r}\right) \mathbf{1}_{\left( \min \left( \overline{\mathbf{\alpha }}%
_{r},\overline{\mathbf{\beta }}_{r}\right) <\overline{\mathbf{\xi }}%
_{m}<\max \left( \overline{\mathbf{\alpha }}_{r},\overline{\mathbf{\beta }}%
_{r}\right) \right) }l^{m}\left( d\mathbf{\xi }_{m}\right) \right] \times  \\
&&\text{ \ \ \ \ \ \ \ \ \ \ \ \ \ \ \ \ \ }\left[ \int_{A\times C}f\left(
\mathbf{\gamma }_{m},\mathbf{\alpha }_{r}\right) f\left( \mathbf{\gamma }%
_{m},\mathbf{\beta }_{r}\right) \mathbf{1}_{\left( \overline{\mathbf{\gamma }%
}_{m}<\min \left( \overline{\mathbf{\alpha }}_{r},\overline{\mathbf{\beta }}%
_{r}\right) \right) }l^{m}\left( d\mathbf{\gamma }_{m}\right) \right]
l^{r}\left( d\mathbf{\alpha }_{r}\right) l^{r}\left( d\mathbf{\beta }%
_{r}\right) .
\end{eqnarray*}%
The last relation implies that%
\begin{equation*}
L\left( 2\right) +L\left( 3\right) =\int_{A\times D}\int_{A\times D}\left[
\int_{A\times C}f\left( \mathbf{\xi }_{m},\mathbf{\alpha }_{r}\right)
f\left( \mathbf{\xi }_{m},\mathbf{\beta }_{r}\right) \mathbf{1}_{\left(
\overline{\mathbf{\xi }}_{m}<\max \left( \overline{\mathbf{\alpha }}_{r},%
\overline{\mathbf{\beta }}_{r}\right) \right) }l^{m}\left( d\mathbf{\xi }%
_{m}\right) \right] ^{2}l^{r}\left( d\mathbf{\alpha }_{r}\right) l^{r}\left(
d\mathbf{\beta }_{r}\right) \text{,}
\end{equation*}%
hence proving (\ref{SqEq}).
\end{proof}

\bigskip

\textbf{Remark }-- With the notation of Lemma \ref{L : Sq. Eq.}, suppose
that the sequence $f_{n}\in L_{s}^{2}\left( l^{m+r}\right) $, $n\geq 1$, is
such that $\left\{ \left\Vert f_{n}\right\Vert _{L_{s}^{2}\left( \zeta
^{m+r}\right) }:n\geq 1\right\} $ is bounded and, as $n\rightarrow +\infty $%
,
\begin{eqnarray}
&&\int_{A\times C}\int_{A\times C}\left[ \int_{A\times D}f_{n}\left( \mathbf{%
\xi }_{m},\mathbf{\alpha }_{r}\right) f_{n}\left( \mathbf{\gamma }_{m},%
\mathbf{\alpha }_{r}\right) l^{r}\left( d\mathbf{\alpha }_{r}\right) \right]
^{2}l^{m}\left( d\mathbf{\xi }_{m}\right) l^{m}\left( d\mathbf{\gamma }%
_{m}\right)   \notag \\
&=&\int_{A\times D}\int_{A\times D}\left[ \int_{A\times C}f_{n}\left(
\mathbf{\xi }_{m},\mathbf{\alpha }_{r}\right) f_{n}\left( \mathbf{\xi }_{m},%
\mathbf{\beta }_{r}\right) l^{m}\left( d\mathbf{\xi }_{m}\right) \right]
^{2}l^{r}\left( d\mathbf{\alpha }_{r}\right) l^{r}\left( d\mathbf{\beta }%
_{r}\right)   \label{i} \\
&\rightarrow &0  \notag
\end{eqnarray}%
(note that the equality in (\ref{i}) derives from a standard Fubini
theorem). Then, by (\ref{nnat}) and (\ref{SqEq}), Lemma \ref{L : Sq. Eq.}
implies that the sequence $Q_{i}\left( n\right) $, defined for $i=1,2,3,4$ by%
\begin{eqnarray*}
Q_{1}\left( n\right)  &=&\int_{A\times C}\int_{A\times C}\left[
\int_{A\times D}f\left( \mathbf{\xi }_{m},\mathbf{\alpha }_{r}\right)
f\left( \mathbf{\gamma }_{m},\mathbf{\alpha }_{r}\right) \mathbf{1}_{\left(
\overline{\mathbf{\alpha }}_{r}<\max \left( \overline{\mathbf{\xi }}_{m},%
\overline{\mathbf{\gamma }}_{m}\right) \right) }l^{r}\left( d\mathbf{\alpha }%
_{r}\right) \right] ^{2}l^{m}\left( d\mathbf{\xi }_{m}\right) l^{m}\left( d%
\mathbf{\gamma }_{m}\right) \text{,} \\
Q_{2}\left( n\right)  &=&\int_{A\times D}\int_{A\times D}\left[
\int_{A\times C}f\left( \mathbf{\xi }_{m},\mathbf{\alpha }_{r}\right)
f\left( \mathbf{\xi }_{m},\mathbf{\beta }_{r}\right) \mathbf{1}_{\left(
\overline{\mathbf{\xi }}_{m}<\max \left( \overline{\mathbf{\alpha }}_{r},%
\overline{\mathbf{\beta }}_{r}\right) \right) }l^{m}\left( d\mathbf{\xi }%
_{m}\right) \right] ^{2}l^{r}\left( d\mathbf{\alpha }_{r}\right) l^{r}\left(
d\mathbf{\beta }_{r}\right) , \\
Q_{3}\left( n\right)  &=&\int_{A\times C}\int_{A\times C}\left[
\int_{A\times D}f\left( \mathbf{\xi }_{m},\mathbf{\alpha }_{r}\right)
f\left( \mathbf{\gamma }_{m},\mathbf{\alpha }_{r}\right) \mathbf{1}_{\left(
\overline{\mathbf{\alpha }}_{r}>\max \left( \overline{\mathbf{\xi }}_{m},%
\overline{\mathbf{\gamma }}_{m}\right) \right) }l^{r}\left( d\mathbf{\alpha }%
_{r}\right) \right] ^{2}l^{m}\left( d\mathbf{\xi }_{m}\right) l^{m}\left( d%
\mathbf{\gamma }_{m}\right)  \\
&=&\int_{A\times D}\int_{A\times D}\left[ \int_{A\times C}f\left( \mathbf{%
\xi }_{m},\mathbf{\alpha }_{r}\right) f\left( \mathbf{\xi }_{m},\mathbf{%
\beta }_{r}\right) \mathbf{1}_{\left( \overline{\mathbf{\xi }}_{m}<\min
\left( \overline{\mathbf{\alpha }}_{r},\overline{\mathbf{\beta }}_{r}\right)
\right) }l^{m}\left( d\mathbf{\xi }_{m}\right) \right] ^{2}l^{r}\left( d%
\mathbf{\alpha }_{r}\right) l^{r}\left( d\mathbf{\beta }_{r}\right) , \\
Q_{4}\left( n\right)  &=&\int_{A\times D}\int_{A\times D}\left[
\int_{A\times C}f\left( \mathbf{\xi }_{m},\mathbf{\alpha }_{r}\right)
f\left( \mathbf{\xi }_{m},\mathbf{\beta }_{r}\right) \mathbf{1}_{\left(
\overline{\mathbf{\xi }}_{m}>\max \left( \overline{\mathbf{\alpha }}_{r},%
\overline{\mathbf{\beta }}_{r}\right) \right) }l^{m}\left( d\mathbf{\xi }%
_{m}\right) \right] ^{2}l^{r}\left( d\mathbf{\alpha }_{r}\right) l^{r}\left(
d\mathbf{\beta }_{r}\right)  \\
&=&\int_{A\times C}\int_{A\times C}\left[ \int_{A\times D}f\left( \mathbf{%
\xi }_{m},\mathbf{\alpha }_{r}\right) f\left( \mathbf{\gamma }_{m},\mathbf{%
\alpha }_{r}\right) \mathbf{1}_{\left( \overline{\mathbf{\alpha }}_{r}<\min
\left( \overline{\mathbf{\xi }}_{m},\overline{\mathbf{\gamma }}_{m}\right)
\right) }l^{r}\left( d\mathbf{\alpha }_{r}\right) \right] ^{2}l^{m}\left( d%
\mathbf{\xi }_{m}\right) l^{m}\left( d\mathbf{\gamma }_{m}\right)
\end{eqnarray*}%
(the equalities after the definitions of $Q_{3}\left( n\right) $ and $%
Q_{4}\left( n\right) $ are again a consequence of the Fubini theorem)
converges to $0$ as $n\rightarrow +\infty $. This fact will be used in the
proof of Theorem \ref{T : Main}-2.

\bigskip

(\textit{Proof of Theorem \ref{T : Main}-1}) By using a standard
multiplication formula for multiple stochastic integrals (see e.g. \cite[%
Proposition 1.5.1]{Nualart}), we obtain that
\begin{equation*}
F_{n}^{2}=d!\left\Vert f_{n}\right\Vert _{\mathfrak{H}_{n}^{\otimes
d}}^{2}+\sum_{r=1}^{d}\left( d-r\right) !\dbinom{d}{r}^{2}I_{2r}^{X_{n}}%
\left[ f_{n}\otimes _{d-r}f_{n}\right] ,\text{ \ \ }n\geq 1\text{,}
\end{equation*}%
and consequently, for $n\geq 1$,%
\begin{eqnarray}
\mathbb{E}\left[ F_{n}^{2}\mid \mathcal{F}_{t_{n}}^{\pi _{n}}\left(
X_{n}\right) \right] &=&d!\left\Vert f_{n}\right\Vert _{\mathfrak{H}%
_{n}^{\otimes d}}^{2}  \label{Multdev} \\
&&+\sum_{r=1}^{d}\left( d-r\right) !\dbinom{d}{r}^{2}I_{2r}^{X_{n}}\left[
\pi _{n,t_{n}}^{\otimes 2r}\left( f_{n}\otimes _{d-r}f_{n}\right) \right] .
\notag
\end{eqnarray}

Now observe that, for $r=1,...,d$,
\begin{eqnarray}
I_{2r}^{X_{n}}\left[ \pi _{n,t_{n}}^{\otimes 2r}\left( f_{n}\otimes
_{d-r}f_{n}\right) \right] &=&I_{2r}^{X_{n}}\left[ \pi _{n,t_{n}}^{\otimes
2r}\left( f_{n}\otimes _{d-r}^{\pi _{n},t_{n}}f_{n}\right) \right]
\label{M2} \\
&&+I_{2r}^{X_{n}}\left[ \left( \pi _{n,t_{n}}^{\otimes d}f_{n}\right)
\otimes _{d-r}\left( \pi _{n,t_{n}}^{\otimes d}f_{n}\right) \right]  \notag
\end{eqnarray}%
and, in particular, for $r=d$%
\begin{equation}
I_{2d}^{X_{n}}\left[ \pi _{n,t_{n}}^{\otimes 2d}\left( f_{n}\otimes
_{0}f_{n}\right) \right] =I_{2d}^{X_{n}}\left[ \left( \pi
_{n,t_{n}}^{\otimes d}f_{n}\right) \otimes _{0}\left( \pi
_{n,t_{n}}^{\otimes d}f_{n}\right) \right] .  \label{M3}
\end{equation}%
\

It follows from formulae (\ref{Multdev}), (\ref{M2}) and (\ref{M3}), that
Theorem \ref{T : Main}-1 is proved, once it is shown that $I_{2r}^{X_{n}}%
\left[ \left( \pi _{n,t_{n}}^{\otimes d}f_{n}\right) \otimes _{d-r}\left(
\pi _{n,t_{n}}^{\otimes d}f_{n}\right) \right] \rightarrow 0$, in $%
L^{2}\left( \mathbb{P}\right) $, for every $r=1,...,d$. But
\begin{equation*}
\mathbb{E}\left\{ I_{2r}^{X_{n}}\left[ \left( \pi _{n,t_{n}}^{\otimes
d}f_{n}\right) \otimes _{d-r}\left( \pi _{n,t_{n}}^{\otimes d}f_{n}\right) %
\right] ^{2}\right\} \leq \left( d-r\right) !\left\Vert \pi
_{n,t_{n}}^{\otimes d}f_{n}\right\Vert _{\mathfrak{H}^{\otimes d}}^{4}%
\underset{n\rightarrow +\infty }{\rightarrow }0
\end{equation*}%
due to assumption (\ref{NEG}), hence yielding the desired conclusion.

\bigskip

(\textit{Proof of Theorem \ref{T : Main}-2}) For $m\geq 1$, we write $%
\mathbf{x}_{m}$ to indicate a vector $\mathbf{x}_{m}=\left(
x_{1},...,x_{m}\right) \in \left[ 0,1\right] ^{m}$, and also $\widehat{%
\mathbf{x}}_{m}=\max_{i=1,...,n}\left( x_{i}\right) $. Moreover, $d\mathbf{x}%
_{m}$ indicates the restriction of the Lebesgue measure to $\left[ 0,1\right]
^{m}$. We first prove Theorem \ref{T : Main}-2 when the following
assumptions (a) and (b) are verified: (a) for every $n\geq 1$,
\begin{equation}
\mathfrak{H}_{n}=L^{2}\left( A_{n}\times \left[ 0,1\right] ,\mu _{n},%
\mathcal{A}_{n}\otimes \mathcal{B}\left( \left[ 0,1\right] \right) \right)
\text{,}  \label{H1}
\end{equation}%
where $\left( A_{n},\mathcal{A}_{n}\right) $ is a measurable space, $\nu
_{n} $ is a $\sigma $-finite (positive) measure on $\left( A_{n},\mathcal{A}%
_{n}\right) $, and
\begin{equation}
\mu _{n}\left( da,dx\right) =k_{n}\left( a,x\right) \left\{ \nu _{n}\left(
da\right) \otimes dx\right\} \text{,}  \label{H2}
\end{equation}%
where $k_{n}\in L^{1}\left( A_{n}\times \left[ 0,1\right] ,\nu _{n},\mathcal{%
A}_{n}\otimes \mathcal{B}\left( \left[ 0,1\right] \right) \right) $ and $%
k_{n}\geq 0$; (b) for every $n$, for every $\left( a,x\right) \in
A_{n}\times \left[ 0,1\right] $ and for every $h\in \mathfrak{H}_{n}$,%
\begin{equation}
\pi _{n,t}h\left( a,x\right) =h\left( a,x\right) \mathbf{1}_{\left[ 0,t%
\right] }\left( x\right) \text{, \ \ }\forall t\in \left[ 0,1\right] .
\label{projecz}
\end{equation}

Note that $\mu _{n}$ is non-atomic, and also that, in this setting, $%
\mathfrak{H}_{n}^{\otimes d}=L^{2}\left( \left( A_{n}\times \left[ 0,1\right]
\right) ^{d},\mu _{n}^{\otimes d}\right) $ for every $d\geq 2$ and
\begin{equation*}
\mathfrak{H}_{n}^{\odot d}=L_{s}^{2}\left( \left( A_{n}\times \left[ 0,1%
\right] \right) ^{d},\mu _{n}^{\otimes d}\right) .
\end{equation*}%
\

It follows that every $f\in \mathfrak{H}_{n}^{\odot d}$ can be identified
with a (square integrable) function
\begin{equation*}
f\left( a_{1},...,a_{d};x_{1},...,x_{d}\right) =f\left( \mathbf{a}_{d};%
\mathbf{x}_{d}\right) \text{, \ \ }\mathbf{a}_{d}\in A_{n}^{d}\text{, \ }%
\mathbf{x}_{d}\in \left[ 0,1\right] ^{d}\text{,}
\end{equation*}%
which is symmetric in the variables $\left( a_{1},x_{1}\right) ,...,\left(
a_{d},x_{d}\right) $. Moreover, by using the notation introduced in formulae
(\ref{op1})-(\ref{op3}), for every $f\in \mathfrak{H}_{n}^{\otimes d}$ and
every $t\in \left[ 0,1\right] $,%
\begin{eqnarray}
\pi _{n,t}^{\otimes d}f\left( \mathbf{a}_{d};\mathbf{x}_{d}\right)
&=&f\left( \mathbf{a}_{d};\mathbf{x}_{d}\right) \mathbf{1}_{\left[ 0,t\right]
^{d}}\left( \mathbf{x}_{d}\right) \text{ \ \ and}  \label{projectors} \\
\left( \pi _{n,1}^{\otimes d}-\pi _{n,t}^{\otimes d}\right) f\left( \mathbf{a%
}_{d};\mathbf{x}_{d}\right) &=&f\left( \mathbf{a}_{d};\mathbf{x}_{d}\right)
\mathbf{1}_{\left[ 0,1\right] ^{d}\backslash \left[ 0,t\right] ^{d}}\left(
\mathbf{x}_{d}\right) ,\text{ \ \ }\mathbf{a}_{d}\in A_{n}^{d}\text{, \ }%
\mathbf{x}_{d}\in \left[ 0,1\right] ^{d}\text{.}  \label{projectors2}
\end{eqnarray}

Finally, we observe that (by using the notation introduced before the
statement of Lemma \ref{L : Sq. Eq.}), for every $m\geq 1$, every $\mathbf{x}%
_{m}=\left( x_{1},...,x_{m}\right) \in \left[ 0,1\right] ^{m}$ and every $%
\mathbf{\xi }_{m}\in \left( A_{n}\times \left[ 0,1\right] \right) ^{m}$ with
the form $\mathbf{\xi }_{m}=\left( \left( a_{1},x_{1}\right) ,...,\left(
a_{m},x_{m}\right) \right) $,
\begin{equation}
\overline{\mathbf{\xi }}_{m}=\overline{\left( \left( a_{1},x_{1}\right)
,...,\left( a_{m},x_{m}\right) \right) }=\widehat{\left(
x_{1},...,x_{m}\right) }=\widehat{\mathbf{x}}_{m}.  \label{maxmax}
\end{equation}

\bigskip

For $\mathfrak{H}_{n}$ and $\pi _{n}\in \mathcal{R}\left( \mathfrak{H}%
_{n}\right) $ ($n\geq 1$) as in (\ref{H1}), (\ref{H2}) and (\ref{projecz}),
consider the sequence of isonormal Gaussian processes $X_{n}=X_{n}\left(
\mathfrak{H}_{n}\right) $, $n\geq 1$, appearing in the statement of Theorem %
\ref{T : Main}. Since, according to (\ref{MWIintegr}), $F_{n}=I_{d}^{X_{n}}%
\left( f_{n}\right) $, we obtain immediately that, for $n\geq 1$,%
\begin{equation}
D_{X_{n}}F_{n}\left( a,x\right) =d\times I_{d-1}^{X_{n}}\left( f_{n}\left(
a,\cdot ;x,\cdot \right) \right) \text{, }  \label{derive}
\end{equation}%
where, for every fixed $\left( a,x\right) \in A\times \left[ 0,1\right] $, $%
f_{n}\left( a,\cdot ;x,\cdot \right) $ stands for the (symmetric) function,
from $(A$ $\times $ $\left[ 0,1\right] )^{d-1}$ to $\mathbb{R}$,
\begin{equation*}
\left( a_{1},...,a_{d-1};x_{1},...,x_{d-1}\right) \mapsto f_{n}\left(
a,a_{1},...,a_{d-1};x,x_{1},...,x_{d-1}\right) \text{.}
\end{equation*}

In this framework, the sequence $proj\left\{ D_{X_{n}}F_{n}\mid L_{\pi
_{n}}^{2}\left( \mathfrak{H}_{n},X_{n}\right) \right\} $, $n\geq 1$, can be
easily made explicit by means of the following result.

\begin{lemma}
\label{L : simpePD}If $\mathfrak{H}_{n}$ and $\pi _{n}\in \mathcal{R}\left(
\mathfrak{H}_{n}\right) $, $n\geq 1$, satisfy relations (\ref{H1}), (\ref{H2}%
) and (\ref{projecz}), for every $u=u\left( a,x\right) \in L^{2}\left(
\mathfrak{H}_{n},X_{n}\right) $, $\mathbb{P}$-a.s.,
\begin{equation}
proj\left\{ u\mid L_{\pi _{n}}^{2}\left( \mathfrak{H}_{n},X_{n}\right)
\right\} \left( a,x\right) =\mathbb{E}\left[ u\left( a,x\right) \mid
\mathcal{F}_{x}^{\pi _{n}}\left( X_{n}\right) \right] \text{,}
\label{identity}
\end{equation}%
for $\mu _{n}$-a.e. $\left( a,x\right) $, where the filtration $\mathcal{F}%
_{x}^{\pi _{n}}\left( X_{n}\right) $, $x\in \left[ 0,1\right] $, is defined
according to (\ref{resfiltration}).
\end{lemma}

\begin{proof}
Denote by $u^{\ast }$ the process appearing on the right hand side of (\ref%
{identity}). To show that $u^{\ast }$ is an element of $L_{\pi
_{n}}^{2}\left( \mathfrak{H}_{n},X_{n}\right) $ we need to show that it is a
$\pi _{n}$-adapted element of $L^{2}\left( \mathfrak{H}_{n},X_{n}\right) .$
Since $u$ belongs to $L^{2}\left( \mathfrak{H}_{n},X_{n}\right) $, so does $%
u^{\ast }$. Moreover, $u^{\ast }$ is $\pi _{n}$-adapted because, for every $%
h\in \mathfrak{H}_{n}$ and every $t\in \left[ 0,1\right] $,%
\begin{eqnarray*}
\left( u^{\ast },\pi _{n,t}h\right) _{\mathfrak{H}_{n}} &=&\int_{A_{n}\times %
\left[ 0,1\right] }u^{\ast }\left( a,x\right) \pi _{n,t}h\left( a,x\right)
\mu _{n}\left( da,dx\right) \\
&=&\int_{A_{n}\times \left[ 0,t\right] }\mathbb{E}\left[ u\left( a,x\right)
\mid \mathcal{F}_{x}^{\pi _{n}}\left( X_{n}\right) \right] h\left(
a,x\right) \mu _{n}\left( da,dx\right) \in \mathcal{F}_{t}^{\pi _{n}}\left(
X_{n}\right) ,
\end{eqnarray*}%
by (\ref{projecz}). Now consider an element of $\mathcal{E}_{\pi _{n}}\left(
\mathfrak{H}_{n},X_{n}\right) $ with the form $g=\Phi \left( t_{1}\right)
\left( \pi _{n,t_{2}}-\pi _{n,t_{1}}\right) f$ where $t_{2}>t_{1}$, $f\in
\mathfrak{H}_{n}$ and $\Phi \left( t_{1}\right) \in \mathcal{F}_{t_{1}}^{\pi
_{n}}\left( X_{n}\right) $ is square-integrable. Then,
\begin{eqnarray*}
\left( u,g\right) _{L^{2}\left( \mathfrak{H}_{n},X_{n}\right) } &=&\mathbb{E}%
\int_{A_{n}\times \left[ 0,1\right] }u\left( a,t\right) g\left( a,x\right)
\mu _{n}\left( da,dx\right) \\
&=&\int_{A_{n}\times \left( t_{1},t_{2}\right] }\mathbb{E}\left( \Phi \left(
t_{1}\right) u\left( a,x\right) f\left( a,x\right) \right) \mu _{n}\left(
da,dx\right) \\
&=&\int_{A_{n}\times \left( t_{1},t_{2}\right] }\mathbb{E}\left( \Phi \left(
t_{1}\right) \mathbb{E}\left[ u\left( a,x\right) \mid \mathcal{F}_{x}^{\pi
_{n}}\left( X_{n}\right) \right] \right) f\left( a,x\right) \mu _{n}\left(
da,dx\right) \\
&=&\left( u^{\ast },g\right) _{L^{2}\left( \mathfrak{H}_{n},X_{n}\right) }%
\text{,}
\end{eqnarray*}%
where we have used a Fubini theorem and the fact that $\Phi \left(
t_{1}\right) \in \mathcal{F}_{t_{1}}^{\pi _{n}}\left( X_{n}\right) $. Since $%
\mathcal{E}_{\pi _{n}}\left( \mathfrak{H}_{n},X_{n}\right) $ is total in $%
L_{\pi _{n}}^{2}\left( \mathfrak{H}_{n},X_{n}\right) $, we deduce that $%
\left( u,g\right) _{L^{2}\left( \mathfrak{H}_{n},X_{n}\right) }=\left(
u^{\ast },g\right) _{L^{2}\left( \mathfrak{H}_{n},X_{n}\right) }$, for every
$g\in L_{\pi _{n}}^{2}\left( \mathfrak{H}_{n},X_{n}\right) $, hence $u^{\ast
}=proj\left\{ u\mid L_{\pi _{n}}^{2}\left( \mathfrak{H}_{n},X_{n}\right)
\right\} $ as required.
\end{proof}

\bigskip

In particular, thanks to the classic properties of multiple Wiener-It\^{o}
integral and conditional expectations (see e.g. \cite{Nualart2}), we deduce
from (\ref{projectors}) and (\ref{derive}) that, for $F_{n}=I_{d}^{X_{n}}%
\left( f_{n}\right) $ as in (\ref{MWIintegr}) and $x\in \left[ 0,1\right] $,%
\begin{eqnarray}
proj\left\{ D_{X_{n}}F_{n}\mid L_{\pi _{n}}^{2}\left( \mathfrak{H}%
_{n},X_{n}\right) \right\} \left( a,x\right) &=&d\mathbb{E}\left[
I_{d-1}^{X_{n}}\left( f_{n}\left( a,\cdot ;x,\cdot \right) \right) \mid
\mathcal{F}_{x}^{\pi _{n}}\left( X_{n}\right) \right]  \label{derMWI} \\
&=&dI_{d-1}^{X_{n}}\left( f_{n}\left( a,\cdot ;x,\cdot \right) \mathbf{1}%
_{A_{n}^{d-1}\times \left[ 0,x\right] ^{d-1}}\left( \cdot ,\cdot \right)
\right) ,  \notag
\end{eqnarray}%
and consequently
\begin{equation}
\left\Vert proj\left\{ D_{X_{n}}F_{n}\mid L_{\pi _{n}}^{2}\left( \mathfrak{H}%
_{n},X_{n}\right) \right\} \right\Vert _{\mathfrak{H}_{n}}^{2}=d^{2}%
\int_{A_{n}\times \left[ 0,1\right] }I_{d-1}^{X_{n}}\left( f_{n}\left(
a,\cdot ;x,\cdot \right) \mathbf{1}_{A_{n}^{d-1}\times \left[ 0,x\right]
^{d-1}}\right) ^{2}\mu _{n}\left( da,dx\right) .  \label{norm}
\end{equation}

Now note that, thanks to (\ref{isoskorohod}), (\ref{ACO}) and the fact that $%
\mathbb{E}\left( F_{n}\right) =0$,%
\begin{eqnarray*}
\mathbb{E}\left\Vert proj\left\{ D_{X_{n}}F_{n}\mid L_{\pi _{n}}^{2}\left(
\mathfrak{H}_{n},X_{n}\right) \right\} \right\Vert _{\mathfrak{H}_{n}}^{2}
&=&\mathbb{E}\left[ \delta \left( proj\left\{ D_{X_{n}}F_{n}\mid L_{\pi
_{n}}^{2}\left( \mathfrak{H}_{n},X_{n}\right) \right\} \right) ^{2}\right] \\
&=&\mathbb{E}\left[ F_{n}^{2}\right] =d!\left\Vert f_{n}\right\Vert _{%
\mathfrak{H}_{n}^{\otimes d}}^{2}.
\end{eqnarray*}

Moreover, the chaotic expansion of the right hand side of (\ref{norm}) can
be made explicit thanks to the standard multiplication formula (see again
\cite[Proposition 1.5.1]{Nualart})%
\begin{equation*}
I_{d-1}^{X_{n}}\left( g\right)^{2} =\left( d-1\right) !\left\Vert
g\right\Vert _{\mathfrak{H}^{\otimes d-1}}^{2}+\sum_{q=0}^{d-1}q!\dbinom{d-1%
}{q}^{2}I_{2\left( d-1-q\right) }^{X_{n}}\left( g\otimes _{q}g\right) \text{,%
}
\end{equation*}%
applied to $g=f_{n}\left( a,\cdot ;x,\cdot \right) \mathbf{1}%
_{A_{n}^{d-1}\times \left[ 0,x\right] ^{d-1}}$ (for every fixed $\left(
a,x\right) $), from which we obtain%
\begin{eqnarray}
&&d^{2}\int_{A_{n}\times \left[ 0,1\right] }I_{d-1}^{X_{n}}\left(
f_{n}\left( a,\cdot ;x,\cdot \right) \mathbf{1}_{A_{n}^{d-1}\times \left[ 0,x%
\right] ^{d-1}}\right) ^{2}\mu _{n}\left( da,dx\right)  \notag \\
&=&d!\left\Vert f_{n}\right\Vert _{\mathfrak{H}_{n}^{\otimes
d}}^{2}+d^{2}\sum_{q=0}^{d-2}q!\dbinom{d-1}{q}^{2}\times  \notag \\
&&\text{ \ \ }\times \int_{A_{n}\times \left[ 0,1\right] }I_{2\left(
d-1-q\right) }^{X_{n}}\left( \int_{\left( A_{n}\times \left[ 0,x\right]
\right) ^{q}}f_{n}\left( a,\mathbf{a}_{q},\cdot ;x,\mathbf{x}_{q},\cdot
\right) \right. \times  \notag \\
&&\text{ \ \ \ \ \ \ \ \ \ }\times \left. f_{n}\left( a,\mathbf{a}_{q},\cdot
\cdot ;x,\mathbf{x}_{q},\cdot \cdot \right) \mathbf{1}_{\left( A_{n}\times %
\left[ 0,x\right] \right) ^{2\left( d-1-q\right) }}\left( \cdot ;\cdot \cdot
\right) \mu _{n}^{\otimes q}\left( d\mathbf{a}_{q},d\mathbf{x}_{q}\right)
^{^{^{^{^{^{{}}}}}}}\right) \mu _{n}\left( da,dx\right)  \notag \\
&=&d!\left\Vert f_{n}\right\Vert _{\mathfrak{H}_{n}^{\otimes
d}}^{2}+\sum_{r=1}^{d-1}\left( d-r\right) !\dbinom{d}{r}^{2}\left(
d-r\right) \times  \label{fubeq} \\
&&\times I_{2r}^{X_{n}}\left( \int_{A_{n}\times \left[ 0,1\right]
}\int_{\left( A_{n}\times \left[ 0,x\right] \right) ^{d-1-r}}f_{n}\left( a,%
\mathbf{a}_{d-1-r},\cdot ;x,\mathbf{x}_{d-1-r},\cdot \right) \right. \times
\notag \\
&&\times \left. f_{n}\left( a,\mathbf{a}_{d-1-r},\cdot \cdot ;x,\mathbf{x}%
_{d-1-r},\cdot \cdot \right) \mathbf{1}_{\left( A_{n}\times \left[ 0,x\right]
\right) ^{2r}}\left( \cdot ,\cdot \cdot \right) \mu _{n}^{\otimes
d-1-r}\left( d\mathbf{a}_{d-1-r},d\mathbf{x}_{d-1-r}\right) \mu _{n}\left(
da,dx\right) ^{^{^{^{^{^{{}}}}}}}\right)  \notag
\end{eqnarray}%
where the last term is obtained by putting $r=d-q-1$, and by using the
identity
\begin{equation*}
d^{2}q!\dbinom{d-1}{d-q-1}^{2}=d^{2}\left( d-1-r\right) !\dbinom{d-1}{r}%
^{2}=\left( d-r\right) !\dbinom{d}{r}^{2}\left( d-r\right) \text{,}
\end{equation*}%
and where we also applied, to obtain (\ref{fubeq}), a standard stochastic
Fubini theorem (which is a consequence of the linearity of multiple
stochastic integrals--see e.g. \cite[Chapter 1]{Nualart}). We shall now use
the symmetry of the function $f_{n}$ in its first $d-r$ variables, as well
as the relation $\left[ 0,1\right] ^{d-r}\overset{a.s.}{=}\cup
_{i=1,...,d-r}\left\{ \left( a_{1},...,a_{d-r}\right) :a_{i}>a_{j}\text{, \
\ }\forall j\neq i\right\} $, where the union is disjoint, and the symbol ` $%
\overset{a.s.}{=}$ ' means that the equality is true up to sets of zero
Lebesgue measure. Thus, for $r=1,...,d-1$, and for any pair $\left( \mathbf{b%
}_{r},\mathbf{z}_{r}\right) ,\left( \mathbf{b}_{r}^{\prime },\mathbf{z}%
_{r}^{\prime }\right) \in \left( A_{n}\times \left[ 0,1\right] \right) ^{r}$%
\begin{eqnarray*}
&&\left( d-r\right) \int_{A_{n}\times \left[ 0,1\right] }\int_{\left(
A_{n}\times \left[ 0,x\right] \right) ^{d-1-r}}f_{n}\left( a,\mathbf{a}%
_{d-1-r},\mathbf{b}_{r};x,\mathbf{x}_{d-1-r},\mathbf{z}_{r}\right)
f_{n}\left( a,\mathbf{a}_{d-1-r},\mathbf{b}_{r}^{\prime };x,\mathbf{x}%
_{d-1-r},\mathbf{z}_{r}^{\prime }\right) \times \\
&&\times \mathbf{1}_{\left( A_{n}\times \left[ 0,x\right] \right)
^{r}}\left( \mathbf{b}_{r},\mathbf{z}_{r}\right) \mathbf{1}_{\left(
A_{n}\times \left[ 0,x\right] \right) ^{r}}\left( \mathbf{b}_{r}^{\prime },%
\mathbf{z}_{r}^{\prime }\right) \mu _{n}^{\otimes d-1-r}\left( d\mathbf{a}%
_{d-1-r},d\mathbf{x}_{d-1-r}\right) \mu _{n}\left( da,dx\right) \\
&=&\int_{\left( A_{n}\times \left[ 0,1\right] \right) ^{d-r}}f_{n}\left(
\mathbf{a}_{d-r},\mathbf{b}_{r};\mathbf{x}_{d-r},\mathbf{z}_{r}\right)
f_{n}\left( \mathbf{a}_{d-r},\mathbf{b}_{r}^{\prime };\mathbf{x}_{d-r},%
\mathbf{z}_{r}^{\prime }\right) \times \\
&&\times \mathbf{1}_{\left( A_{n}\times \left[ 0,\widehat{\mathbf{x}}_{d-r}%
\right] \right) ^{r}}\left( \mathbf{b}_{r},\mathbf{z}_{r}\right) \mathbf{1}%
_{\left( A_{n}\times \left[ 0,\widehat{\mathbf{x}}_{d-r}\right] \right)
^{r}}\left( \mathbf{b}_{r}^{\prime },\mathbf{z}_{r}^{\prime }\right) \mu
_{n}^{\otimes d-r}\left( d\mathbf{a}_{d-r},d\mathbf{x}_{d-r}\right) \\
&=&\int_{\left( A_{n}\times \left[ 0,1\right] \right) ^{d-r}}f_{n}\left(
\mathbf{a}_{d-r},\mathbf{b}_{r};\mathbf{x}_{d-r},\mathbf{z}_{r}\right)
f_{n}\left( \mathbf{a}_{d-r},\mathbf{b}_{r}^{\prime };\mathbf{x}_{d-r},%
\mathbf{z}_{r}^{\prime }\right) \mathbf{1}_{\left\{ \max \left( \widehat{%
\mathbf{z}}_{r},\widehat{\mathbf{z}}_{r}^{\prime }\right) \leq \widehat{%
\mathbf{x}}_{d-r}\right\} }\mu _{n}^{\otimes d-r}\left( d\mathbf{a}_{d-r},d%
\mathbf{x}_{d-r}\right) \text{;}
\end{eqnarray*}%
(recall that $\widehat{\mathbf{x}}_{d-r}=\max_{i=1,...,d-r}x_{i}$). Observe
that the last integral would be the contraction $f_{n}\otimes _{d-r}f_{n}$,
if there was no indicator functions inside the integral. Now denote by%
\begin{equation*}
\int_{\left( A_{n}\times \left[ 0,1\right] \right) ^{d-r}}f_{n}\left(
\mathbf{a}_{d-r},\ast \mathbf{\ast }_{r};\mathbf{x}_{d-r},\mathbf{\ast }%
_{r}\right) f_{n}\left( \mathbf{a}_{d-r},\ast \mathbf{\ast }_{r}^{\prime };%
\mathbf{x}_{d-r},\mathbf{\ast }_{r}^{\prime }\right) \mathbf{1}_{\left\{
\max \left( \widehat{\mathbf{\ast }}_{r},\widehat{\mathbf{\ast }}%
_{r}^{\prime }\right) \leq \widehat{\mathbf{x}}_{d-r}\right\} }\mu
_{n}^{\otimes d-r}\left( d\mathbf{a}_{d-r},d\mathbf{x}_{d-r}\right)
\end{equation*}%
the application, from $\left( A_{n}\times \left[ 0,x\right] \right)
^{r}\times \left( A_{n}\times \left[ 0,x\right] \right) ^{r}$ to $\mathbb{R}$%
, given by%
\begin{eqnarray*}
\left( \left( \mathbf{b}_{r},\mathbf{z}_{r}\right) ,\left( \mathbf{b}%
_{r}^{\prime },\mathbf{z}_{r}^{\prime }\right) \right) &\mapsto
&\int_{\left( A_{n}\times \left[ 0,1\right] \right) ^{d-r}}f_{n}\left(
\mathbf{a}_{d-r},\mathbf{b}_{r};\mathbf{x}_{d-r},\mathbf{z}_{r}\right) \times
\\
&&\times f_{n}\left( \mathbf{a}_{d-r},\mathbf{b}_{r}^{\prime };\mathbf{x}%
_{d-r},\mathbf{z}_{r}^{\prime }\right) \mathbf{1}_{\left\{ \max \left(
\widehat{\mathbf{z}}_{r},\widehat{\mathbf{z}}_{r}^{\prime }\right) \leq
\widehat{\mathbf{x}}_{d-r}\right\} }\mu _{n}^{\otimes d-r}\left( d\mathbf{a}%
_{d-r},d\mathbf{x}_{d-r}\right) .
\end{eqnarray*}

Relation (\ref{norm}) and the preceding computation imply that%
\begin{eqnarray}
&&\left\Vert proj\left\{ D_{X_{n}}F_{n}\mid L_{\pi ^{\left( n\right)
}}^{2}\left( \mathfrak{H}_{n},X_{n}\right) \right\} \right\Vert _{\mathfrak{H%
}_{n}}^{2}  \notag \\
&=&d!\left\Vert f_{n}\right\Vert _{\mathfrak{H}_{n}^{\odot
d}}^{2}+\sum_{r=1}^{d-1}\left( d-r\right) !\dbinom{d}{r}^{2}I_{2r}^{X_{n}}%
\left( \int_{\left( A_{n}\times \left[ 0,1\right] \right) ^{d-r}}f_{n}\left(
\mathbf{a}_{d-r},\ast \mathbf{\ast }_{r};x,\mathbf{x}_{d-r},\mathbf{\ast }%
_{r}\right) \right. \times  \label{a} \\
&&\times \left. f_{n}\left( \mathbf{a}_{d-r},\ast \mathbf{\ast }_{r}^{\prime
};\mathbf{x}_{d-r},\mathbf{\ast }_{r}^{\prime }\right) \mathbf{1}_{\left\{
\max \left( \widehat{\mathbf{\ast }}_{r},\widehat{\mathbf{\ast }}%
_{r}^{\prime }\right) \leq \widehat{\mathbf{x}}_{d-r}\right\} }\mu
_{n}^{\otimes d-r}\left( d\mathbf{a}_{d-r},d\mathbf{x}_{d-r}\right)
^{^{^{^{{}}}}}\right) .  \notag
\end{eqnarray}

Now, for $r=1,...,d-1$ and $t\in \left[ 0,1\right] $,%
\begin{eqnarray*}
&&\int_{\left( A_{n}\times \left[ 0,1\right] \right) ^{d-r}}f_{n}\left(
\mathbf{a}_{d-r},\ast \mathbf{\ast }_{r};\mathbf{x}_{d-r},\mathbf{\ast }%
_{r}\right) f_{n}\left( \mathbf{a}_{d-r},\ast \mathbf{\ast }_{r}^{\prime };%
\mathbf{x}_{d-r},\mathbf{\ast }_{r}^{\prime }\right) \mathbf{1}_{\left\{
\max \left( \widehat{\mathbf{\ast }}_{r},\widehat{\mathbf{\ast }}%
_{r}^{\prime }\right) \leq \widehat{\mathbf{x}}_{d-r}\right\} }\mu
_{n}^{\otimes d-r}\left( d\mathbf{a}_{d-r},d\mathbf{x}_{d-r}\right) \\
&=&\int_{\left( A_{n}\times \left[ 0,t\right] \right) ^{d-r}}f_{n}\left(
\mathbf{a}_{d-r},\ast \mathbf{\ast }_{r};\mathbf{x}_{d-r},\mathbf{\ast }%
_{r}\right) f_{n}\left( \mathbf{a}_{d-r},\ast \mathbf{\ast }_{r}^{\prime };%
\mathbf{x}_{d-r},\mathbf{\ast }_{r}^{\prime }\right) \mathbf{1}_{\left\{
\max \left( \widehat{\mathbf{\ast }}_{r},\widehat{\mathbf{\ast }}%
_{r}^{\prime }\right) \leq \widehat{\mathbf{x}}_{d-r}\right\} }\mu
_{n}^{\otimes d-r}\left( d\mathbf{a}_{d-r},d\mathbf{x}_{d-r}\right) \\
&&+\int_{A_{n}^{d-r}\times \left( \left[ 0,1\right] ^{d-r}\backslash \left[
0,t\right] ^{d-r}\right) }f_{n}\left( \mathbf{a}_{d-r},\ast \mathbf{\ast }%
_{r};\mathbf{x}_{d-r},\mathbf{\ast }_{r}\right) \times \\
&&\text{ \ \ \ \ \ \ \ \ \ \ \ \ \ \ \ \ \ \ \ }\times f_{n}\left( \mathbf{a}%
_{d-r},\ast \mathbf{\ast }_{r}^{\prime };\mathbf{x}_{d-r},\mathbf{\ast }%
_{r}^{\prime }\right) \mathbf{1}_{\left\{ \max \left( \widehat{\mathbf{\ast }%
}_{r},\widehat{\mathbf{\ast }}_{r}^{\prime }\right) \leq \widehat{\mathbf{x}}%
_{d-r}\right\} }\mu _{n}^{\otimes d-r}\left( d\mathbf{a}_{d-r},d\mathbf{x}%
_{d-r}\right) \\
&\triangleq &G_{n,t}^{r}\left( 1\right) +G_{n,t}^{r}\left( 2\right) ,
\end{eqnarray*}%
(plainly, $G_{n,t}^{r}\left( 1\right) ,G_{n,t}^{r}\left( 2\right) \in
\mathfrak{H}_{n}^{\otimes 2r}$) and observe that, by bounding the indicator
function by $1$ and using the Cauchy-Schwarz inequality,%
\begin{equation*}
\mathbb{E}\left[ I_{2r}^{X_{n}}\left( G_{n,t}^{r}\left( 1\right) \right) ^{2}%
\right] =\left( 2r\right) !\left\Vert \left( G_{n,t}^{r}\left( 1\right)
\right) _{s}\right\Vert _{\mathfrak{H}_{n}^{\otimes 2r}}^{2}\leq \left(
2r\right) !\left\Vert f_{n}\mathbf{1}_{\left[ 0,t\right] ^{d}}\right\Vert _{%
\mathfrak{H}_{n}^{\otimes d}}^{4}.
\end{equation*}

Now, if $t_{n}$ is the sequence in the statement of Theorem \ref{T : Main},
one has
\begin{equation*}
d!\left\Vert f_{n}\mathbf{1}_{\left[ 0,t_{n}\right] ^{d}}\right\Vert _{%
\mathfrak{H}_{n}^{\otimes d}}^{2}=\mathbb{E}\left( I_{d}^{X_{n}}\left( \pi
_{n,t_{n}}^{\otimes d}f_{n}\right) ^{2}\right) \rightarrow 0.
\end{equation*}

Thus, (\ref{NEG}) implies%
\begin{equation*}
\lim_{n\rightarrow +\infty }\mathbb{E}\left[ I_{2r}^{X_{n}}\left(
G_{n,t}^{r}\left( 1\right) \right) ^{2}\right] \leq \left( 2r\right)
!\lim_{n\rightarrow +\infty }\left\Vert f_{n}\mathbf{1}_{\left[ 0,t_{n}%
\right] ^{d}}\right\Vert _{\mathfrak{H}_{n}^{\otimes d}}^{4}=0.
\end{equation*}

We now deal with $G_{n,t}^{r}\left( 2\right) $. For every $t\in \left[ 0,1%
\right] $, we may write%
\begin{equation}
G_{n,t}^{r}\left( 2\right) =G_{n,t}^{r}\left( 2\right) \mathbf{1}%
_{A_{n}^{2r}\times \left[ 0,t\right] ^{2r}}+G_{n,t}^{r}\left( 2\right)
\mathbf{1}_{A_{n}^{2r}\times \left( \left[ 0,1\right] ^{2r}\backslash \left[
0,t\right] ^{2r}\right) }\triangleq H_{n,t}^{r}\left( 1\right)
+H_{n,t}^{r}\left( 2\right) .  \label{GiAcca}
\end{equation}

Consider first $H_{n,t}^{r}\left( 1\right) $. Because of the presence of the
indicator function $\mathbf{1}_{A_{n}^{2r}\times \left[ 0,t\right] ^{2r}}$,
the indicator function in the integral defining $G_{n,t}^{r}\left( 2\right) $
is always equal to 1, and one gets, for every $r=1,...,d-1$,
\begin{eqnarray*}
H_{n,t}^{r}\left( 1\right) &=&\left\{ \int_{A_{n}^{d-r}\times \left( \left[
0,1\right] ^{d-r}\backslash \left[ 0,t\right] ^{d-r}\right) }f_{n}\left(
\mathbf{a}_{d-r},\ast \mathbf{\ast }_{r};\mathbf{x}_{d-r},\mathbf{\ast }%
_{r}\right) \times \right. \\
&&\times \left. f_{n}\left( \mathbf{a}_{d-r},\ast \mathbf{\ast }_{r}^{\prime
};\mathbf{x}_{d-r},\mathbf{\ast }_{r}^{\prime }\right) \mu _{n}^{\otimes
d-r}\left( d\mathbf{a}_{d-r},d\mathbf{x}_{d-r}\right) ^{^{^{^{{}}}}}\right\}
\times \mathbf{1}_{A_{n}^{2r}\times \left[ 0,t\right] ^{2r}} \\
&=&\pi _{n,t_{n}}^{\otimes 2r}\left( f_{n}\otimes _{d-r}^{\pi
_{n},t_{n}}f_{n}\right) ,
\end{eqnarray*}%
which appears in (\ref{DEV}). Then, in view of (\ref{NORM}), (\ref{DEV}) and
(\ref{a}), we have that%
\begin{equation*}
\left\Vert proj\left\{ D_{X_{n}}F_{n}\mid L_{\pi _{n}}^{2}\left( \mathfrak{H}%
_{n},X_{n}\right) \right\} \right\Vert _{\mathfrak{H}_{n}}^{2}=Y+o_{\mathbb{P%
}}\left( 1\right) +\sum_{r=1}^{d-1}\left( d-r\right) !\dbinom{d}{r}%
^{2}I_{2r}^{X_{n}}\left( H_{n,t_{n}}^{r}\left( 2\right) \right) ,
\end{equation*}%
where $H_{n,t_{n}}^{r}\left( 2\right) \in \mathfrak{H}_{n}^{\otimes 2r}$ is
defined by (\ref{GiAcca}). We shall now show that (\ref{asCV}) implies $%
H_{n,t_{n}}^{r}\left( 2\right) \rightarrow 0$ in $\mathfrak{H}_{n}^{\otimes
2r}$, for every $r=1,...,d-1$. Now observe that, because of (\ref{H1}), (\ref%
{H2}) and (\ref{projecz}), condition (\ref{asCV}) can be rewritten as
follows: for every $r=1,...,d-1$, the sequence $Z_{r}\left( n\right) \in
\mathfrak{H}_{n}^{\otimes 2r}$, $n\geq 1$, defined as%
\begin{eqnarray}
Z_{r}\left( n\right) &=&\left( \int_{A_{n}^{d-r}\times \left( \left[ 0,1%
\right] ^{d-r}\backslash \left[ 0,t_{n}\right] ^{d-r}\right) }f_{n}\left(
\mathbf{a}_{d-r},\ast \mathbf{\ast }_{r};\mathbf{x}_{d-r},\mathbf{\ast }%
_{r}\right) \right. \times  \label{Zerren} \\
&&\times \left. f_{n}\left( \mathbf{a}_{d-r},\ast \mathbf{\ast }_{r}^{\prime
};\mathbf{x}_{d-r},\mathbf{\ast }_{r}^{\prime }\right) \mu _{n}^{\otimes
d-r}\left( d\mathbf{a}_{d-r},d\mathbf{x}_{d-r}\right)
^{^{^{^{^{^{{}}}}}}}\right) \mathbf{1}_{A_{n}^{2r}\times \left( \left[ 0,1%
\right] ^{2r}\backslash \left[ 0,t_{n}\right] ^{2r}\right) }.  \notag
\end{eqnarray}%
is such that
\begin{equation}
\lim_{n\rightarrow +\infty }\left\Vert Z_{r}\left( n\right) \right\Vert _{%
\mathfrak{H}_{n}^{\otimes 2r}}^{2}=0.  \label{consequenceZ}
\end{equation}

As a consequence, in this case the statement is proved once it is shown
that, for $r=1,...,d-1$, (\ref{consequenceZ}) implies necessarily that $%
\lim_{n\rightarrow +\infty }\left\Vert H_{n,t_{n}}^{r}\left( 2\right)
\right\Vert _{\mathfrak{H}_{n}^{\otimes 2r}}^{2}=0$ (recall that $%
H_{n,t_{n}}^{r}\left( 2\right) $ is given by (\ref{GiAcca})). To this end,
introduce the notation: for every $q\geq 1$, every $p=0,...,q$ and $t\in %
\left[ 0,1\right] $%
\begin{equation*}
T_{n}^{t}\left( q,p\right) \triangleq \left\{ \left( \mathbf{a}_{q},\mathbf{x%
}_{q}\right) \in A_{n}^{q}\times \left[ 0,1\right] ^{q}:\text{ there are
exactly }p\text{ indices }i\text{ such that }x_{i}\leq t\right\}
\end{equation*}%
and note that, for $q\geq 1$, up to sets of zero $\mu _{n}$ -- measure,
\begin{eqnarray}
S_{n}^{t}\left( q\right) &\triangleq &A_{n}^{q}\times \left( \left[ 0,1%
\right] ^{q}\backslash \left[ 0,t\right] ^{q}\right)
=\bigcup\limits_{p=0}^{q-1}T_{n}^{t}\left( q,p\right) \text{,}  \label{dec1}
\\
S_{n}^{t}\left( 2q\right) &=&A_{n}^{2q}\times \left( \left[ 0,1\right]
^{2q}\backslash \left[ 0,t\right] ^{2q}\right) =\bigcup\limits_{\substack{ %
p,s\geq 0  \\ p\wedge s\leq q-1}}T_{n}^{t}\left( q,p\right) \times
T_{n}^{t}\left( q,s\right) ,  \label{dec2}
\end{eqnarray}%
where the unions are disjoint. With this notation, by (\ref{dec1}), (\ref%
{dec2}) and the fact that $\mu _{n}$ is non-atomic (so that we can write $<$
instead of $\leq $ in the indicator function defining $H_{n,t_{n}}^{r}\left(
2\right) $), we therefore obtain that, for each $n$,%
\begin{eqnarray}
&&\left\Vert H_{n,t_{n}}^{r}\left( 2\right) \right\Vert _{\mathfrak{H}%
_{n}^{\otimes 2r}}^{2}  \notag \\
&=&\int_{A_{n}^{2r}\times \left( \left[ 0,1\right] ^{2r}\backslash \left[
0,t_{n}\right] ^{2r}\right) }\left( \int_{A_{n}^{d-r}\times \left( \left[ 0,1%
\right] ^{d-r}\backslash \left[ 0,t_{n}\right] ^{d-r}\right) }f_{n}\left(
\mathbf{a}_{d-r},\mathbf{b}_{r};\mathbf{x}_{d-r},\mathbf{z}_{r}\right)
f_{n}\left( \mathbf{a}_{d-r},\mathbf{b}_{r}^{\prime };\mathbf{x}_{d-r},%
\mathbf{z}_{r}^{\prime }\right) \right. \times  \notag \\
&&\text{ \ \ \ \ \ \ \ \ \ \ }\times \left. \mathbf{1}_{\left\{ \max \left(
\widehat{\mathbf{z}}_{r},\widehat{\mathbf{z}}_{r}^{\prime }\right) <\widehat{%
\mathbf{x}}_{d-r}\right\} }\mu _{n}^{\otimes d-r}\left( d\mathbf{a}_{d-r},d%
\mathbf{x}_{d-r}\right) ^{^{^{^{^{{}}}}}}\right) ^{2}\mu _{n}^{\otimes
2r}\left( d\mathbf{b}_{r},d\mathbf{z}_{r}\right) \mu _{n}^{\otimes 2r}\left(
d\mathbf{b}_{r}^{\prime },d\mathbf{z}_{r}^{\prime }\right)  \label{h2} \\
&=&\sum_{_{p,s\geq 0\text{, \ }p\wedge s\leq r-1}}\int_{T_{n}^{t_{n}}\left(
r,p\right) \times T_{n}^{t_{n}}\left( r,s\right) }\left(
\sum_{q=0}^{d-r-1}\int_{T_{n}^{t_{n}}\left( d-r,q\right) }f_{n}\left(
\mathbf{a}_{d-r},\mathbf{b}_{r};\mathbf{x}_{d-r},\mathbf{z}_{r}\right)
\right. \times  \label{h21} \\
&&\times \left. f_{n}\left( \mathbf{a}_{d-r},\mathbf{b}_{r}^{\prime };%
\mathbf{x}_{d-r},\mathbf{z}_{r}^{\prime }\right) \mathbf{1}_{\left\{ \max
\left( \widehat{\mathbf{z}}_{r},\widehat{\mathbf{z}}_{r}^{\prime }\right) <%
\widehat{\mathbf{x}}_{d-r}\right\} }\mu _{n}^{\otimes d-r}\left( d\mathbf{a}%
_{d-r},d\mathbf{x}_{d-r}\right) ^{^{^{^{^{{}}}}}}\right) ^{2}\mu
_{n}^{\otimes 2r}\left( d\mathbf{b}_{r},d\mathbf{z}_{r}\right) \mu
_{n}^{\otimes 2r}\left( d\mathbf{b}_{r}^{\prime },d\mathbf{z}_{r}^{\prime
}\right)  \notag
\end{eqnarray}%
To prove that $\left\Vert H_{n,t_{n}}^{r}\left( 2\right) \right\Vert _{%
\mathfrak{H}_{n}^{\otimes 2r}}^{2}\rightarrow 0$, it is therefore sufficient
to show that, for every $r=1,...,d-1$, $q=0,...,d-r-1$ and every $p,s\geq 0$
with $p\wedge s\leq r-1$, the sequence%
\begin{eqnarray}
&&\int_{T_{n}^{t_{n}}\left( r,p\right) }\int_{T_{n}^{t_{n}}\left( r,s\right)
}\left( \int_{T_{n}^{t_{n}}\left( d-r,q\right) }f_{n}\left( \mathbf{a}_{d-r},%
\mathbf{b}_{r};\mathbf{x}_{d-r},\mathbf{z}_{r}\right) f_{n}\left( \mathbf{a}%
_{d-r},\mathbf{b}_{r}^{\prime };\mathbf{x}_{d-r},\mathbf{z}_{r}^{\prime
}\right) \right.  \label{p1} \\
&&\text{ \ \ \ \ \ \ \ \ \ \ }\times \left. \mathbf{1}_{\left\{ \max \left(
\widehat{\mathbf{z}}_{r},\widehat{\mathbf{z}}_{r}^{\prime }\right) <\widehat{%
\mathbf{x}}_{d-r}\right\} }\mu _{n}^{\otimes d-r}\left( d\mathbf{a}_{d-r},d%
\mathbf{x}_{d-r}\right) ^{^{^{^{^{{}}}}}}\right) ^{2}\mu _{n}^{\otimes
2r}\left( d\mathbf{b}_{r},d\mathbf{z}_{r}\right) \mu _{n}^{\otimes 2r}\left(
d\mathbf{b}_{r}^{\prime },d\mathbf{z}_{r}^{\prime }\right)  \notag
\end{eqnarray}%
converges to zero, as $n\rightarrow +\infty $. To prove this result, write
(for $n\geq 1$ and $r=1,...,d-1$) $\left\Vert Z_{r}\left( n\right)
\right\Vert _{\mathfrak{H}_{n}^{\otimes 2r}}^{2}$ by means of (\ref{Zerren}%
), decompose the set $S_{n}^{t_{n}}\left( 2r\right) $ according to (\ref%
{dec2}), and apply a standard Fubini argument to obtain that (\ref%
{consequenceZ}) implies that, for every $r=1,...,d-1$ and $q=0,...,r-1$, the
quantity%
\begin{eqnarray}
&&\int_{S_{n}^{t_{n}}\left( d-r\right) }\int_{S_{n}^{t_{n}}\left( d-r\right)
}\left[ \int_{T_{n}^{t_{n}}\left( r,q\right) }f_{n}\left( \mathbf{a}_{r},%
\mathbf{b}_{d-r};\mathbf{x}_{r},\mathbf{z}_{d-r}\right) \right.  \label{EL}
\\
&&\text{ \ \ \ }\times \left. f_{n}\left( \mathbf{a}_{r},\mathbf{b}%
_{d-r}^{\prime };\mathbf{x}_{r},\mathbf{z}_{d-r}^{\prime }\right) \mu
_{n}^{\otimes r}\left( d\mathbf{a}_{r},d\mathbf{x}_{r}\right) ^{^{^{{}}}}%
\right] ^{2}\mu _{n}^{\otimes d-r}\left( d\mathbf{b}_{d-r},d\mathbf{z}%
_{d-r}\right) \mu _{n}^{\otimes d-r}\left( d\mathbf{b}_{d-r}^{\prime },d%
\mathbf{z}_{d-r}^{\prime }\right)  \notag \\
&=&\sum_{p=0}^{d-r-1}\sum_{s=0}^{d-r-1}\int_{T_{n}^{t_{n}}\left(
d-r,p\right) }\int_{T_{n}^{t_{n}}\left( d-r,s\right) }\left[
\int_{T_{n}^{t_{n}}\left( r,q\right) }f_{n}\left( \mathbf{a}_{r},\mathbf{b}%
_{d-r};\mathbf{x}_{r},\mathbf{z}_{d-r}\right) \right.  \label{EL1} \\
&&\text{ \ \ \ }\times \left. f_{n}\left( \mathbf{a}_{r},\mathbf{b}%
_{d-r}^{\prime };\mathbf{x}_{r},\mathbf{z}_{d-r}^{\prime }\right) \mu
_{n}^{\otimes r}\left( d\mathbf{a}_{r},d\mathbf{x}_{r}\right) ^{^{^{{}}}}%
\right] ^{2}\mu _{n}^{\otimes d-r}\left( d\mathbf{b}_{d-r},d\mathbf{z}%
_{d-r}\right) \mu _{n}^{\otimes d-r}\left( d\mathbf{b}_{d-r}^{\prime },d%
\mathbf{z}_{d-r}^{\prime }\right)  \notag \\
&\rightarrow &0\text{, \ \ as }n\rightarrow +\infty \text{,}  \label{CVzero}
\end{eqnarray}%
where the equality (\ref{EL1}) is a consequence of (\ref{dec2}). Now fix $%
p\in \left\{ 0,...,d-r-1\right\} $ and $q\in \left\{ 0,...,r-1\right\} $.
For every $n\geq 1$, we can apply Lemma \ref{L : Sq. Eq.} in the case $%
A=A_{n}$, $l=\mu _{n}$, $f_{n}=f$, $m=d-r$ and $C\subseteq \left[ 0,1\right]
^{d-r}$ and $D\subseteq \left[ 0,1\right] ^{r}$ such that $A_{n}\times
C=T_{n}^{t_{n}}\left( d-r,p\right) $ and $A_{n}\times D=T_{n}^{t_{n}}\left(
r,q\right) $ to obtain that
\begin{eqnarray}
&&\int_{T_{n}^{t_{n}}\left( d-r,p\right) }\int_{T_{n}^{t_{n}}\left(
d-r,p\right) }\left[ \int_{T_{n}^{t_{n}}\left( r,q\right) }f_{n}\left(
\mathbf{a}_{r},\mathbf{b}_{d-r};\mathbf{x}_{r},\mathbf{z}_{d-r}\right)
\right.  \label{acc1} \\
&&\text{ \ \ \ }\times \left. f_{n}\left( \mathbf{a}_{r},\mathbf{b}%
_{d-r}^{\prime };\mathbf{x}_{r},\mathbf{z}_{d-r}^{\prime }\right) \mu
_{n}^{\otimes r}\left( d\mathbf{a}_{r},d\mathbf{x}_{r}\right) ^{^{^{{}}}}%
\right] ^{2}\mu _{n}^{\otimes d-r}\left( d\mathbf{b}_{d-r},d\mathbf{z}%
_{d-r}\right) \mu _{n}^{\otimes d-r}\left( d\mathbf{b}_{d-r}^{\prime },d%
\mathbf{z}_{d-r}^{\prime }\right)  \notag \\
&=&\int_{T_{n}^{t_{n}}\left( d-r,p\right) }\int_{T_{n}^{t_{n}}\left(
d-r,p\right) }\left[ \int_{T_{n}^{t_{n}}\left( r,q\right) }f_{n}\left(
\mathbf{a}_{r},\mathbf{b}_{d-r};\mathbf{x}_{r},\mathbf{z}_{d-r}\right)
f_{n}\left( \mathbf{a}_{r},\mathbf{b}_{d-r}^{\prime };\mathbf{x}_{r},\mathbf{%
z}_{d-r}^{\prime }\right) \right.  \label{acc2} \\
&&\times \left. \mathbf{1}_{\left( \widehat{\mathbf{x}}_{r}<\max \left(
\widehat{\mathbf{z}}_{d-r},\widehat{\mathbf{z}}_{d-r}^{\prime }\right)
\right) }\mu _{n}^{\otimes r}\left( d\mathbf{a}_{r},d\mathbf{x}_{r}\right)
^{^{^{{}}}}\right] ^{2}\mu _{n}^{\otimes d-r}\left( d\mathbf{b}_{d-r},d%
\mathbf{z}_{d-r}\right) \mu _{n}^{\otimes d-r}\left( d\mathbf{b}%
_{d-r}^{\prime },d\mathbf{z}_{d-r}^{\prime }\right)  \notag \\
&&+\int_{T_{n}^{t_{n}}\left( r,q\right) }\int_{T_{n}^{t_{n}}\left(
r,q\right) }\left[ \int_{T_{n}^{t_{n}}\left( d-r,p\right) }f_{n}\left(
\mathbf{a}_{r},\mathbf{b}_{d-r};\mathbf{x}_{r},\mathbf{z}_{d-r}\right)
f_{n}\left( \mathbf{a}_{r}^{\prime },\mathbf{b}_{d-r};\mathbf{x}_{r}^{\prime
},\mathbf{z}_{d-r}\right) \right.  \notag \\
&&\times \left. \mathbf{1}_{\left( \widehat{\mathbf{z}}_{d-r}<\max \left(
\widehat{\mathbf{x}}_{r},\widehat{\mathbf{x}}_{r}^{\prime }\right) \right)
}\mu _{n}^{\otimes r}\left( d\mathbf{b}_{d-r},d\mathbf{z}_{d-r}\right)
^{^{^{{}}}}\right] ^{2}\mu _{n}^{\otimes r}\left( d\mathbf{a}_{r},d\mathbf{x}%
_{r}\right) \mu _{n}^{\otimes r}\left( d\mathbf{a}_{r}^{\prime },d\mathbf{x}%
_{r}^{\prime }\right) .  \notag
\end{eqnarray}

Now observe that the sequence $\left\{ \left\Vert f_{n}\right\Vert _{%
\mathfrak{H}_{n}^{\otimes d}}:n\geq 1\right\} $ is bounded by assumption
(indeed, relation (\ref{NORM}) holds). We can therefore argue as in the
remark following the proof of Lemma \ref{L : Sq. Eq.}, and deduce that,
since for every $p\in \left\{ 0,...,d-r-1\right\} $ and $q\in \left\{
0,...,r-1\right\} $ the sequence in (\ref{acc1}) converges to $0$ (by (\ref%
{CVzero})), then the sequence in (\ref{p1}) converges to $0$ for every $%
r=1,...,d-1$, whenever $p=s\leq r-1$ and $q=0,...,d-r-1$. To prove that (\ref%
{p1}) converges to $0$ for every $r=1,...,d-1$, $q=0,...,d-r-1$ and every $%
p,s\geq 0$ such that $p\wedge s\leq r-1$, thus concluding the proof of
Theorem \ref{T : Gen CV} in this special setting, observe that, due to the
Fubini theorem, the quantity%
\begin{eqnarray*}
&&\int_{T_{n}^{t_{n}}\left( r,p\right) }\int_{T_{n}^{t_{n}}\left( r,s\right)
}\left( \int_{T_{n}^{t_{n}}\left( d-r,q\right) }f_{n}\left( \mathbf{a}_{d-r},%
\mathbf{b}_{r};\mathbf{x}_{d-r},\mathbf{z}_{r}\right) f_{n}\left( \mathbf{a}%
_{d-r},\mathbf{b}_{r}^{\prime };\mathbf{x}_{d-r},\mathbf{z}_{r}^{\prime
}\right) \right. \\
&&\text{ \ \ \ \ \ \ \ \ \ \ }\times \left. \mathbf{1}_{\left\{ \max \left(
\widehat{\mathbf{z}}_{r},\widehat{\mathbf{z}}_{r}^{\prime }\right) <\widehat{%
\mathbf{x}}_{d-r}\right\} }\mu _{n}^{\otimes d-r}\left( d\mathbf{a}_{d-r},d%
\mathbf{x}_{d-r}\right) ^{^{^{^{^{{}}}}}}\right) ^{2}\mu _{n}^{\otimes
2r}\left( d\mathbf{b}_{r},d\mathbf{z}_{r}\right) \mu _{n}^{\otimes 2r}\left(
d\mathbf{b}_{r}^{\prime },d\mathbf{z}_{r}^{\prime }\right)
\end{eqnarray*}%
can be rewritten as
\begin{eqnarray*}
&&\int_{T_{n}^{t_{n}}\left( d-r,q\right) }\int_{T_{n}^{t_{n}}\left(
d-r,q\right) }\left( \int_{T_{n}^{t_{n}}\left( r,p\right) }f_{n}\left(
\mathbf{a}_{d-r},\mathbf{b}_{r};\mathbf{x}_{d-r},\mathbf{z}_{r}\right)
f_{n}\left( \mathbf{a}_{d-r}^{\prime },\mathbf{b};\mathbf{x}_{d-r}^{\prime },%
\mathbf{z}_{r}\right) \right. \\
&&\text{ \ \ \ \ \ \ \ \ \ \ }\times \left. \mathbf{1}_{\left\{ \widehat{%
\mathbf{z}}_{r}<\min \left( \widehat{\mathbf{x}}_{d-r},\widehat{\mathbf{x}}%
_{d-r}^{\prime }\right) \right\} }\mu _{n}^{\otimes r}\left( d\mathbf{b}%
_{r},d\mathbf{z}_{r}\right) ^{^{^{^{^{{}}}}}}\right) \times \\
&&\left( \int_{T_{n}^{t_{n}}\left( r,s\right) }f_{n}\left( \mathbf{a}_{d-r},%
\mathbf{b}_{r};\mathbf{x}_{d-r},\mathbf{z}_{r}\right) f_{n}\left( \mathbf{a}%
_{d-r}^{\prime },\mathbf{b};\mathbf{x}_{d-r}^{\prime },\mathbf{z}\right)
\right. \times \\
&&\text{ \ \ \ \ \ \ \ \ \ \ \ \ \ }\times \left. \mathbf{1}_{\left\{
\widehat{\mathbf{z}}_{r}<\min \left( \widehat{\mathbf{x}}_{d-r},\widehat{%
\mathbf{x}}_{d-r}^{\prime }\right) \right\} }\mu _{n}^{\otimes r}\left( d%
\mathbf{b}_{r},d\mathbf{z}_{r}\right) ^{^{^{^{^{{}}}}}}\right) \mu
_{n}^{\otimes d-r}\left( d\mathbf{a}_{d-r},d\mathbf{x}_{d-r}\right) \mu
_{n}^{\otimes d-r}\left( d\mathbf{a}_{d-r}^{\prime },d\mathbf{x}%
_{d-r}^{\prime }\right) ,
\end{eqnarray*}%
so that the conclusion is obtained by a further application of (\ref{acc2}),
as well as a standard version of the Cauchy-Schwarz inequality.

\bigskip

To prove Theorem \ref{T : Main} in the general case, we start by showing
that, for every real separable Hilbert space $\mathfrak{H}$, and for every
absolutely continuous resolution of the identity $\pi =\left\{ \pi _{t}:t\in %
\left[ 0,1\right] \right\} \in \mathcal{R}_{AC}\left( \mathfrak{H}\right) $,
there exists a Hilbert space $\mathfrak{H}^{\sharp }$ with the form (\ref{H1}%
) and (\ref{H2}) (the dependence on $n$ has been momentarily dropped), and a
resolution $\pi ^{\sharp }=\left\{ \pi _{t}^{\sharp }:t\in \left[ 0,1\right]
\right\} $ on $\mathfrak{H}^{\sharp }$ as in (\ref{projecz}), such that the
following property is verified: there exists a unitary transformation%
\begin{equation}
T:\mathfrak{H}^{\sharp }\mapsto \mathfrak{H},  \label{ISO}
\end{equation}%
from $\mathfrak{H}^{\sharp }\ $onto $\mathfrak{H}$, such that, for every $%
t\in \left[ 0,1\right] $,%
\begin{equation}
\pi _{t}T=T\pi _{t}^{\sharp }\text{.}  \label{EQUI}
\end{equation}%
In the language of \cite[Definition 5.1]{UZ}, (\ref{EQUI}) implies that the
two pairs $\left( \mathfrak{H},\pi \right) $ and $\left( \mathfrak{H}%
^{\sharp },\pi ^{\sharp }\right) $ are \textit{equivalent}. Note that (\ref%
{EQUI}) holds if, and only if, the following condition is verified: for
every $t\in \left[ 0,1\right] $, $\pi _{t}=T\pi _{t}^{\sharp }T^{-1}$.
Moreover, since $T$ is a unitary transformaton, $T^{-1}=T^{\ast }$. To prove
the existence of such a $T$, let $\pi $ be absolutely continuous on $%
\mathfrak{H}$, set $q=rank\left( \pi \right) $, and consider a fully
orthogonal $\pi $-reproducing subset
\begin{equation}
S=\left\{ g_{j}:1\leq j\leq q\right\} \subseteq \mathfrak{H}.  \label{FO}
\end{equation}

Note that the full orthogonality of $S$ implies that, for every $s,t\in %
\left[ 0,1\right] $ and for every $i\neq j$, $\left( \pi _{s}g_{i},\pi
_{t}g_{j}\right) _{\mathfrak{H}}=0$. Moreover, since $\pi $ is absolutely
continuous, for every $j\geq 1$ there exists a function $m_{j}\left(
t\right) $, $t\in \left[ 0,1\right] $, such that $m_{j}\left( \cdot \right)
\geq 0$, and
\begin{equation}
\left\Vert \pi _{t}g_{j}\right\Vert _{\mathfrak{H}}^{2}=\int_{0}^{t}m_{j}%
\left( x\right) dx\text{, \ \ }t\in \left[ 0,1\right] \text{.}  \label{Acont}
\end{equation}

Note that (\ref{Acont}) implies that $m_{j}\left( \cdot \right) \in
L^{1}\left( \left[ 0,1\right] ,dx\right) $, and that we can always define
the set $S$ in (\ref{FO}) to be such that
\begin{equation}
\sum_{j=1}^{q}\left\Vert g_{j}\right\Vert _{\mathfrak{H}}^{2}=\sum_{j=1}^{q}%
\int_{0}^{1}m_{j}\left( x\right) dx<+\infty .  \label{AcontConv}
\end{equation}%
Now define $A=\mathbb{N}=\left\{ 1,2,...\right\} $, set $\nu $ equal to the
counting measure on $A$, and $\mu $ equal to the measure on $A\times \left[
0,1\right] $ given by $\mu \left( da,dx\right) =k\left( a,x\right) \left\{
\nu \left( da\right) \otimes dx\right\} $, where%
\begin{equation*}
k\left( a,x\right) =\sum_{j=1}^{q}\mathbf{1}_{\left\{ j\right\} }\left(
a\right) m_{j}\left( x\right) \text{, \ \ }\left( a,x\right) \in A\times %
\left[ 0,1\right] \text{,}
\end{equation*}%
$\mathbf{1}_{\left\{ j\right\} }$ stands for the indicator of the singleton $%
\left\{ j\right\} $, and $dx$ is once again Lebesgue measure. Finally, we
define
\begin{eqnarray}
\mathfrak{H}^{\sharp } &=&L^{2}\left( A\times \left[ 0,1\right] ,\mu ,%
\mathcal{A}\otimes \mathcal{B}\left[ 0,1\right] \right)  \label{NHS} \\
\pi _{t}^{\sharp }h^{\sharp }\left( a,x\right) &=&h^{\sharp }\left(
a,x\right) \mathbf{1}_{\left[ 0,t\right] }\left( x\right) \text{,}  \notag
\end{eqnarray}%
for every $h^{\sharp }\in \mathfrak{H}^{\sharp }$, every $t\in \left[ 0,1%
\right] $, and every $\left( a,x\right) \in A\times \left[ 0,1\right] $. We
now introduce a transformation $T$ defined on a dense subset of $\mathfrak{H}%
^{\sharp }$: for every $h^{\sharp }\in \mathfrak{H}^{\sharp }$ with the form
\begin{equation}
h^{\sharp }\left( a,x\right) =\sum_{j=1}^{M}c_{j}\mathbf{1}_{\left\{
k_{j}\right\} }\left( a\right) \mathbf{1}_{\left[ 0,u_{j}\right] }\left(
x\right)  \label{basic}
\end{equation}%
($\left( a,x\right) \in A\times \left[ 0,1\right] $), where $M,k_{j}\geq 1$,
$u_{j}\in \left[ 0,1\right] $ and $c_{j}\in \mathbb{R}$ ($j=1,...,M$),%
\begin{equation*}
Th^{\sharp }=\sum_{j=1}^{M}c_{j}\pi _{u_{j}}g_{k_{j}}\text{,}
\end{equation*}%
where the $g_{k}$'s are the elements of the full orthogonal set $S$, as
defined in (\ref{FO}). By using the relation,
\begin{equation*}
\left( Th^{\sharp },Th^{\sharp }\right) _{\mathfrak{H}}=\left( h^{\sharp
},h^{\sharp }\right) _{\mathfrak{H}^{\sharp }}\text{, }
\end{equation*}%
which is verified for every $h^{\sharp }$ as in (\ref{basic}), one
immediately sees that $T$ can be extended by density to a unitary
transformation, from $\mathfrak{H}^{\sharp }$ onto $\mathfrak{H}$, and
moreover, since, for every $t\in \left[ 0,1\right] $,
\begin{equation*}
\pi _{t}Th^{\sharp }=\sum_{j=1}^{M}c_{j}\pi _{u_{j}\wedge t}g_{k_{j}}=T\pi
_{t}^{\sharp }h^{\sharp }\text{,}
\end{equation*}%
condition (\ref{EQUI}) is verified. We note again $T$ this extended
isomorphism, and, for $d\geq 2$, we write $T^{d}\triangleq T^{\otimes d}$,
and also $T^{1}\triangleq T^{\otimes 1}=T$. Observe that $T^{d}$ is an
isomorphism\ from $\left( \mathfrak{H}^{\sharp }\right) ^{\otimes d}$ onto $%
\mathfrak{H}^{\otimes d}$, $\left( T^{-1}\right) ^{d}=\left( T^{d}\right)
^{-1}$. Also, for $t\in \left[ 0,1\right] $ and due to (\ref{EQUI}),%
\begin{equation}
T^{d}\left( \pi _{t}^{\sharp }\right) ^{\otimes d}=\left( T\pi _{t}^{\sharp
}\right) ^{\otimes d}=\left( \pi _{t}T\right) ^{\otimes d}=\pi _{t}^{\otimes
d}T^{d}\text{.}  \label{EQUI2}
\end{equation}%
Now, for an absolutely continuous resolution $\pi $ on $\mathfrak{H}$, and
for $h^{\sharp }$ and $\mathfrak{H}^{\sharp }$ as in (\ref{NHS}), we define $%
X=\left\{ X\left( f\right) :f\in \mathfrak{H}\right\} $ to be an isonormal
Gaussian process over $\mathfrak{H}$, and set
\begin{equation}
X_{T}=X_{T}\left( \mathfrak{H}^{\sharp }\right) =\left\{ X_{T}\left(
h^{\sharp }\right) :h^{\sharp }\in \mathfrak{H}^{\sharp }\right\} \text{,}
\label{isodiesis}
\end{equation}%
where $X_{T}\left( h^{\sharp }\right) \triangleq X\left( Th^{\sharp }\right)
$, $\forall h^{\sharp }\in \mathfrak{H}^{\sharp }$. It is clear that, due to
the isometric property of $T$, $X_{T}$ is an isonormal Gaussian process over
$\mathfrak{H}^{\sharp }$. The proof of the following useful lemma is
deferred to the Appendix.

\begin{lemma}
\label{L : opIdentities}Under the above notation,

\begin{enumerate}
\item For every $d\geq 1$, $f\in \mathfrak{H}^{\otimes d}$, $I_{d}^{X}\left(
f\right) =I_{d}^{X_{T}}\left( \left( T^{d}\right) ^{-1}f\right) ;$

\item $\mathbb{D}_{X}^{1,2}=\mathbb{D}_{X_{T}}^{1,2}$, and, for every $F\in
\mathbb{D}_{X}^{1,2}$,
\begin{equation}
D_{X}F=T\left( D_{X_{T}}F\right) \text{;}  \label{Id.Diff.}
\end{equation}

\item For every $t\in \left[ 0,1\right] $,
\begin{equation*}
\mathcal{F}_{t}^{\pi }\left( X\right) =\sigma \left\{ X\left( \pi
_{t}f\right) :f\in \mathfrak{H}\right\} =\sigma \left\{ X_{T}\left( \pi
_{t}^{\sharp }h^{\sharp }\right) :h^{\sharp }\in \mathfrak{H}^{\sharp
}\right\} =\mathcal{F}_{t}^{\pi ^{\sharp }}\left( X_{T}\right) \text{;}
\end{equation*}

\item For every $u\in L^{2}\left( \mathfrak{H}^{\sharp },X_{T}\right) $, $%
u\in L_{\pi ^{\sharp }}^{2}\left( \mathfrak{H}^{\sharp },X_{T}\right) $ if,
and only if, $Tu\in L_{\pi }^{2}\left( \mathfrak{H},X\right) $;

\item For every $F\in \mathbb{D}_{X}^{1,2}$, a.s.-$\mathbb{P}$,
\begin{eqnarray}
proj\left\{ D_{X}F\mid L_{\pi }^{2}\left( \mathfrak{H},X\right) \right\}
&=&T\circ proj\left\{ D_{X_{T}}F\mid L_{\pi ^{\sharp }}^{2}\left( \mathfrak{H%
}^{\sharp },X_{T}\right) \right\}  \label{id.proj} \\
&=&T\circ proj\left\{ T^{-1}D_{X}F\mid L_{\pi ^{\sharp }}^{2}\left(
\mathfrak{H}^{\sharp },X_{T}\right) \right\} ;  \label{id.proj2}
\end{eqnarray}

\item For every $d\geq 2$, $f\in \mathfrak{H}^{\otimes d}$ (therefore, $f$
need not be a symmetric tensor), $r=1,...,d-1$ and $t\in \left[ 0,1\right] $%
,
\begin{eqnarray}
&&\left\Vert \left( \pi _{1}^{\otimes 2r}-\pi _{t}^{\otimes 2r}\right)
\left( f\otimes _{d-r}^{\pi ,t}f\right) \right\Vert _{\mathfrak{H}^{\otimes
2r}}  \label{Id.Norm} \\
&=&\left\Vert \left( \left( \pi _{1}^{\sharp }\right) ^{\otimes 2r}-\left(
\pi _{t}^{\sharp }\right) ^{\otimes 2r}\right) \left( \left( T^{d}\right)
^{-1}f\otimes _{d-r}^{\pi ^{\sharp },t}\left( T^{d}\right) ^{-1}f\right)
\right\Vert _{\left( \mathfrak{H}^{\sharp }\right) ^{\otimes 2r}}\text{.}
\notag
\end{eqnarray}
\end{enumerate}
\end{lemma}

\bigskip

Now adopt the assumptions and notation of Theorem \ref{T : Main}. If $\pi
_{n}$ is absolutely continuous on $\mathfrak{H}_{n}$, for every $n\geq 1$
there exists an isomorphism $T_{n}$, from $\mathfrak{H}_{n}$ onto some space
$\mathfrak{H}_{n}^{\sharp }$, endowed with a resolution $\pi _{n}^{\sharp }$
as in (\ref{NHS}) and such that properties (\ref{EQUI}) and (\ref{EQUI2})
(with $T_{n}$ substituting $T$) are verified. We also note $X_{T_{n}}\left(
h^{\sharp }\right) =X_{n}\left( T_{n}h^{\sharp }\right) $, for every $%
h^{\sharp }\in \mathfrak{H}_{n}^{\sharp }$. It follows from Lemma \ref{L :
opIdentities}-6 above that, if for every $r=1,...,d-1$, relation (\ref{asCV}%
) is verified, then%
\begin{equation*}
\left\Vert \left( \left( \pi _{n,1}^{\sharp }\right) ^{\otimes 2r}-\left(
\pi _{n,t}^{\sharp }\right) ^{\otimes 2r}\right) \left( \left(
T_{n}^{d}\right) ^{-1}f_{n}\otimes _{d-r}^{\pi _{n}^{\sharp },t}\left(
T_{n}^{d}\right) ^{-1}f_{n}\right) \right\Vert _{\left( \mathfrak{H}%
_{n}^{\sharp }\right) ^{\otimes 2r}}\underset{n\rightarrow +\infty }{%
\rightarrow }0\text{.}
\end{equation*}

Moreover, thanks to Points 1 and 3 of Lemma \ref{L : opIdentities},%
\begin{equation*}
\mathbb{E}\left[ I_{d}^{X_{T_{n}}}\left( \left( T_{n}^{d}\right)
^{-1}f_{n}\right) \mid \mathcal{F}_{t_{n}}^{\pi _{n}^{\sharp }}\left(
X_{T_{n}}\right) \right] =\mathbb{E}\left[ F_{n}\mid \mathcal{F}%
_{t_{n}}^{\pi _{n}}\left( X_{n}\right) \right] \overset{\mathbb{P}}{%
\rightarrow }0
\end{equation*}%
and
\begin{equation*}
\mathbb{E}\left[ I_{d}^{X_{T_{n}}}\left( \left( T_{n}^{d}\right)
^{-1}f_{n}\right) ^{2}\mid \mathcal{F}_{t_{n}}^{\pi _{n}^{\sharp }}\left(
X_{T_{n}}\right) \right] =\mathbb{E}\left[ F_{n}^{2}\mid \mathcal{F}%
_{t_{n}}^{\pi _{n}}\left( X_{n}\right) \right] \overset{\mathbb{P}}{%
\rightarrow }Y\in \mathcal{F}_{\ast },
\end{equation*}%
from which, by using the first part of the proof, we deduce that
\begin{equation}
\left\Vert proj\left\{ D_{X_{T_{n}}}F_{n}\mid L_{\pi _{n}^{\sharp
}}^{2}\left( \mathfrak{H}_{n}^{\sharp },X_{T_{n}}\right) \right\}
\right\Vert _{\mathfrak{H}_{n}^{\sharp }}^{2}\overset{\mathbb{P}}{\underset{%
n\rightarrow +\infty }{\rightarrow }}Y\text{.}  \label{CVCV}
\end{equation}%
The proof of Theorem \ref{T : Main} is now concluded by using (\ref{CVCV})
and Theorem \ref{T : Gen CV} since, due to Lemma \ref{L : opIdentities}-5
above and the fact that $T$ is an isomorphism,%
\begin{eqnarray*}
\left\Vert proj\left\{ D_{X_{T_{n}}}F_{n}\mid L_{\pi _{n}^{\sharp
}}^{2}\left( \mathfrak{H}_{n}^{\sharp },X_{T_{n}}\right) \right\}
\right\Vert _{\mathfrak{H}_{n}^{\sharp }}^{2} &=&\left\Vert T_{n}\circ
proj\left\{ D_{X_{T_{n}}}F_{n}\mid L_{\pi _{n}^{\sharp }}^{2}\left(
\mathfrak{H}_{n}^{\sharp },X_{T_{n}}\right) \right\} \right\Vert _{\mathfrak{%
H}_{n}}^{2} \\
&=&\left\Vert proj\left\{ D_{X}F_{n}\mid L_{\pi _{n}}^{2}\left( \mathfrak{H}%
_{n},X_{n}\right) \right\} \right\Vert _{\mathfrak{H}_{n}}^{2}.
\end{eqnarray*}

\bigskip

\textbf{Remark }(\textit{Concrete realizations of Wiener spaces}) \textbf{--
}For the sake of completeness, we establish some connections between the
unitary transformation $T:\mathfrak{H}^{\sharp }\mapsto \mathfrak{H}$ used
in the last part of the preceding proof (see (\ref{ISO})) and the concept of
\textit{concrete (filtered) Wiener space} introduced in \cite[Section 5]{UZ}%
. In particular, we point out that every \textquotedblleft
filtered\textquotedblright\ isonormal Gaussian process such as the pair $%
\left( X_{T}\left( \mathfrak{H}^{\sharp }\right) ,\pi ^{\sharp }\right) $
introduced in (\ref{NHS}) and (\ref{isodiesis}), is equivalent (in a sense
analogous to \cite[Definition 5.1]{UZ}) to a concrete Wiener space whose
dimension equals the rank of $\pi ^{\sharp }$. To do this, fix $q\in \left\{
1,2,...,+\infty \right\} $ and define $C_{0}\left( \left[ 0,1\right] \right)
$ to be set of continuous functions on $\left[ 0,1\right] $ that are
initialized at zero. We define $\mathbb{W}_{\left( q\right) }$ to be the set
of all $q$-dimensional vectors of the type $\mathbf{w}_{\left( q\right)
}=\left( w_{1},w_{2},...,w_{q}\right) $ (plainly, if $q=+\infty $, $\mathbf{w%
}_{\left( q\right) }$ is an infinite sequence) where $\forall i$, $w_{i}\in
C_{0}\left( \left[ 0,1\right] \right) $. The set $\mathbb{W}_{\left(
q\right) }$ is endowed with the norm $\left\Vert \mathbf{w}_{\left( q\right)
}\right\Vert _{\left( q\right) }=\sup_{i\leq q}\left\vert w_{i}\right\vert $%
, where $\left\vert w_{i}\right\vert =\sup_{t\in \left[ 0,1\right]
}\left\vert w_{i}\left( t\right) \right\vert $. Under $\left\Vert \mathbf{%
\cdot }\right\Vert _{\left( q\right) }$, $\mathbb{W}_{\left( q\right) }$ is
a Banach space. Now consider an Hilbert space $\mathfrak{H}$, as well as a
resolution $\pi \in \mathcal{R}_{AC}\left( \mathfrak{H}\right) $ such that $%
rank\left( \pi \right) =q$. We define $S=\left\{ g_{j}:1\leq j\leq q\right\}
$ to be the fully orthogonal $\pi $-reproducing subset of $\mathfrak{H}$
appearing in formula (\ref{FO}), and associate to each $g_{j}\in S$ the
function $m_{j}\in L^{1}\left( \left[ 0,1\right] ,dx\right) $ satisfying (%
\ref{Acont}), in such a way that (\ref{AcontConv}) is verified. To the pair $%
\left( \mathfrak{H},\pi \right) $ we associate the Hilbert space $\mathbb{H}%
_{\left( q\right) }$ and a resolution of the identity $\pi ^{\left( q\right)
}=\left\{ \pi _{s}^{\left( q\right) }:s\in \left[ 0,1\right] \right\} \in
\mathcal{R}\left( \mathbb{H}_{\left( q\right) }\right) $ defined as follows:
(i) $\mathbb{H}_{\left( q\right) }$ is the collection of all vectors of the
kind $\mathbf{h}_{\left( q\right) }=\left( h_{1},h_{2}...,h_{q}\right) $,
where, for each $j\leq q$, $h_{j}$ is a function of the form $h_{j}\left(
t\right) =\int_{0}^{t}h_{j}^{\prime }\left( x\right) dx$, for some $%
h_{j}^{\prime }\in L^{2}\left( \left[ 0,1\right] ,m_{j}\left( x\right)
dx\right) $, and also%
\begin{equation}
\sum_{j=1}^{q}\int_{0}^{1}\left( h_{j}^{\prime }\left( x\right) \right)
^{2}m_{j}\left( x\right) dx<+\infty \text{;}  \label{kvnorm}
\end{equation}%
(ii) $\mathbb{H}_{\left( q\right) }$ is endowed with the inner product
\begin{equation}
\left( \mathbf{h}_{\left( q\right) },\mathbf{k}_{\left( q\right) }\right)
_{\left( q\right) }=\sum_{j=1}^{q}\int_{0}^{1}h_{j}^{\prime }\left( x\right)
k_{j}^{\prime }\left( x\right) m_{j}\left( x\right) dx\text{,}  \label{opq}
\end{equation}%
whereas $\left\vert \cdot \right\vert _{\left( q\right) }=\left( \mathbf{%
\cdot },\mathbf{\cdot }\right) _{\left( q\right) }^{1/2}$ is the
corresponding norm; (iii) for every $s\in \left[ 0,1\right] $ and every $%
\mathbf{h}_{\left( q\right) }=\left( h_{1},h_{2}...,h_{q}\right) \in \mathbb{%
H}_{\left( q\right) }$,
\begin{equation}
\pi _{s}^{\left( q\right) }\mathbf{h}_{\left( q\right) }=\left(
h_{1}^{s},...,h_{q}^{s}\right) \text{, where \ }h_{j}^{s}\left( t\right)
\triangleq \int_{0}^{t\wedge s}h_{j}^{\prime }\left( x\right) dx\text{.}
\label{absPr}
\end{equation}%
Note that $\mathbb{H}_{\left( q\right) }\subset \mathbb{W}_{\left( q\right) }
$, and therefore $\mathbb{W}_{\left( q\right) }^{\ast }\subset \mathbb{H}%
_{\left( q\right) }^{\ast }=\mathbb{H}_{\left( q\right) }$. Moreover, from
relation (\ref{kvnorm}) it follows that the restriction of $\left\Vert
\mathbf{\cdot }\right\Vert _{\left( q\right) }$ to $\mathbb{H}_{\left(
q\right) }$ is a measurable seminorm, in the sense of \cite[Definition 4.4]%
{Kuo}. Also, $\mathbb{W}_{\left( q\right) }$ is the completion of $\mathbb{H}%
_{\left( q\right) }$ with respect to $\left\Vert \mathbf{\cdot }\right\Vert
_{\left( q\right) }$, and $\mathbb{W}_{\left( q\right) }^{\ast }\ $is dense
in $\mathbb{H}_{\left( q\right) }$ with respect to the norm $\left\vert
\cdot \right\vert _{\left( q\right) }$. As a consequence (see again \cite[%
Theorem 4.1]{Kuo}), there exists a canonical Gaussian measure $\mu _{\left(
q\right) }$ on $\left( \mathbb{W}_{\left( q\right) },\mathcal{B}\left(
\mathbb{W}_{\left( q\right) }\right) \right) $, such that, for every $\left(
\mathbf{l}_{1},...,\mathbf{l}_{m}\right) \in \left( \mathbb{W}_{\left(
q\right) }^{\ast }\right) ^{m}$, the mapping $\mathbf{w}_{\left( q\right)
}\mapsto \left( \mathbf{l}_{1}\left( \mathbf{w}_{\left( q\right) }\right)
,...,\mathbf{l}_{m}\left( \mathbf{w}_{\left( q\right) }\right) \right) $
defines a centered Gaussian vector such that, for every $j=1,...,m$ and
every $\lambda \in \mathbb{R}$,
\begin{equation}
\mathbb{E}_{\mu _{\left( q\right) }}\left[ \exp \left( i\lambda \mathbf{l}%
_{j}\right) \right] \triangleq \int_{\mathbb{W}_{\left( q\right) }^{\ast
}}\exp \left( i\lambda \mathbf{l}_{j}\left( \mathbf{w}_{\left( q\right)
}\right) \right) d\mu _{\left( q\right) }\left( \mathbf{w}_{\left( q\right)
}\right) =\exp \left( -\frac{\lambda ^{2}}{2}\left\vert \mathbf{l}%
_{j}\right\vert _{\left( q\right) }^{2}\right) \text{.}  \label{Gross}
\end{equation}%
Following \cite[p. 26]{UZ}, the triple $\left( \mathbb{W}_{\left( q\right) },%
\mathbb{H}_{\left( q\right) },\mu _{\left( q\right) }\right) $ (endowed with
the resolution $\pi ^{\left( q\right) }$ defined in (\ref{absPr})) is called
a \textit{concrete Wiener space }of dimension $q$. Note that, since $\mathbb{%
W}_{\left( q\right) }^{\ast }\ $is dense in $\mathbb{H}_{\left( q\right) }$,
there exists a unique collection of centered Gaussian random variables
defined on $\left( \mathbb{W}_{\left( q\right) },\mathcal{B}\left( \mathbb{W}%
_{\left( q\right) }\right) \right) $, denoted
\begin{equation}
X_{\left( q\right) }=X_{\left( q\right) }\left( \mathbb{H}_{\left( q\right)
}\right) =\left\{ X_{\left( q\right) }\left( \mathbf{h}_{\left( q\right)
}\right) :\mathbf{h}_{\left( q\right) }\in \mathbb{H}_{\left( q\right)
}\right\} \text{,}  \label{isoCWS}
\end{equation}%
such that $X_{\left( q\right) }\left( \mathbf{l}\right) \left( \mathbf{w}%
_{\left( q\right) }\right) =\mathbf{l}\left( \mathbf{w}_{\left( q\right)
}\right) $ for every $\mathbf{l}\in \mathbb{W}_{\left( q\right) }^{\ast }$
and $\mathbb{E}_{\mu _{\left( q\right) }}\left[ \exp \left( i\lambda
X_{\left( q\right) }\left( \mathbf{h}_{\left( q\right) }\right) \right) %
\right] $ $=$ $\exp \left( -\frac{\lambda ^{2}}{2}\left\vert \mathbf{h}%
_{\left( q\right) }\right\vert _{\left( q\right) }^{2}\right) $, $\forall
\mathbf{h}_{\left( q\right) }\in \mathbb{H}_{\left( q\right) }$. In
particular, $X_{\left( q\right) }\left( \mathbb{H}_{\left( q\right) }\right)
$ is an isonormal Gaussian process over $\mathbb{H}_{\left( q\right) }$. Now
consider the Hilbert space $\mathfrak{H}^{\sharp }$ and the resolution $\pi
^{\sharp }$ defined in (\ref{NHS}), and define the application $T_{\circ }:%
\mathfrak{H}^{\sharp }\mapsto \mathbb{H}_{\left( q\right) }$ as follows: for
every $h^{\sharp }\left( a,x\right) \in $ $\mathfrak{H}^{\sharp }$,
\begin{equation*}
T_{\circ }h^{\sharp }=\left( \int_{0}^{\cdot }h^{\sharp }\left( 1,x\right)
dx,...,\int_{0}^{\cdot }h^{\sharp }\left( q,x\right) dx\right) .
\end{equation*}%
It is easily seen that $T_{\circ }$ is a unitary transformation such that $%
T_{\circ }\pi _{t}^{\sharp }=\pi _{t}^{\left( q\right) }T_{\circ }$ for
every $t$, thus implying that the two pairs $\left( \mathfrak{H}^{\sharp
},\pi ^{\sharp }\right) $ and $\left( \mathbb{H}_{\left( q\right) },\pi
^{\left( q\right) }\right) $, and hence the two filtered isonormal processes
$\left( X_{\left( q\right) }\left( \mathbb{H}_{\left( q\right) }\right) ,\pi
^{\left( q\right) }\right) $ and $\left( X_{T}\left( \mathfrak{H}^{\sharp
}\right) ,\pi ^{\sharp }\right) $, are equivalent in the sense of \cite[%
Definition 5.1]{UZ}.

\subsection{Proof of Theorem \protect\ref{T : Main Double}}

The implications (i) $\Longrightarrow $ (ii) and (iv) $\Longrightarrow $ (i)
(in which assumption (\ref{L2/PRO}) is immaterial) are consequences,
respectively, of Theorem \ref{T : Gen CV} and Theorem \ref{T : Main}. Now
suppose (ii) is verified. Since $\mathbb{E}\left[ F_{n}^{2}\mid \mathcal{F}%
_{t_{n}}^{\pi _{n}}\left( X_{n}\right) \right] \overset{\mathbb{P}}{%
\rightarrow }Y$ by assumption, we may use the second part of Theorem \ref{T
: GenConv 2} to deduce that for every sequence $n\left( k\right) $, there
exists a subsequence $n\left( k_{r}\right) $, $r\geq 1$, s.t., a.s.-$\mathbb{%
P}$,%
\begin{equation*}
\mathbb{E}\left[ \exp \left( i\lambda F_{n\left( k_{r}\right) }\right) \mid
\mathcal{F}_{t_{n\left( k_{r}\right) }}^{\pi _{n\left( k_{r}\right) }}\left(
X_{n\left( k_{r}\right) }\right) \right] \underset{r\rightarrow +\infty }{%
\rightarrow }\exp \left( -\frac{\lambda ^{2}}{2}Y\right) ,\ \ \forall
\lambda \in
\mathbb{R}
\text{.}
\end{equation*}

Moreover, since the usual properties of multiple Wiener-It\^{o} integrals
(see e.g. \cite[Chapter VI]{Jansson}) imply that, a.s.-$\mathbb{P}$ and due
to (\ref{NORM}),%
\begin{equation*}
\sup_{r\geq 1}\mathbb{E}\left[ \left\vert F_{n\left( k_{r}\right)
}\right\vert ^{M}\mid \mathcal{F}_{t_{n\left( k_{r}\right) }}^{\pi _{n\left(
k_{r}\right) }}\left( X_{n\left( k_{r}\right) }\right) \right] <+\infty
\text{, \ \ }\forall M\geq 1\text{,}
\end{equation*}%
we conclude that, a.s.-$\mathbb{P}$,%
\begin{equation*}
\mathbb{E}\left[ \left( F_{n\left( k_{r}\right) }\right) ^{4}\mid \mathcal{F}%
_{t_{n\left( k_{r}\right) }}^{\pi _{n\left( k_{r}\right) }}\left( X_{n\left(
k_{r}\right) }\right) \right] \underset{r\rightarrow +\infty }{\rightarrow }%
3Y^{2},
\end{equation*}%
and therefore that (iii) holds. To conclude, assume that the two conditions
(iii) and (\ref{L2/PRO}) are verified, and write%
\begin{equation}
F_{n}=I_{2}^{X_{n}}\left( \pi _{n,t_{n}}^{\otimes 2}f_{n}\right)
+2I_{2}^{X_{n}}\left( \pi _{n,t_{n}}\otimes \left( \pi _{n,1}-\pi
_{n,t_{n}}\right) f_{n}\right) +I_{2}^{X_{n}}\left( \left( \pi _{n,1}-\pi
_{n,t_{n}}\right) ^{\otimes 2}f_{n}\right) \triangleq
F_{n,0}+F_{n,1}+F_{n,2}.  \label{dec}
\end{equation}

Due to (\ref{NEG}), $F_{n,0}\rightarrow 0$ in $L^{2}$. Also, for every $%
n\geq 1$, $F_{n,2}$ is independent of $\mathcal{F}_{t_{n}}^{\pi _{n}}\left(
X_{n}\right) $ and, conditionally on $\mathcal{F}_{t_{n}}^{\pi _{n}}\left(
X_{n}\right) $, $F_{n,1}$ is a centered Gaussian random variable. Moreover%
\begin{equation*}
\mathbb{E}\left[ \left( F_{n}\right) ^{2}\mid \mathcal{F}_{t_{n}}^{\pi
_{n}}\left( X_{n}\right) \right] =\mathbb{E}\left[ \left( F_{n,0}\right)
^{2}+\left( F_{n,1}\right) ^{2}+\left( F_{n,2}\right) ^{2}\mid \mathcal{F}%
_{t_{n}}^{\pi _{n}}\left( X_{n}\right) \right] .
\end{equation*}%
By writing $A_{n}\sim B_{n}$ to indicate that $A_{n}-B_{n}\overset{\mathbb{P}%
}{\rightarrow }0$, we have therefore%
\begin{eqnarray*}
\mathbb{E}\left[ \left( F_{n}\right) ^{4}\mid \mathcal{F}_{t_{n}}^{\pi
_{n}}\left( X_{n}\right) \right] &\sim &\mathbb{E}\left[ \left(
F_{n,1}\right) ^{4}\mid \mathcal{F}_{t_{n}}^{\pi _{n}}\left( X_{n}\right) %
\right] +\mathbb{E}\left[ \left( F_{n,2}\right) ^{4}\mid \mathcal{F}%
_{t_{n}}^{\pi _{n}}\left( X_{n}\right) \right] \\
&&+6\mathbb{E}\left[ \left( F_{n,1}F_{n,2}\right) ^{2}\mid \mathcal{F}%
_{t_{n}}^{\pi _{n}}\left( X_{n}\right) \right] \\
&=&3\mathbb{E}\left[ \left( F_{n,1}\right) ^{2}\mid \mathcal{F}_{t_{n}}^{\pi
_{n}}\left( X_{n}\right) \right] ^{2}+\mathbb{E}\left[ \left( F_{n,2}\right)
^{4}\right] \\
&&+6\mathbb{E}\left[ \left( F_{n,1}F_{n,2}\right) ^{2}\mid \mathcal{F}%
_{t_{n}}^{\pi _{n}}\left( X_{n}\right) \right] .
\end{eqnarray*}

By reasoning as in \cite[pp. 182-183]{NuPe}, and noting $f_{n,0}=\left( \pi
_{n,1}-\pi _{n,t_{n}}\right) ^{\otimes 2}f_{n}\in \mathfrak{H}^{\odot 2}$,%
\begin{eqnarray*}
\mathbb{E}\left[ \left( F_{n,2}\right) ^{4}\right] &=&3\left\Vert
f_{n,0}\right\Vert _{\mathfrak{H}^{\odot 2}}^{4}+48\left\Vert f_{n,0}\otimes
_{1}f_{n,0}\right\Vert _{\mathfrak{H}^{\otimes 2}}^{2} \\
&=&3\mathbb{E}\left[ \left( F_{n,2}\right) ^{2}\mid \mathcal{F}_{t_{n}}^{\pi
_{n}}\left( X_{n}\right) \right] ^{2}+48\left\Vert \left( \pi _{n,1}-\pi
_{n,t_{n}}\right) ^{\otimes 2}f_{n}\otimes _{1}^{\pi
_{n},t_{n}}f_{n}\right\Vert _{\mathfrak{H}^{\otimes 2}}^{2}.
\end{eqnarray*}

Standard calculations yield finally that, since (\ref{NEG}) and (\ref{NORM})
hold, there exist constants $c_{1},c_{2}>0$ such that
\begin{eqnarray*}
\mathbb{E}\left[ \left( \mathbb{E}\left[ \left( F_{n}\right) ^{4}\mid
\mathcal{F}_{t_{n}}^{\pi _{n}}\left( X_{n}\right) \right] -3Y^{2}\right) ^{2}%
\right]  &=&c_{1}\left\Vert \left( \pi _{n,1}-\pi _{n,t_{n}}\right)
^{\otimes 2}f_{n}\otimes _{1}^{\pi _{n},t_{n}}f_{n}\right\Vert _{\mathfrak{H}%
^{\otimes 2}}^{2} \\
&&+c_{2}\left\Vert \left( \pi _{n,1}-\pi _{n,t_{n}}\right) \otimes \pi
_{n,t_{n}}\left( f_{n}\otimes _{1}^{\pi _{n},t_{n}}f_{n}\right) \right\Vert
_{\mathfrak{H}^{\otimes 2}}^{2}\text{,}
\end{eqnarray*}%
and, since (\ref{L2/PRO}) is verified and%
\begin{equation*}
\left( \pi _{n,1}-\pi _{n,t_{n}}\right) ^{\otimes 2}+\left( \pi _{n,1}-\pi
_{n,t_{n}}\right) \otimes \pi _{n,t_{n}}+\pi _{n,t_{n}}\otimes \left( \pi
_{n,1}-\pi _{n,t_{n}}\right) =\pi _{n,1}^{\otimes 2}-\pi _{n,t_{n}}^{\otimes
2}\text{,}
\end{equation*}%
we obtain immediately the desired implication (iii)$\Longrightarrow $(iv).$%
\blacksquare $

\section{Appendix}

\textbf{Proof of Lemma \ref{L : opIdentities}} -- (\textit{Point 1}) Let $%
\left\{ e_{j}:j\geq 1\right\} $ be an orthonormal basis of $\mathfrak{H}$,
and define, for $d\geq 1$, $\mathbb{A}\left[ d\right] $ to be the set of
sequences $\left( a_{1},a_{2},...\right) $ with values in $\mathbb{N}$, and
such that $\sum_{j\geq 1}a_{j}=d$ (note that this implies that there are
only finitely many $a_{j}$ that are different from zero). Then, a total set
in $\mathfrak{H}^{\odot d}$ is given by%
\begin{equation*}
\mathcal{A}_{d}=\left\{ (\bigotimes\limits_{j=1}^{\infty }e_{j}^{\otimes
a_{j}})_{s}:\left( a_{1},a_{2},...\right) \in \mathbb{A}\left[ d\right]
\right\} ,
\end{equation*}%
where $e^{\otimes 0}=1$ by definition, and $\left( \cdot \right) _{s}$
indicates symmetrization. Moreover, a classic characterization of multiple
stochastic integrals (see \cite[Ch. 1]{Nualart}) as well as the fact that $%
X\left( e_{i}\right) =X\left( TT^{-1}e_{i}\right) =X_{T}\left(
T^{-1}e_{i}\right) $ by definition, imply the following relations: for every
$\left( a_{1},a_{2},...\right) \in \mathbb{A}\left[ d\right] $,%
\begin{eqnarray*}
I_{d}^{X}\left( (\bigotimes\limits_{j=1}^{\infty }e_{j}^{\otimes
a_{j}})_{s}\right) &=&d!\prod_{j=1}^{\infty }H_{a_{j}}\left( X\left(
e_{j}\right) \right) =d!\prod_{j=1}^{\infty }H_{a_{j}}\left( X_{T}\left(
T^{-1}e_{j}\right) \right) \\
&=&I_{d}^{X_{T}}\left( (\bigotimes\limits_{j=1}^{\infty }\left(
T^{-1}e_{j}\right) ^{\otimes a_{j}})_{s}\right) =I_{d}^{X_{T}}\left( \left(
T^{d}\right) ^{-1}(\bigotimes\limits_{j=1}^{\infty }\left( e_{j}\right)
^{\otimes a_{j}})_{s}\right) ,
\end{eqnarray*}%
where $\left\{ H_{a}:a\geq 1\right\} $ is the family of Hermite polynomials
defined e.g. in \cite[p. 4]{Nualart}. It is therefore clear that $%
I_{d}^{X}\left( f\right) =I_{d}^{X_{T}}\left( \left( T^{d}\right)
^{-1}f\right) $ is true for every $f$ that is a linear combination of
elements of $\mathcal{A}_{d}$, and the general result is achieved by a
standard density argument. (\textit{Point 2}) For $m\geq 1$, let $%
C_{b}^{\infty }\left( \mathbb{R}^{m}\right) $ denote the class of bounded
and infinitely differentiable functions on $\mathbb{R}^{m}$, whose
derivatives are also bounded. We start by observing that, since $T$ is a
one-to-one unitary transformation, random variables of the type
\begin{equation}
F=f\left( X\left( h_{1}\right) ,...,X\left( h_{m}\right) \right) =f\left(
X_{T}\left( T^{-1}h_{1}\right) ,...,X_{T}\left( T^{-1}h_{m}\right) \right)
\label{smooth}
\end{equation}%
(the equality is again a consequence of the relation $X\left( h_{i}\right)
=X\left( TT^{-1}e_{i}\right) =X_{T}\left( T^{-1}h_{i}\right) $), where $%
m\geq 1$, $f\in C_{b}^{\infty }\left( \mathbb{R}^{m}\right) $ and $%
h_{1},...,h_{m}\in C_{b}^{\infty }\left( \mathbb{R}^{m}\right) $, are dense
both in $\mathbb{D}_{X}^{1,2}$ and $\mathbb{D}_{X_{T}}^{1,2}$. To conclude,
use a density argument, as well as the fact that, for $F$ as in (\ref{smooth}%
),%
\begin{eqnarray*}
D_{X}F &=&\sum_{j=1}^{m}\frac{\partial }{\partial x_{j}}f\left( X\left(
h_{1}\right) ,...,X\left( h_{m}\right) \right) h_{j} \\
&=&\sum_{j=1}^{m}\frac{\partial }{\partial x_{j}}f\left( X_{T}\left(
T^{-1}h_{1}\right) ,...,X_{T}\left( T^{-1}h_{m}\right) \right)
TT^{-1}h_{j}=TD_{X_{T}}F,
\end{eqnarray*}%
hence proving (\ref{Id.Diff.}). (\textit{Point 3}) This is a consequence of
the relations $X\left( \pi _{t}h\right) =X_{T}\left( T^{-1}\pi _{t}h\right) $
$=$ $X_{T}\left( T^{-1}\pi _{t}TT^{-1}h\right) =X_{T}\left( \pi _{t}^{\sharp
}T^{-1}h\right) $, that are verified for every $t\in \left[ 0,1\right] $,
since $T^{-1}\pi _{t}T=\pi _{t}^{\sharp }$, due to (\ref{EQUI}). (\textit{%
Point 4}) Suppose $u\in L_{\pi ^{\sharp }}^{2}\left( \mathfrak{H}^{\sharp
},X_{T}\right) $. Then, since $T$ is an isometry, $\mathbb{E}\left[
\left\Vert Tu\right\Vert _{\mathfrak{H}}^{2}\right] =\mathbb{E}\left[
\left\Vert u\right\Vert _{\mathfrak{H}^{\sharp }}^{2}\right] <+\infty $, and
therefore $Tu\in L_{\pi }^{2}\left( \mathfrak{H},X\right) $. To prove that $%
Tu$ is also $\pi $-adapted, use the fact that, since $T$ is an isometry and (%
\ref{EQUI}) holds, for every $t\in \left[ 0,1\right] $ and every $h\in
\mathfrak{H}$,
\begin{eqnarray*}
\left( Tu,\pi _{t}h\right) _{\mathfrak{H}} &=&\left( Tu,TT^{-1}\pi
_{t}h\right) _{\mathfrak{H}}=\left( u,T^{-1}\pi _{t}h\right) _{\mathfrak{H}%
^{\sharp }} \\
&=&\left( u,\pi _{t}^{\sharp }T^{-1}h\right) _{\mathfrak{H}^{\sharp }}\in
\mathcal{F}_{t}^{\pi ^{\sharp }}\left( X_{T}\right) =\mathcal{F}_{t}^{\pi
}\left( X\right) \text{,}
\end{eqnarray*}%
due to Point 3, thus yielding $u\in L_{\pi }^{2}\left( \mathfrak{H},X\right)
$. The opposite implication is obtained analogously. (\textit{Point 5})
Consider first an elementary random variable $\eta ^{\sharp }\in \mathcal{E}%
_{\pi ^{\sharp }}\left( \mathfrak{H}^{\sharp },X_{T}\right) $ with the form $%
\eta ^{\sharp }=\Phi \left( t\right) \left( \pi _{t+s}^{\sharp }-\pi
_{t}^{\sharp }\right) h^{\sharp }$, where $\Phi \left( t\right) \in \mathcal{%
F}_{t}^{\pi ^{\sharp }}\left( X_{T}\right) $ ($=\mathcal{F}_{t}^{\pi }\left(
X\right) $), $h^{\sharp }\in \mathfrak{H}^{\sharp }$ and $s,t\geq 0$. Then,
due to (\ref{EQUI}), $T\eta ^{\sharp }$ $=$ $\Phi \left( t\right) T\left(
\pi _{t+s}^{\sharp }-\pi _{t}^{\sharp }\right) h^{\sharp }$ $=$ $\Phi \left(
t\right) \left( \pi _{t+s}-\pi _{t}\right) Th^{\sharp }$, and therefore $%
Th^{\sharp }\in \mathcal{E}_{\pi }\left( \mathfrak{H},X\right) $. Now, for $%
F\in \mathbb{D}_{X_{T}}^{1,2}$, observe that a variable $P\in L_{\pi
^{\sharp }}^{2}\left( \mathfrak{H}^{\sharp },X_{T}\right) $ is equal to $%
proj\left\{ D_{X_{T}}F\mid L_{\pi ^{\sharp }}^{2}\left( \mathfrak{H}^{\sharp
},X_{T}\right) \right\} $ if, and only if, for every $\eta ^{\sharp }\in
\mathcal{E}_{\pi ^{\sharp }}\left( \mathfrak{H}^{\sharp },X_{T}\right) $ as
before%
\begin{equation}
\mathbb{E}\left[ \left( P,\eta ^{\sharp }\right) _{\mathfrak{H}^{\sharp }}%
\right] =\mathbb{E}\left[ \left( D_{X_{T}}F,\eta ^{\sharp }\right) _{%
\mathfrak{H}^{\sharp }}\right] .  \label{jo}
\end{equation}%
But, since $T$ is an isometry, (\ref{jo}) and (\ref{Id.Diff.}) imply also
that
\begin{equation*}
\mathbb{E}\left[ \left( TP,T\eta ^{\sharp }\right) _{\mathfrak{H}}\right] =%
\mathbb{E}\left[ \left( TD_{X_{T}}F,T\eta ^{\sharp }\right) _{\mathfrak{H}}%
\right] =\mathbb{E}\left[ \left( D_{X}F,T\eta ^{\sharp }\right) _{\mathfrak{H%
}}\right] .
\end{equation*}%
Hence, since $TP\in L_{\pi }^{2}\left( \mathfrak{H},X\right) $ due to Point
4,
\begin{equation*}
TP=T\circ proj\left\{ D_{X_{T}}F\mid L_{\pi ^{\sharp }}^{2}\left( \mathfrak{H%
}^{\sharp },X_{T}\right) \right\} =proj\left\{ D_{X}F\mid L_{\pi ^{\sharp
}}^{2}\left( \mathfrak{H}^{\sharp },X_{T}\right) \right\} ,
\end{equation*}%
thus proving (\ref{id.proj}). To prove (\ref{id.proj2}), just observe that (%
\ref{Id.Diff.}) implies that $D_{X_{T}}F=T^{-1}D_{X}F$. (\textit{Point 6})
Let again $\left\{ e_{j}:j\geq 1\right\} $ be an ONB of $\mathfrak{H}$. Note
first that, for every $d\geq 2$, $f\in \mathfrak{H}^{\otimes d}$, $%
r=1,...,d-1$, $t\in \left[ 0,1\right] $, and $i_{1},...,i_{d-r}\geq 1$%
\begin{eqnarray*}
&&\left( \left( \pi _{1}^{\otimes d-r}-\pi _{t}^{\otimes d-r}\right)
f,e_{i_{1}}\otimes \cdot \cdot \cdot \otimes e_{i_{d-r}}\right) _{\mathfrak{H%
}^{\otimes d-r}} \\
&=&\left( \left( T^{d-r}\right) ^{-1}\left( \pi _{1}^{\otimes d-r}-\pi
_{t}^{\otimes d-r}\right) f,T^{-1}e_{i_{1}}\otimes \cdot \cdot \cdot \otimes
T^{-1}e_{i_{d-r}}\right) _{\left( \mathfrak{H}^{\sharp }\right) ^{\otimes
d-r}} \\
&=&\left( \left( \left( \pi _{1}^{\sharp }\right) ^{\otimes d-r}-\left( \pi
_{t}^{\sharp }\right) ^{\otimes d-r}\right) \left( T^{d-r}\right)
^{-1}f,T^{-1}e_{i_{1}}\otimes \cdot \cdot \cdot \otimes
T^{-1}e_{i_{d-r}}\right) _{\left( \mathfrak{H}^{\sharp }\right) ^{\otimes
d-r}}.
\end{eqnarray*}%
Thanks to (\ref{op2}), it follows that%
\begin{eqnarray*}
&&\left( T^{2r}\right) ^{-1}f\otimes _{d-r}^{\pi ,t}f \\
&=&\sum_{i_{1},\ldots ,i_{d-r}=1}^{\infty }\ \left\langle \left( \left(
T^{d}\right) ^{-1}f\right) ,\left( \left( \pi _{1}^{\sharp }\right)
^{\otimes d-r}-\left( \pi _{t}^{\sharp }\right) ^{\otimes d-r}\right)
T^{-1}e_{i_{1}}\otimes \cdot \cdot \cdot \otimes
T^{-1}e_{i_{d-r}}\right\rangle _{\left( \mathfrak{H}^{\sharp }\right)
^{\otimes d-r}} \\
&&\text{ \ \ \ \ \ \ \ \ \ \ \ \ \ \ \ \ \ \ \ \ \ \ \ \ \ \ }\left\langle
\left( \left( T^{d}\right) ^{-1}f\right) ,\left( \left( \pi _{1}^{\sharp
}\right) ^{\otimes d-r}-\left( \pi _{t}^{\sharp }\right) ^{\otimes
d-r}\right) T^{-1}e_{i_{1}}\otimes \cdot \cdot \cdot \otimes
T^{-1}e_{i_{d-r}}\right\rangle _{\left( \mathfrak{H}^{\sharp }\right)
^{\otimes d-r}} \\
&=&\left( \left( T^{d}\right) ^{-1}f\otimes _{d-r}^{\pi ^{\sharp },t}\left(
T^{d}\right) ^{-1}f\right)
\end{eqnarray*}

As a consequence, by using (\ref{EQUI2}) and the fact that $T^{d}$ and $%
\left( T^{d}\right) ^{-1}$ are isometries,
\begin{eqnarray*}
&&\left\Vert \left( \left( \pi _{1}^{\sharp }\right) ^{\otimes 2r}-\left(
\pi _{t}^{\sharp }\right) ^{\otimes 2r}\right) \left( \left( T^{d}\right)
^{-1}f\otimes _{d-r}^{\pi ^{\sharp },t}\left( T^{d}\right) ^{-1}f\right)
\right\Vert _{\left( \mathfrak{H}^{\sharp }\right) ^{\otimes 2r}} \\
&=&\left\Vert \left( \left( \pi _{1}^{\sharp }\right) ^{\otimes 2r}-\left(
\pi _{t}^{\sharp }\right) ^{\otimes 2r}\right) \left( T^{2r}\right)
^{-1}f\otimes _{d-r}^{\pi ,t}f\right\Vert _{\left( \mathfrak{H}^{\sharp
}\right) ^{\otimes 2r}} \\
&=&\left\Vert \left( T^{2r}\right) ^{-1}\left( \pi _{1}^{\otimes 2r}-\pi
_{t}^{\otimes 2r}\right) f\otimes _{d-r}^{\pi ,t}f\right\Vert _{\left(
\mathfrak{H}^{\sharp }\right) ^{\otimes 2r}} \\
&=&\left\Vert \left( \pi _{1}^{\otimes 2r}-\pi _{t}^{\otimes 2r}\right)
\left( f\otimes _{d-r}^{\pi ,t}f\right) \right\Vert _{\mathfrak{H}^{\otimes
2r}}\text{,}
\end{eqnarray*}%
which proves (\ref{Id.Norm}). \ $\blacksquare $

\end{document}